\documentclass[11pt]{article}
\usepackage{amsfonts, latexsym, amssymb,amsmath,amstext}
\usepackage[dvips]{graphics,color}
\usepackage[mathscr]{eucal}
\newenvironment{proof}{\par\noindent\emph{Proof.}\ }{\hfill$\Box$\par}
\newenvironment{derivation}{\par\noindent\emph{Derivation.}\ }{\hfill$\Box$\par}
\newtheorem{theorem}{Theorem}[section]

\newtheorem{remark}{Remark}[section]
\newtheorem{corollary}{Corollary}[section]

\newtheorem{proposition}{Proposition}[section]
\makeatletter
\@addtoreset{equation}{section}

\makeatother
\begin{document}
\title{An Isomonodromy Cluster of Two Regular Singularities}
\author{A.~V.~Kitaev
\thanks{E-mail: kitaev@pdmi.ras.ru}\\
Steklov Mathematical Institute, Fontanka 27, St.Petersburg, 191023, Russia\\
and\\
School of Mathematics and Statistics, University of Sydney,\\
Sydney, NSW 2006, Australia}
\date{\today}
\maketitle
\begin{abstract}
We consider a linear $2\times2$ matrix ODE with two coalescing regular singularities.
This coalescence is restricted with an isomonodromy condition with respect to the
distance between the merging singularities in a way consistent with the ODE.
In particular, a zero-distance limit for the ODE exists. The monodromy group of the limiting ODE is
calculated in terms of the original one. This coalescing process generates a limit for
the corresponding nonlinear systems of isomonodromy deformations. In our main example
the latter limit reads as $P_6\to P_5$, where $P_n$ is the $n$-th Painlev\'e equation.
We also discuss some general problems which arise while studying the above-mentioned
limits for the Painlev\'e equations.
\vspace{24pt}\vspace{24pt}\\
{\bf Short title}: {Isomonodromy Cluster}
\end{abstract}
\newpage
\setcounter{page}2
\section{Prelude}
  \label{sec:prelude}
This work appeared first as Sfb 288 Preprint No.~149 (Teschnishe Universit\"at, Berlin)
in December 1994. It was distributed to all leading mathematical institutions and many
researchers. It was submitted to Comm. Math. Phys. (November 1994), but was not accepted
as "not interesting for the readership" of that Journal.
The copy of this preprint in various formats can be downloaded from
CERN Document Server\footnote[1]{http://preprints.cern.ch/cgi-bin/
setlink?base=preprint\&categ=.\&id=SCAN-9501196}. During the past (more than) 11 years
I presented this work in a number of talks on various conferences, colloquium talks,
and many private discussions. During this period appeared also a number of works by
different authors that discuss questions closely related with this work. A review
of those works actually requires considerable space and time. So, I was forced to
avoid this issue. The only changes I made in the bibliography is the inclusion of the
brief announcement~\cite{IKF} of the works~\cite{FIK1, FIK2} and also an appropriate
renumbering of the references.

In spite of the development that have been achieved during the past decade,
the main results obtained in my preprint were not reproduced. Moreover, the
intensive development of the Random Matrix Theory and related topics in the
theory of Orthogonal Polynomials reveals that various double scaling limits
of the Painlev\'e equations, one of those we study in this work, play a significant
role in various questions of these theories and applications.

So, I believe that the publication of this work will
be helpful as presenting an important property of the nonlinear special function -
the sixth Painlev\'e equation ($P_6$) and for future studies of isomonodromy deformations.

I decided not to make any changes into the mathematical content of this work, just
correction of English and misprints in formulae. Since it past so much time from the date
the work was written, it is important to mention that the discretization procedure of the
transition limit that studied in this work was motivated by the author's joint work with
Alexander Its and Athanasis Fokas on matrix models and orthogonal polynomials\cite{FIK1, FIK2}.

I am very grateful to the Guest Editors of this volume, Nalini Joshi and Frank
Nijhoff, for giving me an opportunity to publish this work. I am also
indebt to Fedor Andreev who typed the original version of this
manuscript\footnotemark[1] in \AmS-\TeX\; and Arthur Vartanian for
helping me in correction of English.
\section{Introduction}
 \label{sec:introduction}
By the isomonodromy cluster of two regular singularities we mean two points
in the complex plane ${\mathbb C}$ with respective positions $\lambda_i(t)$
$(i=0,1)$, considered as given functions of the parameter $t$ defined in
a neighborhood ${\cal O}(t_0)\subset\mathbb C$ (or $\mathbb R$) of some fixed
point $t_0\in\mathbb C$ (or $\mathbb R$), satisfying the following conditions:
(1) $\lambda_0(t_0)=\lambda_1(t_0)$; (2) $\lambda_0(t)\neq\lambda_1(t)$ for
$t\in{\cal O}(t_0)\setminus t_0$; and (3) they are the simple poles of the
linear ODE
\begin{equation}
 \label{eq:cluster-definition}
\frac{d}{d\lambda}\Psi=\left(\frac{A_0(t)}{\lambda-\lambda_0(t)}+
\frac{A_1(t)}{\lambda-\lambda_1(t)}+ A(\lambda,t)\right)\Psi,
\end{equation}
where $A_0(t), A_1(t), A(\lambda, t)\in\text{Mat}(n,\mathbb C)$ are analytic
functions of $t$, and $A(\lambda, t)$ is a rational function of $\lambda$.
It is also assumed that there exists a fundamental solution of
Eq.~(\ref{eq:cluster-definition}) with the manifold of monodromy
data\footnote[2]{The manifold of monodromy data is defined in Section~\ref{sec:2}}
independent of $t$ ({\it hard isomonodromy condition}).\\
Thus Eq.~(\ref{eq:cluster-definition}) is equipped with a system of nonlinear
ODEs (with respect to $t$) governing isomonodromy deformations of its coefficients.

For the simplest nontrivial backgrounds (the functions $A(\lambda,t)$), the systems of
isomonodromy deformations can be reduced to classical Painlev\'e equations; in particular,
in this work we consider the fifth Painlev\'e equation $P_5$
($n=2$, $\lambda-0=0$, $\lambda_1=t$, $A(\lambda,t)={\rm const}\cdot\sigma_3\neq0$)
and the sixth one $P_6$ ($n=2$, $\lambda_0=0$, $\lambda_1=t$, $A(\lambda, t)=A(t)/(\lambda-1)$,
$A_0(t)+A_1(t)+A(t)={\rm const}\cdot\sigma_3\neq0$).

The main idea in studying the cluster system~(\ref{eq:cluster-definition}), equipped with the
isomonodromy condition, is to substitute into Eq.~(\ref{eq:cluster-definition}), instead of the
cluster entries, some singularity of a regular or irregular type. Carrying out this procedure in
accordance with the isomonodromy condition, we obtain, instead of (\ref{eq:cluster-definition}),
a new equation, whose fundamental solution solution is denoted by $\Psi_{\rm new}$, and a novel
system of isomonodromy deformations with respect to a parameter which is proportional to $t-t_0$.
As a result we obtain a formal limit (or formal asymptotics) as $t\to t_0$ of the initial
Eq.~(\ref{eq:cluster-definition}) to a new one, and of the initial system of isomonodromy
deformations to a novel one. To translate this formal limit into the nonformal result, i.e., to
understand the solutions of the new system of isomonodromy deformations as defining the master term
of asymptotic expansions as $t\to t_0$ for the solutions of the initial system of isomonodromy
deformations, one has to solve the following nontrivial problem: to calculate the manifold of
monodromy data for $\Psi_{\rm new}$ in terms of the monodromy manifold for $\Psi$. This is the
problem addressed here by taking the background in Eq.~(\ref{eq:cluster-definition}) corresponding
to $P_6$. The latter result is easy to generalize to an arbitrary number of regular singularities.

It was M. Jimbo~\cite{J} who considered isomonodromy cluster of two regular singularities for the
linear ODEs associated with $P_5$ and $P_6$. In accordance with the aforementioned, Jimbo's work can be
interpreted as follows: Jimbo considered Eq.~(\ref{eq:cluster-definition}) with $P_5$ and $P_6$
backgrounds and substituted for the cluster terms a regular singularity. The resulting ODE for
$\Psi_{\rm new}$, by means of simple gauge and scaling transformations, can be reduced to matrix
versions of the standard ODEs for the Gauss hypergeometric function in the case
of $P_6$, and the Whittaker confluent hypergeometric function in the case of $P_5$. From this, one
observes:
\begin{enumerate}
\item
The coefficients of the standard ODEs, which are constant (with respect to $t$) parameters,
define the manifolds of monodromy data for the Gauss hypergeometric and Whittaker functions
and, thereby, define the manifold of monodromy data for $\Psi_{\rm new}$;
\item
The coefficients of the standard ODEs also define coefficients of the ODE for $\Psi_{\rm new}$,
the latter coefficients, which define solutions of the new system of isomonodromy deformations,
asymptotically, as $t\to0$ behave as linear combinations of power-like functions, $t^\alpha$.
\end{enumerate}
Because the deformations are isomonodromic, it is easy to relate the monodromy manifolds of
the initial function $\Psi$ with the monodromy manifolds  of the standard equations.
Thus, having done items 1 and 2 for each Painlev\'e equation, $P_5$ and $P_6$,
Jimbo was able to parameterize the small-$t$ asymptotics for $P_5$ and $P_6$ by points of
the corresponding manifolds of monodromy data; for $P_6$, i.e., when Eq.~(\ref{eq:cluster-definition})
has, along with the cluster $\{0,t\}$, two additional singularities at $1$ and $\infty$, Jimbo
considered the group of fractional-linear transformations of $\mathbb C$ which interchange the cluster with the singularities
$1$ and $\infty$. Having found the action of this group on the manifold of monodromy data, he
obtained connection formulae for asymptotic expansions of general solutions of $P_6$ as $t\to0$,
$t\to1$, and $t\to\infty$. Subsequently, McCoy and Tang~\cite{MT1,MT2,MT3} used Jimbo's monodromy
parametrization for small-$t$ asymptotics of $P_5$ to connect them with the large-$t$ asymptotic
expansions for the same equation.

Like Jimbo~\cite{J}, we consider here the cluster on the $P_6$ background; however, unlike the
{\it hard} isomonodromy condition used by Jimbo, we introduce a {\it weak} one, i.e., the
deformations which preserve the generators of the monodromy group, forming an essential part of
the manifold of monodromy data, rather than the whole manifold as for the hard isomonodromy
deformations. As we see later, the condition that isomonodromy deformations are weak but
not hard implies their discretization. At the same time the weak isomonodromy condition gives
us an opportunity to substitute for the cluster not only the regular singularity (as in the
Jimbo case) but also an irregular one. In this work we consider two different clusters on the
$P_6$ background: the first one described by an irregular singularity and the other by a
regular one. Even the ``regular'' case studied here differs from that in \cite{J}; in particular,
it is related with a formal limit $P_6\to P_5$, the``irregular'' case leads to another formal
limit $P_6\to P_5$. The main goal of this paper is to give a proper asymptotic interpretation
of the formal limits in the manner explained above.

We will study the first limit passage of the following oriented graph,

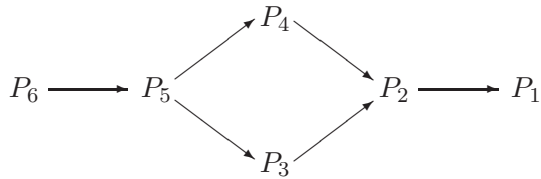
\begin{figure}[ht]
\begin{center}
\begin{picture}(260,65)
\put(30,35){$P_6$}
\put(45,38){\vector(1,0){30}}
\put(80,35){$P_5$}
\put(93,42){\vector(4,3){30}}
\put(93,34){\vector(4,-3){30}}
\put(125,62){$P_4$}
\put(125,07){$P_3$}
\put(138,64){\vector(4,-3){30}}
\put(138,11){\vector(4,3){30}}
\put(170,35){$P_2$}
\put(185,38){\vector(1,0){30}}
\put(220,35){$P_1$}
\end{picture}
\end{center}
\caption{The Painlev\'e - Okamoto degeneration scheme}
\label{fig:1}
\end{figure}
which represents the successive degeneration of the Painlev\'e equations ($P_n$). This
degeneration scheme was obtained by Painlev\'e~\cite{P},  who showed that by substituting
into the canonical  Painlev\'e equations a parameter $\epsilon$ and transforming  the
subsequent equations via appropriate $\epsilon\to0$ limits, simpler Painlev\'e equations
can be obtained\footnote[3]{The scheme was extensively studied by Okamoto, who
also presented it as the graph on Figure 1.}. Okamoto~\cite{O} pointed out that these formal limits can be viewed  as
infinitesimal canonical transformations of the Hamiltonian  systems associated with the
Painlev\'e equations; thus he equipped the scheme of Figure 1 with the analogous scheme
for Hamiltonians. Furthermore, Okamoto associated  with this step-by-step degeneration process
a step-by-step confluence scheme for certain scalar linear ODEs, whose isomonodromy deformations
are governed by the Painlev\'e equations. These scalar equations are known to be intimately
connected with matrix equations of the form~(\ref{eq:cluster-definition}). The degeneration
scheme for the linear ODEs can be formulated in a manner analogous to that shown in Figure 1,
in terms of a special symbol which represents the corresponding transformation of the
Poincar\'e ranks of singularities of the linear ODEs.

The authors works~\cite{K1,K2,KK} in conjunction with the present one can be understood as
the isomonodromy regularization of the degeneration scheme of Figure~\ref{fig:1}. In addition to the
scheme discussed above, our regularization suggests the following scheme of transformation
for the manifolds of monodromy data of the associated linear matrix ODEs:

\begin{figure}[ht]
\begin{center}
\begin{picture}(260,65)
\put(19,35){${\cal M}_6$}
\put(39,38){\vector(1,0){30}}
\put(74,35){${\cal M}_5$}
\put(93,42){\vector(4,3){30}}
\put(93,34){\vector(4,-3){30}}
\put(125,62){${\cal M}_4$}
\put(125,07){${\cal M}_3$}
\put(143,64){\vector(4,-3){30}}
\put(143,11){\vector(4,3){30}}
\put(174,35){${\cal M}_2$}
\put(193,38){\vector(1,0){30}}
\put(228,35){${\cal M}_1$}
\put(110,-10){Figure 2}
\end{picture}
\end{center}
\label{fig:2}
\end{figure}

Let us now discuss the schemes in Figures~\ref{fig:1} and \ref{fig:2} in  more detail. Each oriented
edge in the Figures means the existence of a corresponding limit passage; in fact, there could also
exist some other limits; for example, the limit $P_6\to P_5$ ($t\to1$), found by Painlev\'e~\cite{P},
differs from both of ours (in which $t\to0$). We call the formal limits equivalent if the
diagram in Figure~3 is commutative, where the $x$-arrows denote the action of the
transformation groups on the solutions of the Painlev\'e equations
(for $n=6$, $m=5$, see~\cite{O1,O2}) or the manifolds of monodromy data, the $y$-arrows denote
the same formal limits as on Figures~\ref{fig:1} and \ref{fig:2}, respectively, and where the
$z$-arrows ``enumerate'' the solutions by the points of the manifolds of monodromy data, i.e., they
``point out'' on solutions of the corresponding direct monodromy problems (see Section~\ref{sec:2}).

\begin{figure}[ht]
\begin{center}
\begin{picture}(100,85)
\put(28,65){$P_n$}
\put(42,68){\vector(1,0){29}}
\put(74,65){$P_n$}
\put(24,25){${\cal M}_n$}
\put(42,28){\vector(1,0){29}}
\put(74,25){${\cal M}_n$}
\put(30,63){\vector(0,-1){28}}
\put(80,63){\vector(0,-1){28}}
\put(08,53){$P_m$}
\put(22,56){\vector(1,0){29}}
\put(54,53){$P_m$}
\put(03,13){${\cal M}_m$}
\put(22,16){\vector(1,0){29}}
\put(54,13){${\cal M}_m$}
\put(10,51){\vector(0,-1){28}}
\put(60,51){\vector(0,-1){28}}
\put(26,64){\vector(-3,-2){10}}
\put(28,24){\vector(-3,-2){10}}
\put(78,24){\vector(-3,-2){10}}
\put(75,64){\vector(-3,-2){10}}
\put(10,-10){Figure 3}
\end{picture}
\hspace{4cm}
\begin{picture}(100,85)
\put(28,65){$P_n$}
\put(42,68){\vector(1,0){29}}
\put(74,65){$P_k$}
\put(24,25){${\cal M}_n$}
\put(42,28){\vector(1,0){29}}
\put(74,25){${\cal M}_k$}
\put(30,63){\vector(0,-1){28}}
\put(80,63){\vector(0,-1){28}}
\put(08,53){$P_m$}
\put(22,56){\vector(1,0){29}}
\put(54,53){$P_l$}
\put(03,13){${\cal M}_m$}
\put(22,16){\vector(1,0){29}}
\put(54,13){${\cal M}_l$}
\put(10,51){\vector(0,-1){28}}
\put(60,51){\vector(0,-1){28}}
\put(26,64){\vector(-3,-2){10}}
\put(28,24){\vector(-3,-2){10}}
\put(78,24){\vector(-3,-2){10}}
\put(75,64){\vector(-3,-2){10}}
\put(10,-10){Figure 4}
\end{picture}
\end{center}
\end{figure}
Note that if the transformation $P_n\to P_n$ acts not only between the initial Painlev\'e
equations but also during the entire limiting process, i.e., for arbitrary
$\epsilon\to0$\footnote[4]{More precisely, if we denote by $\varepsilon$ and $\varepsilon_1$ small
parameters in two copies of $P_n$, then the transformation maps $\varepsilon\to\varepsilon_1$, such
that the functions $\varepsilon(\varepsilon_1)\to0$ and $\varepsilon_1(\varepsilon)\to0$, when
$\varepsilon_1$ or, respectively, $\varepsilon\to0$.}, then the entire ${\cal M}$-plane of Figure 3 is not needed , since this
commutativity is valid automatically. In general, there are some transformations of the Painlev\'e
equations which do not have this property, i.e., they cannot be extended for arbitrary $\varepsilon$,
we present the diagram in extended form. Although, as explained above, the limits $P_6\to P_5$
that we study in this work have different asymptotic senses, they appear to be equivalent under the
above definition; moreover, all limits $P_6\to P_5$ known as at the time of the presents are
equivalent. Therefore, the following question arises naturally: {\it are there an $n$ and
$m$ such that nonequivalent formal limits $P_n\to P_m$ exist?}

There are some Painlev\'e equations which are known to be equivalent; for example, $P_5'$
($P_5$ with $\delta_5=0$, see (\ref{eq:P5})) is equivalent to $P_3$ \cite{G} and $P_{34}$ is
equivalent to $P_2$ \cite{I}. This equivalence of the equations means that there is a mapping
which is invertible on transcendental solutions and birationally dependent of them and their
derivatives. Along with the limit $P_4\to P_2$ (see Figure 1), the formal limit $P_4\to P_{34}$
was found in \cite{K1}. Are these limits ``substantially'' different? \\
To formulate the general notion of equivalence of the formal limits of different Painlev\'e
equations, i.e., the equivalence of the limits $P_n\to P_m$ and $P_k\to P_l$, we consider
Figure~4. In this figure the upper $x$-arrows denote the birational mappings (which are assumed to
exist). The mappings ${\cal M}_n\to{\cal M}_k$ and ${\cal M}_m\to{\cal M}_l$ are uniquely defined
by demanding commutativity for the $xz$-walls of the cube. The $y$-arrows denote the equivalence
classes of the limits defined earlier via Figure 3. The $z$-arrows have the same meaning as that
in Figure 3. {\it The limits $P_n\to P_m$ and $P_k\to P_l$ are said to be equivalent if the diagram
in Figure~{\rm4} is commutative. Are there any nonequivalent limits for pairwise equivalent
equations?} In particular, are the limits $P_4\to P_2$ and $P_4 \to P_{34}$ equivalent? Although
according to \cite{JK} these limits show different actions on local expansions for solutions of $P_4$,
they can appear to be equivalent under the above definition.

The cluster point of view makes clear that there must exist limits other than those shown in
Figure~\ref{fig:1}; e.g., along with the chain of successive limits $P_4\to P_2\to P_1$, there exist
direct ones \cite{K1,K2}, which look, at first glance, to be different from the limit obtained
via the successive procedure.
Thus we can expect the existence of a more complicated limit structure (than that shown on
Figure~\ref{fig:1}), which can be presented as the following oriented spatial graph:

\begin{figure}[ht]
\begin{center}
\begin{picture}(270,140)
\put(30,70){\circle*{3}}
\put(15,70){$P_6$}
\put(220,70){\circle*{3}}
\put(225,70){$P_1$}
\put(30,70){\vector(1,0){60}}
\put(90,70){\line(1,0){130}}
\put(30,70){\vector(2,1){45}}
\put(75,92){\line(2,1){50}}
\put(30,70){\vector(2,-1){50}}
\put(77,47){\line(2,-1){50}}
\put(30,70){\vector(3,-1){44}}
\put(63,59){\line(3,-1){36}}
\put(84,81){\vector(4,1){10}}
\put(30,70){\line(5,1){54}}
\put(92,83){\line(5,1){57}}
\put(148,94){\vector(3,-1){36}}
\put(184,82){\line(3,-1){38}}
\put(125,23){\vector(1,3){13}}
\put(148,93){\line(-1,-3){12}}
\put(135,107){\line(1,-1){13}}
\put(125,117){\vector(1,-1){12}}
\put(125,22){\vector(2,1){45}}
\put(159,39){\line(2,1){60}}
\put(220,70){\line(-2,1){45}}
\put(125,117){\vector(2,-1){50}}
\put(101,47){\vector(1,3){12}}
\put(105,59){\line(1,3){19}}
\put(101,46){\line(5,1){60}}
\put(151,56){\vector(4,1){10}}
\put(159,58){\line(5,1){60}}
\put(101,47){\vector(1,-1){12}}
\put(112,36){\line(1,-1){11}}
\put(101,47){\circle*{3}}
\put(88,54){$P_5$}
\put(125,117){\circle*{3}}
\put(125,122){$P_4$}
\put(125,23){\circle*{3}}
\put(125,08){$P_3$}
\put(101,47){\vector(1,1){18}}
\put(113,60){\line(1,1){34}}
\put(148,94){\circle*{3}}
\put(149,80){$P_2$}
\put(-30,-10){Figure 5: A spatial degeneration scheme for the Painlev\'e equations}
\end{picture}
\end{center}
\label{fig:5}
\end{figure}

Furthermore, it seems that the vertices $P_2$ and $P_3$ could have a more complicated
microlocal structure $P_2\sim P_2\rightleftarrows P_{34}$ and
$P_3\sim P_3\rightleftarrows P_5'\rightleftarrows P_3'$ where $P_3'$ is a special case
of $P_3$ defined by the following conditions: $\gamma\delta=0$ and $|\gamma|+|\delta|\neq0$,
where $\gamma$ and $\delta$ are the coefficients of the canonical form of $P_3$ \cite{I}.
If this is the case, then some of the edges which are incident to $P_3$ must be paired off
in order to be incident to $P_3$ and $P_3'$. The answers to our ``equivalence'' questions
(see above) could lead to a further complication of the diagram in Figure 5. The answer to
the first question governs the number of its edges, and the answer to the second question could
result in additional vertices and edges (say, $P_n\to P_{34}$ and $P_n\to P_5'$). In
Figure 5 we observe  that for most of the pairs of vertices there are several routes
connecting them.  The natural question therefore is: {\it are there any nonequivalent
routes for some pair of vertices?} The notion of equivalent routes, $P_n\to P_k\to P_l$
and $P_n\to P_m\to P_l$, can also be formulated by means of the commutativity of the diagram
in Figure~4. For this purpose we must identify the roles of the $x$- and $y$-arrows,
understanding them as denoting the equivalence classes of the limits. For routes with
three or four edges this definition must, naturally, be generalized.

As mentioned, the original setting of the formal limits for ODEs of the Painlev\'e
type includes the small parameter $\varepsilon$. The appearance of $\varepsilon$
endows the problem with an infinite-dimensional status, as the constants of integration
can be considered as arbitrary functions of $\varepsilon$. Furthermore, it is possible
that the latter functions do not have a limit as $\varepsilon\to0$, although the corresponding
limiting equations exist and therefore possess no explicit $\varepsilon$ dependence. Thus,
we see that for a proper understanding of such limits on the level of the solutions a
{\it regularization procedure} may be required. The regularization procedure, which we adopt
in this work, means that we will be making such choices of the functions mentioned above so
that they would satisfy the weak isomonodromy condition. This requirement leads us
(see \cite{K1,K2,KK} and Sections~\ref{sec:4} and \ref{sec:5}) to the
{\it discretization} of $\varepsilon$. It is clear that one could propose different
regularization procedures for the formal limits: different procedures depend predominantly
on possible applications and/or the employed mathematical machinery. Of course, some
regularizations could require no $\varepsilon$-discretization.
In particular we mean the approach proposed by Joshi and Kruskal \cite{JK}. It is necessary
to mention at once that not all the details are set out explicitly in \cite{JK}; therefore,
\cite{JK} yields a variety of different possibilities: here, we consider one such
possibility. Consider the Laurent expansion about a neighborhood of a pole for a solution
of some Painlev\'e equation. It has the following parameters: (coefficients of the Painlev\'e
equation; (2) one more constant of integration along with the position of the pole;
and (3) dependent and independent variables, i.e., the Painlev\'e function itself and its
argument. It is possible to define transformation of all the above parameters (at least,
again, as formal Laurent-type $\varepsilon$-expansions near $\varepsilon=0$) which map the
initial Laurent expansion for the solution of the Painlev\'e equation  into a new one. The
latter, in turn, can be considered as a pole expansion of the solution of some other
(limiting) Painlev\'e equation. All the parameters of the new pole expansion are defined by
the coefficients of the leading terms of the formal $\varepsilon$-expansions. These leading
terms are precisely the formulae given by Painlev\'e in his original work \cite{P}. So in \cite{JK}
no $\varepsilon$-discretization  and any associated linear structures and thereby
manifolds of monodromy data are involved.

We can provide Joshi's and Kruskal's work \cite{JK} with, perhaps, a somewhat unexpected
interpretation, which yields an opportunity to relate their method with isomonodromy
deformations. This interpretation lead us, again, to a discretization of $\varepsilon$;
but in this case, in a way different from that proposed in \cite{K1}, and below in
Sections~\ref{sec:4} and \ref{sec:5}. Recall that all the Painlev\'e equations are known
to possess transcendental solutions having an infinite number of poles, $\{t_n\}$, on
the real axis accumulating at infinity. The above $\varepsilon$-transformation can be
interpreted as a procedure for driving the poles of solutions of the Painlev\'e equation
to $\infty$. We fix a particular solution, which has an infinite number of poles accumulating
at infinity, and treat the ``pole-drive'' as a pole-to-pole jumping. This leads to a
discretization of $\varepsilon$ ($\varepsilon\to\varepsilon_n$), because we are aloud now
to do only discrete shifts of the poles $t_n\to t_{n+1}$. Simultaneously iterating the
pole expansion via B\"acklund transformations by the number of times consistent with the
$\varepsilon$-transformation, we arrive at the pole expansion of some particular solution of
the limiting Painlev\'e equation. The consistency means that we take the number of B\"acklund
transformations such that the coefficients of the corresponding Painlev\'e equation increase
by a rate prescribed by the formulae for $\varepsilon$-transformation
(with $\varepsilon\to\varepsilon_n$). In other words, on the level of the pole expansions,
we can interpret the limiting passage between Painlev\'e equations as some special
asymptotics of poles of (fixed) solutions. These asymptotics, in turn, are given in terms of
a pole of some particular solution of the limiting Painlev\'e equation. Note that our
interpretation of the limiting passage is consistent with isomonodromy deformations.
We know that in the standard situation (without B\"acklund iterations) pole-asymptotics of
solutions can be parameterized by the monodromy data (see \cite{IN}). {\it Is that possible
to get an analogous parametrization of the poles in the new situation described above
\footnote[5]{The latter parametrization, in fact, means that the monodromy data of the initial
Painlev\'e equation are related with the monodromy data of the limiting one.}?}

It is important to make a distinction with our previous works~\cite{K1,K2,KK} which were
concerned with the WKB-method and where the clusters of turning points were considered,
and this work, where no WKB-method is involved. Hence, now, we are considering a different
asymptotic process. In general, for Painlev\'e equations all limiting procedures are related
with various processes of merging turning points and singularities in the associated linear
ODEs describing corresponding isomonodromy deformations. For example, the limit $P_4\to P_2$
could be organized by considering the following two possibilities for the asymptotic behavior
of the $2\times2$ matrix linear ODE associated with $P_4$:
(1) merging of four turning points; and (2) merging of the irregular and regular singular points,
i.e., in our terminology, we can consider either an isomonodromy cluster of four turning points or
an isomonodromy cluster of irregular and regular singularities. So it seems that, really, nonequivalent
limits could exist, and the graph on Figure~5 could actually be incomplete.

The discretization of $\varepsilon$ is not too important a feature of the technique we are developing.
Actually, in the isomonodromy systems with more than one continuous variable, the so-called higher
Painlev\'e equations or Garnier systems, we can observe the analogous asymptotic processes and apply
the same technique without any $\varepsilon$-discretization~\cite{K3}. We mention also an interesting
WKB-theory for special solutions of the Painlev\'e equations and its interrelation with the corresponding
WKB objects for the associated scalar linear ODEs~\cite{KT}.

The technique considered in \cite{K1,K2,KK} and here has an interesting scope of applications. One meets
the clusters of turning points while studying via the isomonodromy approach double-scaling limits of
partition functions for matrix models of $2D$ quantum gravity~\cite{FIK1,FIK2}. Clusters of stationary
phase (saddle) points, which are similar to the clusters of turning points, appear in the description of
caustics in $1+1$ systems integrable via the Inverse Scattering Transform~\cite{K3}. The clusters of regular
singularities play an important role in the study of a zero-curvature limit for holonomic quantum field
theory of bosons in the Poincare disc\cite{NT,PBT}\footnote[6]{Cited works are written from a different
perspective}. The limit relates this theory with the original Euclidean Sato-Miwa-Jimbo
theory~\cite{SMJ1,SMJ2,SMJ3,SMJ4,SMJ5,JMMS}. The correspondence of two-point correlators of these theories in
our notation reads as $P_6\to P_3$ or $P_6\to P_5'$. Considering Figure~5 we see that if we add to the
limit $P_6\to P_5$ the next one $P_5\to P_3$, then we obtain the desired limit $P_6\to P_3$. The limit
$P_5\to P_3$ in our approach is also related with the isomonodromy cluster of two regular singularities
but on the ``irregular background'' ($A(\lambda,t)=const\cdot\sigma_3\neq0$). Of course, there is a direct
limit $P_6\to P_3$, and the question of equivalence discussed above should be studied. I hope to return to
this question in a special publication.

This paper is organized as follows:
Section~\ref{sec:2} contains no new material. It is based mainly on the works~\cite{J,JMU,JM}. In this
Section we set the notation related with the description of isomonodromy deformations for $P_6$ and $P_5$.
Most of the formulae presented here are extensively used throughout Sections~\ref{sec:3}-\ref{sec:5}.
Some of our definitions are slightly different from those of the cited works; in particular, I found
it convenient to define two different manifolds of monodromy data for $P_5$ (${\cal M}_5$ and
$\tilde{\cal M}_5$): they are related with the schemes of the calculations in Sections~\ref{sec:4}
and \ref{sec:5}, namely, constructions of semi-infinite sequences of Schlesinger transformations, which
are important issues in our approach.

In Section~\ref{sec:3}, we follow the simple procedure explained at the beginning of the Introduction to
derive two different and novel formal limits $P_6\to P_5$. The corresponding formulae for the canonical
Painlev\'e functions are rather cumbersome, whilst those being rewritten for the corresponding $\tau$-functions
become much simpler.

In Section~\ref{sec:4} we derive a discrete regularized version of the first limit. We calculate the
manifold of the monodromy data for the limiting $P_5$ transcendent in terms of the corresponding manifold
for the initial $P_6$ transcendent. The result is stated in Theorem~\ref{th:1}, and the following part of
the section presents details of the derivation.

Theorem~\ref{th:2} of Section~\ref{sec:5} states results for a discrete regularization of the second limit
analogous to those in Section~\ref{sec:4}; however, a derivation of these results is more complicated than
those in Section~4, because, in Section~\ref{sec:5} we meet an additional problem with the normalization
of {\it a priori} unknown function $\Psi_5$. It causes a number of additional technical detail one of them
is considered in Appendix. As a by-product for these efforts, we get:
(1) a clear picture of how the monodromy matrices of the merging regular singularities produce the
Stokes multipliers of the resulting $\Psi_5$-function which has an irregular singularity describing
the merging process; and
(2) the latter derivation is easy to generalize for a background with an arbitrary number of regular
singularities.

In the case of $P_6$, when Equation~(\ref{eq:cluster-definition}) has four regular singularities, the
cluster and exterior cluster domains can be mapped one into another via a fractional-linear transformation.
The latter means a transformation which maps the cluster of Section~\ref{sec:4} into the one of
Section~\ref{sec:5}. Thus we prove the equivalence of the limits for the Painlev\'e transcendents and the
corresponding $\tau$-functions.

In the Appendix we study the problem of the simultaneous reduction of a pair of ${\rm SL}(2,\mathbb C)$
matrices to the lower and upper triangular forms. The solution is stated in a theorem which refers to a
number of special cases studied in corresponding propositions. These results have an important consequence
on the solvability of the asymptotic problem studied in Section~(\ref{sec:5}). A generalization of the
problem of the simultaneous transformation of $n\geq3$ ${\rm SL}(2,\mathbb C)$ matrices to triangular
forms is important for the investigation of the isomonodromy clusters with $n$ regular singularities.

\section{$P_6$ and $P_5$ as Isomonodromy Deformations}
 \label{sec:2}
In this section we recall some facts from the isomonodromy theory of $P_6$ and $P_5$ following, mainly, the
works~\cite{J,JMU,JM}. To make a distinction between analogous objects related to $P_n$ $(n=5,6)$ for
different $n$, we supply all of them by the subscript $n$. The convenience of this agreement becomes
evident in the following sections.

Consider the $2\times2$ matrix linear ODE with four regular singularities at $\lambda=0,1,t,\infty$:
\begin{equation}
 \label{eq:linear-4-regular}
\frac{d}{d\lambda}\Psi_6=\left(\frac{A_{06}}\lambda+
\frac{A_{16}}{\lambda-1}+\frac{A_{t6}}{\lambda-t_6}\right)\Psi_6.
\end{equation}
It is assumed that
\begin{equation}
 \label{eq:niormalization-for-P6}
 A_{06}+A_{16}+A_{t6}=-\frac{\Theta_{\infty6}}2\sigma_3,\qquad\sigma_3=
 \left(\begin{array}{cc}1&0\\0&-1\end{array}\right),\qquad
 \Theta_{\infty6}\in\mathbb{C}\setminus\mathbb{Z},
\end{equation}
and there exist $R_{\nu6}\in{\rm SL}(2,\mathbb{C})$ such that
\begin{equation}
 \label{eq:Rnu6}
 R_{\nu6}^{-1}A_{\nu6}R_{\nu6}=\frac{\Theta_{\nu6}}2\sigma_3,\qquad
 \Theta_{\nu6}\in\mathbb{C}\setminus\mathbb{Z}.
\end{equation}
In a neighborhood of the regular singularity $\lambda=\nu\;(\nu\neq\infty)$, the $\Psi_6$-function can be
expanded as
\begin{equation}
 \label{eq:psi-expansion-at-nu}
 \Psi_6\,\underset{\lambda\to\nu}=\,\sum\limits_{k=0}^\infty\Psi_{k\nu6}
 (\lambda-\nu)^{k+\frac{\Theta_{\nu6}}2\sigma_3}C_{\nu6},
\end{equation}
where the matrices $\Psi_{k\nu6}$, $C_{\nu6}$ are independent of $\lambda$, because of the normalization
conditions (\ref{eq:niormalization-for-P6}) and (\ref{eq:Rnu6}), we can assume that
\begin{equation}
 \label{eq:condition-det-1}
 \Psi_6, R_{\nu6},\Psi_{0\nu6},C_{\nu6}\in{\rm SL}(2,\mathbb{C}),\qquad\nu=0,1,t.
\end{equation}
Because of Equation~(\ref{eq:niormalization-for-P6}), the expansion of the $\Psi_6$-function can be
normalized at $\lambda=\infty$:
\begin{equation}
 \label{eq:psi6-infty-normalization}
\Psi_6\underset{\lambda\to\infty}=\left(I+
\sum\limits_{k=1}^\infty\Psi_{k\infty6}\lambda^{-k}\right)\lambda^{-\frac{\Theta_{\infty6}}2\sigma_3}.
\end{equation}
The single-valued $\Psi_6$-function with the stated properties can be defined on the $\lambda$-plane
"cut" along the negative imaginary semi-axis $[-i\infty,0]$, the segment $[0,t_6]$, and the positive
real semi-axis $[1,+\infty]$. Henceforth, we assume that $t_6\in(0,1)$. Following \cite{J} we choose
the paths in $\bar{\mathbb{C}}\setminus\{0,t_6,1,\infty\}$ as is shown on Figure~6.
\begin{figure}[ht]
\begin{center}
\begin{picture}(310,90)
\put(70,30){\circle*{3}}
\put(68,36){$0$}
\put(120,30){\circle*{3}}
\put(112,36){$t$}
\put(190,30){\circle*{3}}
\put(176,36){$1$}
\put(300,30){\circle*{3}}
\put(286,36){$\infty$}
\put(60,90){\circle*{3}}
\put(45,85){$\lambda_0$}
\qbezier(60,90)(90,25)(70,20)
\qbezier(60,90)(45,25)(70,20)
\qbezier(60,90)(155,28)(130,20)
\qbezier(60,90)(94,22)(130,20)
\qbezier(60,90)(188,08)(214,20)
\qbezier(60,90)(247,30)(214,20)
\qbezier(60,90)(350,50)(320,25)
\qbezier(60,90)(310,10)(320,25)
\put(124,42){\vector(-1,1){10}}
\put(189,43){\vector(-2,1){10}}
\put(77,41){\vector(-1,4){3}}
\put(250,58){\vector(-4,1){10}}
\put(140,00){Figure 6}
\end{picture}
\end{center}
\end{figure}
Continuing the $\Psi_6$-function along these paths we find a multi-valued function with the
monodromy matrices $M_{\nu6}$ $(\nu=0,1,t,\infty)$ which, in terms of the local
expansions (\ref{eq:psi-expansion-at-nu}) and (\ref{eq:psi6-infty-normalization}), can be
written as follows:
\begin{equation}
 \label{eq:monodromy-local}
M_{\nu6}=C_{\nu6}^{-1}e^{\pi i\Theta_{\nu6}\sigma_3}C_{\nu6},\qquad\nu=0,1,t,\infty,\quad C_{\infty6}=I.
\end{equation}
One proves the cyclic relation
\begin{equation}
 \label{eq:cyclic6}
M_{\infty6}M_{16}M_{t6}M_{06}=I.
\end{equation}
The monodromy group $({\cal MG}_6)$ is the subgroup of ${\rm SL}(2,\mathbb{C})$ generated by the matrices
$M_{\nu6}$. The set of monodromy data $({\cal M}_6)$ is an ordered set of matrix elements of $M_{\nu6}$
completed with four complex parameters $\Theta_{\nu6}$ $(\nu=0,1,t,\infty)$, satisfying  the equation
\begin{equation}
 \label{eq: trace-M}
{\rm tr}M_{\nu6}=2\cos(\pi\Theta_{\nu6}).
\end{equation}
The numbers $\Theta_{\nu6}$ are called the coefficients of formal monodromy.

Thus we see that for a fixed parameter $t_6$ and matrices $A_{\nu6}$ the above procedure yields the
unique set ${\cal M}_6$. One proves that a different set of $\{A_{\nu6}\}_{\nu=0,1,t}$, with the
same $t_6$, defines a different set ${\cal M}_6$; therefore, considering the Riemann problem of
reconstruction of all differential equations corresponding to the given set ${\cal M}_6$, we have to
vary not only $A_{\nu6}$, but also $t_6$. Following the Schlesinger approach it is convenient to consider
$t_6$ as an independent variable, then $A_{\nu6}=A_{\nu6}(t_6)$. To find these functions explicitly one
has to notice that the {\it hard isomonodromy condition}, i.e., $\partial_{t_6}{\cal M}_6=0$, leads to
the additional ODE with respect to $t_6$ for the function $\Psi_6$, which, under the assumptions on the
formal monodromy (\ref{eq:niormalization-for-P6}) and (\ref{eq:Rnu6}), reads
\begin{equation}
 \label{eq:t-6}
\frac{\partial}{\partial t_6}\Psi_6=-\frac{A_{t6}}{\lambda-t_6}\Psi_6.
\end{equation}
The compatibility condition of Equations~(\ref{eq:linear-4-regular}) and (\ref{eq:t-6}) is
\begin{equation}
 \label{eq:schlesinger-p6}
 \frac{d A_{\nu6}}{dt_6}=\frac{[A_{t6},A_{\nu6}]}{t_6-\nu},\qquad\nu=0,1,
\end{equation}
where $[A,B]\equiv AB-BA$. In the derivation of (\ref{eq:schlesinger-p6}) we took into account the
normalization condition (\ref{eq:niormalization-for-P6}) together with $\partial_{t_6}\Theta_{\infty6}=0$.

The system~(\ref{eq:schlesinger-p6}) is called the {\it system of iso\-mo\-no\-dro\-my deformations}
or the Schle\-sin\-ger system. It consists of eight scalar nonlinear ODEs of the first order for the matrix
elements of $A_{\nu6}$ $(\nu=0,1)$. By using Equations~(\ref{eq:niormalization-for-P6}) and
(\ref{eq:Rnu6}), one finds the following five first integrals:
$$
{\rm tr}\,A_{\nu6}=0,\quad\det\,A_{\nu6}=-\frac{\Theta_{\nu6}^2}4,\;\nu=0,1,\quad
\det(\Theta_{\infty6}\sigma_3/2+A_{06}+A_{16})=-\Theta^2_{t6}/4,
$$
reducing the number of independent scalar equations to three. This number coincides with with the complex
dimension of the {\it manifold of monodromy data},
${\cal M}_6(\Theta_{06},\Theta_{16},\Theta_{t6},\Theta_{\infty6})$. The points of this manifold are in
one-to-one correspondence with the set ${\cal M}_6$ for the fixed parameters $\Theta_{\nu6}$. More
precisely: identifying an ordered set of the matrix elements of ${\cal M}_6$ as a point in
$\mathbb{C}^{14}$, one defines ${\cal M}_6(\Theta_{06},\Theta_{16},\Theta_{t6},\Theta_{\infty6})$ as an
algebraic variety in $\mathbb{C}^{14}$ by Equations~(\ref{eq:cyclic6}), (\ref{eq: trace-M}), and
$\det M_{\nu6}=1$, $(\nu=0,1,t,\infty)$.

Following \cite{JM}, introduce the notation,
\begin{equation}
 \label{eq:A-and-R}
A_{\nu6}=\left(\begin{array}{cc}
z_{\nu6}+\frac{\Theta_{\nu6}}2&-u_{\nu6}z_{\nu6}\\
\frac{z_{\nu6}+\Theta_{\nu6}}{u_{\nu6}}&-z_{\nu6}-\frac{\Theta_{\nu6}}2
\end{array}\right),\qquad
R_{\nu6}=\left(\begin{array}{cc}
\frac1{\Theta_{\nu6}s_{\nu6}}&u_{\nu6}z_{\nu6}s_{\nu6}\\
\frac1{\Theta_{\nu6}u_{\nu6}s_{\nu6}}&(z_{\nu6}+\Theta_{\nu6})s_{\nu6}
\end{array}\right),
\end{equation}
where $\nu=0,1,t$, and $s_{\nu6}\in\mathbb{C}\setminus\{0\}$. Equation~(\ref{eq:niormalization-for-P6})
in the notation~(\ref{eq:A-and-R}) reads:
\begin{gather}
 \label{eq:normalization-scalar-diag}
z_{06}+z_{16}+z_{t6}=-\frac12(\Theta_{06}+\Theta_{16}+\Theta_{t6}+\Theta_{\infty6}),\\
u_{06}z_{06}+u_{16}z_{16}+u_{t6}z_{t6}=0,\qquad\frac{z_{06}+\Theta_{06}}{u_{06}}+\frac{z_{16}+
\Theta_{16}}{u_{16}}+\frac{z_{t6}+\Theta_{t6}}{u_{t6}}=0.
\label{eq:normalization-scalar-off}
\end{gather}
Thus, only three parameters amongst those which define the coefficients of Equation~(\ref{eq:linear-4-regular})
are independent. In the hard isomonodromy case the parameters appear to be functions of $t_6$:
$z_{\nu6}=z_{\nu6}(t_6)$ and $u_{\nu6}=u_{\nu6}(t_6)$ defined by System~(\ref{eq:schlesinger-p6}). As
explained above, one can rewrite (\ref{eq:schlesinger-p6}) as a system of three scalar first-order ODEs.
A convenient form of this system is given in \cite{JM}. Here we define two important functions closely
related with System~(\ref{eq:schlesinger-p6}). The first one,
\begin{equation}
 \label{eq:y6}
y_6=\frac{t_6u_{06}z_{06}}{(t_6+1)u_{06}z_{06}+t_6u_{16}z_{16}+u_{t6}z_{t6}}=
\frac1{1+\left(1-\frac1{t_6}\right)\frac{u_{16}z_{16}}{u_{06}z_{06}}}=
\frac{t_6}{1+(1-t_6)\frac{u_{t6}z_{t6}}{u_{06}z_{06}}},
\end{equation}
is a solution of the sixth Painlev\'e equation,
\begin{align}
\frac{d^2y_6}{dt_6}^2&=\frac12\left(\frac1{y_6}+\frac1{y_6-1}+\frac1{y_6-t_6}\right)\left
(\frac{dy_6}{dt_6}\right)^2-\left(\frac1{t_6}+\frac1{t_6-1}+\frac1{y_6-t_6}\right)\frac{dy_6}{dt_6}
 \nonumber\\
&+\frac{y_6(y_6-1)(y_6-t_6)}{t_6^2(t_6-1)^2}\left(\alpha_6+\frac{\beta_6t_6}{y_6^2}+
\frac{\gamma_6(t_6-1)}{(y_6-1)^2}+\frac{\delta_6t_6(t_6-1)}{(y_6-t_6)^2}\right),
 \label{eq:p6}\\
&\alpha_6=\frac12(\Theta_{\infty6}-1)^2,\quad
\beta_6=-\frac12\Theta_{06}^2,\quad
\gamma_6=\frac12\Theta_{16}^2,\quad
\delta_6=\frac12(1-\Theta_{t6}^2).
 \label{eq:p6coeff}
\end{align}
The second function, which is important for applications, is the $\tau$-function~\cite{JM}:
\begin{equation}
 \label{eq:tau-6}
\frac{d}{dt_6}\log\,\tau_6(t_6)={\rm tr}\,\left(\frac{A_{06}}{t_6}+\frac{A_{16}}{t_6-1}\right)A_{t6}
\end{equation}
The function
\begin{equation}
 \label{eq:sigma-6}
\hat\sigma_6(t_6)=t_6(t_6-1)\frac{d}{dt_6}\log\tau_6(t_6)
\end{equation}
satisfies an ODE of the second order, quadratic with respect to
$\hat\sigma''_6(t_6)$~\cite{J,O,JM}\footnote[7]{In the cited papers, the authors use a definition of
the $\sigma$-function shifted by a linear function of $t_6$. The latter function satisfies a differential
equation symmetric with respect to the formal monodromies. We use a ``hat'' in our notation to make a
difference between these functions.}.
Conversely, starting from an arbitrary solution of Equation~(\ref{eq:p6}) and using the formulae
given in \cite{JM} one can construct a solution of the Schlesinger System~(\ref{eq:schlesinger-p6})
satisfying Equation~(\ref{eq:y6}). Therefore, all solutions of $P_6$ (\ref{eq:p6}) can be obtained from
some solution of the Schlesinger System (\ref{eq:schlesinger-p6}) via Equation~(\ref{eq:y6}).

Let us consider the {\it weak isomonodromy deformations} of Equation~((\ref{eq:linear-4-regular}) as
the deformations of $A_{\nu6}$ ($\nu,0,1,t$) preserving, up to the sign, the generators of ${\cal MG}_6$,
i.e., isomonodromy deformations in the sense of the ${\rm PSL}_2(\mathbb{C})$ monodromy group. Using
Equations~(\ref{eq:cyclic6}) and (\ref{eq: trace-M}) one proves that in the weak case the continuous
Schlesinger deformations can be extended only by the group of discrete Schlesinger transformations acting
on the formal monodromies as
\begin{equation}
 \label{eq:theta-shift}
 \Theta_{\nu6}\to\Theta_{\nu6}+n_\nu,\quad\nu=0,1,t,\infty,\qquad n_0+n_1+n_t+n_\infty=0({\rm mod}\,2),
\end{equation}
and the reflections
\begin{equation}
 \label{eq: theta-reflections}
\Theta_{\nu6}\to-\Theta_{\nu6},\qquad\nu=0,1,t.
\end{equation}
Consider how these transformations act on the coefficients $A_{\nu6}$ of Equation~(\ref{eq:linear-4-regular}).
First, notice that the reflections~(\ref{eq:theta-shift}) are actually related with a certain ambiguity in
the parametrization of the matrices $A_{\nu6}$ rather than with any transformation of these matrices. The
ambiguity in the parametrization is related with the ambiguity in writing the local expansion of the
$\Psi_6$ function near the corresponding singular point: $(\Psi_{k\nu6},\Theta_{\nu6},C_{\nu6})$ $\mapsto$
$(\Psi_{k\nu6}\sigma_1,-\Theta_{\nu6},\sigma_1C_{\nu6})$, where $\sigma_1$ is the Pauli matrix (with the
matrix elements $\sigma_1^{ij}=0$ if $i=j$, else $\sigma_1^{ij}=1$). This nonuniqueness leads, simply, to
the reparametrization of the matrices $A_{\nu6}$: $(u_{\nu6},\Theta_{\nu6},z_{\nu6})$ $\mapsto$
$(\tilde u_{\nu6},\tilde\Theta_{\nu6},\tilde z_{\nu6})$, where
$z_{\nu6}+\frac{\Theta_{\nu6}}2=\tilde z_{\nu6}+\frac{\tilde\Theta_{\nu6}}2,\tilde\Theta_{\nu6}=\Theta_{\nu6}$,
and $u_{\nu6}z_{\nu6}=\tilde u_{\nu6}\tilde z_{\nu6}$, which doesn't affect the functions $y_6(t_6)$
and $\hat\sigma_6(t_6)$ (see Equations (\ref{eq:y6}), (\ref{eq:tau-6}), and (\ref{eq:sigma-6})).

 Note that the reflection
\begin{equation}
 \label{eq:reflection-infinity}
\tilde\Theta_{\infty6}=-\Theta_{\infty6}
\end{equation}
does not preserve the generators of ${\cal MG}_6$. This reflection is related with the change of the
generators of ${\cal MG}_6$ to their inverse: $\widetilde{\cal MG}_6=\sigma_1{\cal MG}_6\sigma_1$.
To see this one can define the following action of the reflection (\ref{eq:reflection-infinity}) on
the $\Psi_6$ function, $\tilde\Psi_6=\sigma_1\Psi_6\sigma_1$; it yields:
\begin{equation}
 \label{eq:action-reflect-infty}
\tilde M_{\nu6}=\sigma_1M_{\nu6}\sigma_1,\qquad
\tilde A_{\nu6}=\sigma_1A_{\nu6}\sigma_1.
\end{equation}
It is an immediate consequence of Equations (\ref{eq:action-reflect-infty}), (\ref{eq:tau-6}), and
(\ref{eq:sigma-6}) that $\tilde{\hat{\sigma}}_6=\hat\sigma_6$. In terms of the matrix elements the last
equation in (\ref{eq:action-reflect-infty}) reads:
\begin{gather}
 \label{eq:trans-theta-reflect-infty}
\tilde\Theta_{\nu6}=\pm\Theta_{\nu6},\quad\nu=0,1,t,\\
\tilde z_{\nu6}=-\left(z_{\nu6}+\frac{\tilde\Theta_{\nu6}+\Theta_{\nu6}}2\right),\quad
\tilde u_{\nu6}=\frac1{u_{\nu6}}\cdot
\frac{z_{\nu6}+\Theta_{\nu6}}{z_{\nu6}+\frac{\tilde\Theta_{\nu6}+\Theta_{\nu6}}2},
\label{eq:trans-coef-reflect-infty}
\end{gather}
where in Equation (\ref{eq:trans-theta-reflect-infty}) for every value of $\nu$ we can make arbitrary
choices of the signs. One can use Equations (\ref{eq:trans-coef-reflect-infty}), (\ref{eq:y6}), and
Equations (C.49), (C.52), and (C.55) of \cite{JM} to find $\tilde y_6$ as a rational function of
$y_6$, $y_6'$, and $t_6$.

Consider the group of the discrete Schlesinger Transformations (\ref{eq:theta-shift}). Evidently, it has
four generators and it is isomorphic to the group of translations along special basis in
$\mathbb{C}^4$. We denote by ${\cal L}_{\nu\nu'}^{\pm\pm}$, where $\nu\neq\nu'$, the elementary
Schlesinger Transformations~\cite{JM}:
$$
{\cal L}_{\nu\nu'}^{\pm\pm}:\quad \tilde\Theta_{\varkappa6}=\Theta_{\varkappa6}\pm1,\quad
{\rm for}\quad\varkappa=\nu,\,\nu';
\qquad\tilde\Theta_{\varkappa6}=\Theta_{\varkappa6}\quad{\rm for}\quad\varkappa\in\{0,1,t,\infty\}
\setminus\{\nu,\,\nu'\}.
$$
In ${\cal L}_{\nu\nu'}^{\pm\pm}$ we choose the signs over $\kappa=\nu,\,\nu'$ in the same way as in the formula
for $\tilde\Theta_{\varkappa6}$. Since $\tilde M_{\varkappa6}=\pm M_{\varkappa6}$ (the sign minus occurs only for
$\varkappa=\nu,\,\nu'$), the action of ${\cal L}_{\nu\nu'}^{\pm\pm}$ on the $\Psi_6$-function can be
defined as the left multiplication
\begin{equation}
 \label{eq:L-dressing}
 \tilde\Psi_6={\cal L}_{\nu\nu'}^{\pm\pm}\Psi_6,
\end{equation}
where
\begin{gather*}
L_{\nu\nu'}^{\pm\pm}=\frac{\sqrt{\lambda-\nu'}}{\sqrt{\lambda-\nu}}J_{\nu\nu'}^{\pm\pm}+
\frac{\sqrt{\lambda-\nu}}{\sqrt{\lambda-\nu'}}J_{\nu'\nu}^{\pm\pm},\qquad\nu,\,\nu'=0,1,t,\\
L_{\nu\infty}^{\pm\pm}=\sqrt{\lambda-\nu}\sigma_\infty^\pm+\frac1{\sqrt{\lambda-\nu}}J_{\nu\infty}^{\pm\pm},
\quad\nu\neq\infty,\quad\sigma_\varkappa^+=\left(\begin{array}{cc}0&0\\0&1\end{array}\right),\;\;
\sigma_\varkappa^-=\left(\begin{array}{cc}1&0\\0&0\end{array}\right).
\end{gather*}
Here, the branches of the roots are fixed as $\sqrt{\lambda-\nu}/\sqrt\lambda\to1$ and
$\sqrt{\lambda-\nu}/\sqrt{\lambda-\nu'}\to1$ as $\lambda\to\infty$ and $\lambda$ belongs to the complex
plane cut as explained above (see the paragraph right after Equation~(\ref{eq:psi6-infty-normalization})).
The choice of the signs over all equal subscripts is the same (upper/lower). The matrices
$J_{\nu\nu'}^{\pm\pm}$ for $\nu,\nu'\neq\infty$ are uniquely defined by the equations:
$$
J_{\nu\nu'}^{\pm\pm}+J_{\nu'\nu}^{\pm\pm}=I,\qquad J_{\nu\nu'}^{\pm\pm}R_{\nu6}\sigma_\nu^{\mp}=
\left(\begin{array}{c}0\\0\end{array}\right).
$$
The result is as follows
\begin{equation}
 \label{eq:J-and-Delta}
J_{\nu'\nu}^{\pm\pm}=\frac1{\Delta_{\nu'\nu}^{\pm}}
\left(\begin{array}{cc}b_\nu^{\pm}&0\\0&-a_\nu^{\pm}\end{array}\right)
\left(\begin{array}{cc}a_{\nu'}^{\pm}&b_{\nu'}^{\pm}\\a_{\nu'}^{\pm}&b_{\nu'}^{\pm}\end{array}\right),
\qquad\Delta_{\nu'\nu}^{\pm}=a_{\nu'}^{\pm}b_\nu^\pm-a_\nu^{\pm}b_{\nu'}^\pm,
\end{equation}
where $a_\nu^{\pm}$ and $b_\nu^{\pm}$ are different notations for the matrix elements of $R_{\nu6}$
(see Equation~(\ref{eq:Rnu6})) which are convenient here, namely,
\begin{equation}
 \label{eq:Rnu6a-b}
R_{\nu6}=\left(\begin{array}{cc}b_{\nu}^+&b_{\nu}^-\\-a_{\nu}^+&-a_{\nu}^-\end{array}\right).
\end{equation}
The matrix $J_{\nu\nu'}^{\pm\pm}$ can be found by the same Equations~(\ref{eq:J-and-Delta}) and
(\ref{eq:Rnu6a-b}) by making the permutation of the subscripts $(\nu'\leftrightarrow\nu)$. Note also
the following useful properties of $J_{\nu\nu'}^{\pm\pm}$:
\begin{equation}
 \label{eq:J-projectors}
\big(J_{\nu\nu'}^{\pm\pm}\big)^2=J_{\nu\nu'}^{\pm\pm},\qquad
J_{\nu\nu'}^{\pm\pm}J_{\nu'\nu}^{\pm\pm}=0,
\end{equation}
and analogous equations with $\nu\leftrightarrow\nu'$.

We see that Transformation~(\ref{eq:L-dressing}) exists iff $\Delta_{\nu\nu'}^\pm\neq0$. Using Equations
(C.51), (C.52), and (C.55) of \cite{JM} one finds that the condition $\Delta_{\nu\nu'}^\pm=0$ is equivalent
to the existence of a one-parameter solution $y_6$ of $P_6$, which solves an ODE of the first order,
\begin{equation}
 \label{eq:first-order-y6}
\frac{dy_6}{dt_6}=R(y_6,t_6),
\end{equation}
where $R$ is a rational function of its arguments with the coefficients defined by $\Theta_{\nu6}$
$(\nu=0,1,t,\infty)$. Thus, for general continuous Schlesinger deformations (\ref{eq:schlesinger-p6})
which correspond to transcendental (nonclassical) solutions of $P_6$, the condition
\begin{equation}
 \label{eq:Delta-ne0}
 \Delta_{\nu\nu'}^\pm\neq0
\end{equation}
is valid. Furthermore, we consider iterations of Transformations~(\ref{eq:L-dressing}). If in some step we
find $\Delta_{\nu\nu'}^\pm=0$, then it means we start from a solution of $P_6$ which is the iteration of
a special (classical) solution of (\ref{eq:p6}), i.e., it can be presented in the form
$\tilde R(y_6,y_6',t_6)$, where $\tilde R$ is a rational function of its arguments with the coefficients
defined by $\Theta_{\nu6}$ and where $y_6$ is a solution of an equation of the type (\ref{eq:first-order-y6}).
Thus we can iterate general (transcendental) solutions of $P_6$ (\ref{eq:p6}) without any restrictions.
The condition (\ref{eq:Delta-ne0}) is assumed throughout this paper. It would be interesting to perform
the complete investigation of our problem including the special (classical) solutions of
Equation~(\ref{eq:p6}). This investigation is in progress now \cite{U,W}.

Let us also find the action of ${\cal L}_{\nu\nu'}^{\pm\pm}$ on the matrices $A_{\nu6}$: substituting
transformation~(\ref{eq:L-dressing}) into Equation~(\ref{eq:p6}) written for $\tilde A_{\nu6}$ and
$\tilde\Psi_6$, one obtains:
\begin{gather}
 \label{eq:Atilde-nu}
J_{\nu'\nu}^{\pm\pm}\left(A_{\nu'6}-\frac12\right)=\tilde A_{\nu'6}J_{\nu'\nu}^{\pm\pm},\qquad
J_{\nu\nu'}^{\pm\pm}\left(A_{\nu6}-\frac12\right)=\tilde A_{\nu6}J_{\nu\nu'}^{\pm\pm},\\
\tilde A_{\mu6}=\left(J_{\nu'\nu}^{\pm\pm}+\frac{\nu-\mu}{\nu'-\mu}J_{\nu\nu'}^{\pm\pm}\right)A_{\mu6}
\left(J_{\nu'\nu}^{\pm\pm}+\frac{\nu'-\mu}{\nu-\mu}J_{\nu\nu'}^{\pm\pm}\right),\qquad\mu\neq\nu,\nu'.
\label{eq:Atilde-mu}
\end{gather}
Using these equations and  noting that the set $\{\nu,\nu',\mu\}$ is the permutation of $\{0,1,t\}$ we find
\begin{equation}
 \label{eq:Atilde-infty}
\tilde A_{\nu6}+\tilde A_{\nu'6}+\tilde A_{\mu6}=-\frac{\Theta_{\infty6}}2\sigma_3.
\end{equation}
Multiplying Equation~(\ref{eq:Atilde-infty}) by $J_{\nu\nu'}^{\pm\pm}$ on the right and using Equations
(\ref{eq:J-projectors}), (\ref{eq:Atilde-nu}), and (\ref{eq:Atilde-mu}), we find
\begin{align}
 \label{eq:Atilde-nu-result}
\tilde A_{\nu'6}=&-\frac{\Theta_{\infty6}}2\sigma_3J_{\nu\nu'}^{\pm\pm}+
J_{\nu'\nu}^{\pm\pm}\left(A_{\nu'6}-\frac12\right)-J_{\nu\nu'}^{\pm\pm}\left(A_{\nu6}-\frac12\right)-\\
&\left(\frac{\nu'-\mu}{\nu-\mu}J_{\nu'\nu}^{\pm\pm}+J_{\nu\nu'}^{\pm\pm}\right)A_{\mu6}J_{\nu\nu'}^{\pm\pm},
\nonumber
\end{align}
where $\tilde A_{\nu6}$ is given by the same Equation~(\ref{eq:Atilde-nu-result}) but with the permutation
$\nu'\leftrightarrow\nu$;  $\tilde A_{\mu6}$ can then be found from Equation~(\ref{eq:Atilde-infty}).

Consider, now, Transformation~(\ref{eq:L-dressing}) for $\nu'=\infty$: we find that
\begin{equation}
 \label{eq:J-nu-infty}
J_{\nu\infty}^{\pm+}=\frac1{a_\nu^\pm}\left(\begin{array}{cc}1&0\\0&-\Psi_{1\infty6}^{21}\end{array}\right)
\left(\begin{array}{cc}a_\nu^\pm&b_\nu^\pm\\a_\nu^\pm&b_\nu^\pm\end{array}\right),\quad
J_{\nu\infty}^{\pm-}=
\frac1{b_\nu^\pm}\left(\begin{array}{cc}-\Psi_{1\infty6}^{12}&0\\0&1\end{array}\right)
\left(\begin{array}{cc}a_\nu^\pm&b_\nu^\pm\\a_\nu^\pm&b_\nu^\pm\end{array}\right),
\end{equation}
where $a_\nu^\pm$ and $b_\nu^\pm$ are defined by Equations~(\ref{eq:Rnu6a-b}) and (\ref{eq:A-and-R}),
$\Psi_{1\infty6}^{ij}$ are the matrix elements of the first coefficient of the
expansion (\ref{eq:psi6-infty-normalization}). These matrix elements can be calculated via the matrix
elements of $A_{\nu6}$, since (recall $\Theta_{\infty6}\notin\mathbb{Z}$):
$$
-\Psi_{1\infty6}+\frac{\Theta_{\infty6}}2\big[\sigma_3,\Psi_{1\infty6}\big]=A_{16}+t_6A_{t6}.
$$
Thus, Transformation~(\ref{eq:L-dressing}) with $\nu'=\infty$ exists iff
\begin{equation}
 \label{eq:L-exiatence}
{\cal L}_{\nu\infty}^{\pm+}:\;\;a_\nu^\pm\neq0,\qquad
{\cal L}_{\nu\infty}^{\pm-}:\;\;b_\nu^\pm\neq0.
\end{equation}
The violation of conditions~(\ref{eq:L-exiatence}) can be discussed in the same manner as the violation of
Condition~(\ref{eq:Delta-ne0}) (see the paragraph between Equations~(\ref{eq:J-projectors})
and (\ref{eq:Atilde-nu})). The only difference is that Equation~(\ref{eq:first-order-y6}) now takes the
form
$$
\left(\frac{dy_6}{dt_6}\right)^2=R(y_6,t_6).
$$
Hereafter we assume that Conditions (\ref{eq:L-exiatence}) are valid (as well as the previously
assumed (\ref{eq:Delta-ne0})).

Consider now the action of ${\cal L}_{\nu\infty}^{\pm\pm}$ on $A_{\mu6}$, namely:
\begin{equation}
 \label{eq:A-tilde-mu-infty}
\tilde A_{mu6}=\left(\sigma_\infty^\pm+\frac1{\mu-\nu}J_{\nu\infty}^{\pm\pm}\right)A_{\mu6}
\left(\sigma_\infty^\pm+\frac1{\mu-\nu}J_{\nu\infty}^{\pm\pm}\right)^{-1},\qquad\mu\neq\nu.
\end{equation}
The action of ${\cal L}_{\nu\infty}^{\pm\pm}$ on $A_{\nu6}$ can be obtained by substituting
Equation~(\ref{eq:A-tilde-mu-infty}) for
$\tilde A_{mu6}$ and (with $\mu\leftrightarrow\mu'$) for $\tilde A_{mu'6}$, where $\mu'$ is
defined from the condition that $\{\nu,\mu,\mu'\}$ is a permutation of $0,1,t$, into the equation
\begin{equation}
 \label{eq:Atilde-infty2}
\tilde A_{\nu6}+\tilde A_{\mu6}+\tilde A_{\mu'6}=-\frac{\tilde\Theta_{\infty6}^\pm}2\sigma_3=
-\frac{\Theta_{\infty6}\pm1}2\sigma_3,\qquad \mu\neq\mu'\neq\nu\neq\mu.
\end{equation}

To summarize, let us fix $t_6^0\in\mathbb{C}\setminus\{0,1\}$ and $A_{\nu6}^0$, $\Theta_{\nu6}^0$ as
demanded by the last conditions in (\ref{eq:niormalization-for-P6}) and (\ref{eq:Rnu6}). The general
Schlesinger (isomonodromy) deformations (GSD) for Equation~(\ref{eq:linear-4-regular}) are the matrices
$A_{\nu6}$ (or their matrix elements) which depend on the continuous variable $t_6$ and the discrete
variables $\Theta_{\nu6}$: $\Theta_{\nu6}-\Theta_{\nu6}^0\in\mathbb{Z},\;
\sum_\nu(\Theta_{\nu6}-\Theta_{\nu6}^0)=0({\rm mod}2)$. The continuous deformations of $A_{\nu6}$ are
governed by (\ref{eq:schlesinger-p6}) and the discrete deformations by Equations~(\ref{eq:Atilde-mu})--
(\ref{eq:Atilde-nu-result}) and (\ref{eq:A-tilde-mu-infty}), (\ref{eq:Atilde-infty2}). The initial
condition is stated as $A_{\nu6}(t_6,\Theta_{06}^0,\Theta_{16}^0,\Theta_{t6}^0,\Theta_{\infty6}^0)=A_{\nu6}^0$.
The continuous and discrete deformations are commuting so that GSD are correctly defined. Any GSD can be
uniquely characterized by a point on
${\cal M}_6(\Theta_{06}^0,\Theta_{16}^0,\Theta_{t6}^0,\Theta_{\infty6}^0)$. {\it The direct monodromy
problem} is: construct ${\cal M}_6(\Theta_{06}^0,\Theta_{16}^0,\Theta_{t6}^0,\Theta_{\infty6}^0)$ for
the given GSD. {\it The inverse monodromy problem} is: construct the GSD for given
${\cal M}_6(\Theta_{06}^0,\Theta_{16}^0,\Theta_{t6}^0,\Theta_{\infty6}^0)$ \cite{B}.

Consider the $2\times2$ matrix linear ODE related with $P_5$:
\begin{equation}
 \label{eq:linear-5}
\frac{d\Psi_5}{d\lambda}=\left(\frac{t_5}2\sigma_3+\frac{A_{05}}\lambda+\frac{A_{15}}{\lambda-1}\right)\Psi_5.
\end{equation}
This equation possess two regular singular points, at $\lambda=0$ and $1$, and an irregular one at
the point of $\infty$.
We require the following conditions:
\begin{equation}
 \label{eq:diag-A5}
{\rm diag}(A_{05}+A_{15})=-\frac{\Theta_{\infty5}}2\sigma_3,\qquad\Theta_{\infty5}\in\mathbb{C};
\end{equation}
and that there exists $R_{\nu5}\in{\rm SL}(2,\mathbb{C})$ such that
\begin{equation}
 \label{eq:R5-def}
R_{\nu5}^{-1}A_{\nu5}R_{\nu5}=\frac{\Theta_{\nu5}}2\sigma_3,\qquad
\Theta_{\nu5}\in\mathbb{C}\setminus\mathbb{Z},\quad\nu=0,1.
\end{equation}
To define the monodromy data let us define the canonical solutions $\Psi_5^k$ of
Equation~(\ref{eq:linear-5}) by setting their asymptotics at the infinity point as
\begin{gather}
 \label{eq:Psi5-asymptotics-infty}
\Psi_5^k=\left(I+\sum\limits_{m=1}^\infty\Psi_{m\infty5}\lambda^{-m}\right)
\exp\!\left(\!\!\Big(\frac{\lambda t_5}2-\frac{\Theta_{\infty5}}2\ln\lambda\Big)\!\sigma_3\!\right),\\
\lambda\to\infty,\qquad
-\frac32\pi+\pi k<\arg(\lambda t_5)<\frac{\pi}2+\pi k,\quad k\in\mathbb{Z}.\nonumber
\end{gather}
The single-valued function $\Psi_5^k$ can be defined in the domain
$\mathbb{C}\setminus\big([0,1]\cup[0,\infty e^{\frac{i\pi}2+\pi ki}\big)$,
where we choose the main branch
of $\ln\lambda$: ${\rm Im}\ln\lambda={\rm arg}\lambda$. Using Asymptotics
(\ref{eq:Psi5-asymptotics-infty}), one proves that
\begin{equation}
 \label{eq:infinity-circle}
\Psi_5^{k+2}\big(\lambda\, e^{2\pi i}\big)=\Psi_5^k(\lambda)e^{-\pi i\Theta_{\infty5}\sigma_3}.
\end{equation}
Now we define the Stokes matrices $S_k$ as
\begin{equation}
 \label{eq:Stokes-definition}
\Psi_5^{k+1}(\lambda)=\Psi_5^k(\lambda)S_k.
\end{equation}
Definitions (\ref{eq:Psi5-asymptotics-infty}) and (\ref{eq:Stokes-definition}) yield
$$
S_{2l}=\left(\begin{array}{cc}1&0\\s_{2l}&1\end{array}\right),\qquad
S_{2l+1}=\left(\begin{array}{cc}1&s_{2l+1}\\0&1\end{array}\right),\qquad l\in\mathbb{Z},
$$
where $s_k$ are called the Stokes multipliers. Using Equations~(\ref{eq:infinity-circle}) and
(\ref{eq:Stokes-definition}) one finds
\begin{equation}
 \label{eq:stokes-shift-2}
 S_{k+2}=e^{\pi i\Theta_{\infty5}\sigma_3}S_ke^{-\pi i\Theta_{\infty5}\sigma_3}.
\end{equation}
The monodromy matrix at the point of infinity, $M_{k\infty5}$, for the $\Psi_5^k$-function is given by the
equation
\begin{equation}
 \label{eq:monodromy-infty}
\Psi_5^k\big(\lambda\,e^{-2\pi i}\big)=\Psi_5^k(\lambda)M_{k\infty5}.
\end{equation}
Comparing Equations~(\ref{eq:infinity-circle}) and (\ref{eq:Stokes-definition}), one arrives at
\begin{equation}
 \label{eq:M-infty-stokes-k}
M_{k\infty5}=S_kS_{k+1}e^{\pi i\Theta_{\infty5}\sigma_3}.
\end{equation}
In the following we set
$$
M_{\infty5}\equiv M_{0\infty5}.
$$
All the others $M_{k\infty5}$ can be expressed in terms of $M_{\infty5}$ via the recurrence formula
$$
M_{k+1\infty5}=S_k^{-1}M_{k\infty5}S_k
$$
and Equation~(\ref{eq:M-infty-stokes-k}). The monodromy matrices $M_{\nu5}$, $\nu=0,1$, at the regular
singularities $\lambda=\nu$ of Equation~(\ref{eq:linear-5}) are defined with the help of the paths
given in Figure~7 by the same formulae (\ref{eq:psi-expansion-at-nu}), (\ref{eq:condition-det-1}),
(\ref{eq:monodromy-local}), and (\ref{eq: trace-M}) with the subscript $6$ changed to $5$.
\begin{figure}[ht]
\begin{center}
\begin{picture}(160,100)
\put(30,80){\circle*{3}}
\put(28,66){$0$}
\put(80,80){\circle*{3}}
\put(72,66){$1$}
\put(150,80){\circle*{3}}
\put(136,70){$\infty$}
\put(20,20){\circle*{3}}
\put(5,15){$\lambda_0$}
\qbezier(20,20)(60,100)(30,100)
\qbezier(20,20)(10,100)(30,100)
\qbezier(20,20)(118,80)(90,98)
\qbezier(20,20)(72,106)(90,98)
\qbezier(20,20)(188,68)(164,90)
\qbezier(20,20)(160,110)(164,90)
\put(80,62){\vector(1,1){10}}
\put(142,62){\vector(2,1){10}}
\put(37,58){\vector(1,3){5}}
\put(60,00){Figure 7}
\end{picture}
\begin{picture}(40,90)
\put(70,30){\circle*{3}}
\put(68,36){$0$}
\put(120,30){\circle*{3}}
\put(112,36){$1$}
\put(190,30){\circle*{3}}
\put(176,36){$\infty$}
\put(60,90){\circle*{3}}
\put(45,85){$\tilde\lambda_0$}
\qbezier(60,90)(90,25)(70,20)
\qbezier(60,90)(45,25)(70,20)
\qbezier(60,90)(155,28)(130,20)
\qbezier(60,90)(94,22)(130,20)
\qbezier(60,90)(188,08)(214,20)
\qbezier(60,90)(247,30)(214,20)
\put(124,42){\vector(-1,1){10}}
\put(189,43){\vector(-2,1){10}}
\put(77,41){\vector(-1,4){3}}
\end{picture}
\begin{picture}(160,100)
\put(60,00){Figure 8}
\end{picture}
\end{center}
\end{figure}

The cyclic relation reads
\begin{equation}
 \label{eq:P5-cyclic-M}
M_{05}M_{15}M_{\infty5}=I.
\end{equation}
The monodromy group $MG_5$ is a subgroup of ${\rm SL}(2,\mathbb{C})$ generated by the matrices $M_{\nu5}$.
So, under the monodromy group we mean a particular monodromy representation of the fundamental group
$\pi(\lambda_0,\bar{\mathbb{C}}\setminus\{0,1,\infty\})$. Together with the group $MG_5$ defined above
we will use another one, $\widetilde{MG}_5$, which is generated by the monodromy matrices
$\widetilde M_{\nu6}\equiv\widetilde M_{0\nu6}$ obtained by the analytic continuation of the same
canonical solution $\Psi_5^0(\lambda)$ but along the paths presented in Figure~8.

These matrices obey the following cyclic relation:
\begin{equation}
 \label{eq:P5-cyclic-M-tilde}
\widetilde M_{\infty5}\widetilde M_{15}\widetilde M_{05}=I
\end{equation}
The relation between both monodromy groups can be obtained by comparing representations
of the fundamental groups in Figures 7 and 8 together with our way of defining a singlevalued
branch of the function $\Psi_5^0$ explained in the paragraph below the
asymptotics at the point of infinity (\ref{eq:Psi5-asymptotics-infty}):
$$
\widetilde M_{15}=M_{05}M_{15}M_{05}^{-1},\quad\widetilde M_{05}=M_{05},\quad
\widetilde M_{\infty5}=M_{\infty5}.
$$
{\it The sets of monodromy data} ${\cal M}_5$ (or $\widetilde{\cal M}_5$) are ordered sets of the matrix elements
$\{M_{\nu5}\}_{\nu=0,1,\infty}$ (or $\{\widetilde{M}_{\nu5}\}_{\nu=0,1,\infty}$) completed with three
complex parameters $\Theta_{\nu5}$ ($\nu=0,1,\infty$) satisfying the following equations
\begin{equation}
 \label{eq:M5-traces}
\begin{aligned}
&\nu=0,1:\;\;{\rm tr}\,M_{\nu5}=2\cos(\pi\Theta_{\nu5}),\quad\Theta_{\nu5}\in\mathbb{C}\setminus\mathbb{Z},\\
&\nu=\infty:\;\;{\rm tr}\,M_{\infty5}=2\cos(\pi\Theta_{\infty5})+e^{-\pi i\Theta_{\infty5}}s_0s_1,\quad
\Theta_{\infty5}\in\mathbb{C}.
\end{aligned}
\end{equation}
{\it The manifold of monodromy data} ${\cal M}_5(\Theta_{05},\Theta_{15},\Theta_{\infty5})$ is an
algebraic variety defined by Equations~(\ref{eq:P5-cyclic-M}), (\ref{eq:M5-traces}) and
$$
\det M_{\nu5}=1,\qquad\nu=0,1,\infty
$$
in $\mathbb{C}^{12}$ (we identify an ordered set of matrix elements $\{M_{\nu5}\}_{\nu=0,1,\infty}$ as
a point in $\mathbb{C}^{12}$). In an analogous way we define
$\widetilde{\cal M}_5(\Theta_{05},\Theta_{15},\Theta_{\infty5})$.
It is easy to see that the complex dimension of ${\cal M}_5(\Theta_{05},\Theta_{15},\Theta_{\infty5})$
(respectively, $\widetilde{\cal M}_5(\Theta_{05},\Theta_{15},\Theta_{\infty5})$) equals $3$.

Consider, following \cite{JM}, the parametrization of $A_{\nu5}$:
\begin{gather}
 \label{eq:A05-parametrization}
A_{05}=\left(\!\!\begin{array}{cc}
z_5+\Theta_{05}/2&-u_5(z_5+\Theta_{05})\\
u_5^{-1}z_5&-z_5-\Theta_{05}/2
\end{array}\!\!\right),\\
 \label{eq:A15-parametrization}
A_{15}=\left(\!\!\begin{array}{cc}
-z_5-(\Theta_{05}+\Theta_{\infty5})/2&u_5y_5(z_5+(\Theta_{05}-\Theta_{15}+\Theta_{\infty5})/2)\\
-\frac1{u_5y_5}(z_5+(\Theta_{05}+\Theta_{15}+\Theta_{\infty5})/2)&z_5+(\Theta_{05}+\Theta_{\infty5})/2
\end{array}\!\!\right).
\end{gather}
We see that for fixed $t_5$ and formal monodromies $\Theta_{\nu5}$ ($\nu=0,1,\infty$) the number of
parameters ($u_5$, $z_5$, $y_5$) in Equation~(\ref{eq:linear-5}) is $3$: it exactly coincides with
${\rm dim}\,{\cal M}_5(\Theta_{05},\Theta_{15},\Theta_{\infty5})$.

The hard isomonodromy condition: $\partial_{t_5}\Theta_{\nu5}=0$ and $\partial_{t_5}M_{\nu5}=0$, for
$\nu=0,1,\infty$ implies an additional ODE for the function $\Psi_5$ with respect to $t_5$:
\begin{equation}
 \label{eqP5-linear-t5}
\frac{d\Psi_5}{dt_5}=\left(\frac\lambda2\sigma_3+
\frac1{t_5}\left(\frac{\Theta_{\infty5}}2\sigma_3+A_{05}+A_{15}\right)\!\!\right)\Psi_5.
\end{equation}
The compatibility condition of Equations~(\ref{eq:linear-5}) and (\ref{eqP5-linear-t5}) implies that
the matrices $A_{\nu5}$ and hence the parameters: $u_5$, $z_5$, and $y_5$, are functions of $t_5$.
These functions are governed by the system of isomonodromy deformations,
\begin{equation}
 \label{eq:system-idm5}
\frac{dA_{\nu5}}{dt_5}=\left[\frac1{t_5}\left(\frac{\Theta_{\infty5}}2\sigma_3+A_{05}+A_{15}\right)+
\frac\nu2\sigma_3,\,A_{\nu5}\right],\qquad\nu=0,1.
\end{equation}
In terms of the parameters $u_5$, $z_5$, and $y_5$, this system is given in \cite{JM}: eliminating the
function $z_5$ from these equations one finds that $y_5(t_5)$ solves the fifth Painlev\'e equation:
\begin{equation}
 \label{eq:P5}
 \begin{gathered}
\frac{d^2y_5}{dt_5^2}=\!\left(\frac1{2y_5}+\frac1{y_5-1}\right)\!\!\left(\frac{dy_5}{dt_5}\right)^2-
\frac1{t_5}\frac{dy_5}{dt_5}+\left(\frac{y_5-1}{t_5}\right)^2\!\!
\left(\alpha_5y_5+\frac{\beta_5}{y_5}\right)\!+\\
\gamma_5\frac{y_5}{t_5}+\delta_5\frac{y_5(y_5+1)}{y_5-1},
\end{gathered}
\end{equation}
\begin{equation}
 \label{eq:P5coeff}
\begin{gathered}
\alpha_5=\frac12\left(\frac{\Theta_{05}-\Theta_{15}+\Theta_{\infty5}}2\right)^2,\quad
\beta_5=-\frac12\left(\frac{\Theta_{05}-\Theta_{15}-\Theta_{\infty5}}2\right)^2,\\
\gamma_5=1-\Theta_{05}-\Theta_{15},\quad\delta_5=-\frac12.
\end{gathered}
\end{equation}
There is an ambiguity in the definition of the formal monodromies $\Theta_{\nu5}\to-\Theta_{\nu5}$,
$\nu=0,1$, as well as in the case for the function $\Psi_6$. The change $\Theta_{\nu5}\to-\Theta_{\nu5}$
leads to a repa\-ra\-met\-rization of the matrices $A_{\nu5}$. Contrary to the above case for $P_6$ this
reparametrization yields a nontrivial transformation of the solution $y_5$. We won't discuss it here.
It is important for us
that we can choose the signs of $\theta_{\nu5}$, $\nu=0,1$, arbitrarily, and then use the corresponding
parametrization of $A_{\nu5}$.

The $\tau$-function for the isomonodromy deformations~(\ref{eq:system-idm5}) is defined in \cite{JM}
as follows:
\begin{equation}
 \label{eq:tau5-via-Psi}
\frac{d}{dt_5}\ln\,\tau_5(t_5)=-\frac12{\rm tr}\,(\Psi_{1\infty5}\sigma_3).
\end{equation}
In terms of the matrix elements Equation~(\ref{eq:system-idm5}) reads:
\begin{equation}
 \label{eq:tau5-yz}
\begin{gathered}
\frac{d}{dt_5}\ln\,\tau_5(t_5)=(A_{15})_{11}+\frac1{t_5}(A_{05}+A_{15})_{21}(A_{05}+A_{15})_{12}=
-z_5-\frac{\Theta_{05}+\Theta_{\infty5}}2\\-\frac1{t_5}\left(z_5-\frac1{y_5}
\left(\!\!z_5+\frac{\Theta_{05}+\Theta_{15}+\Theta_{\infty5}}2\right)\!\!\right)\!\left(\!\!z_5+\Theta_{05}-
y_5\!\left(\!z_5+\frac{\Theta_{05}-\Theta_{15}+\Theta_{\infty5}}2\right)\!\!\right).
\end{gathered}
\end{equation}
We can also rewrite Equation~(\ref{eq:tau5-yz}) in the ``matrix'' form:
\begin{equation}
 \label{eq:tau5-matrix-form}
\frac{d}{dt_5}\ln\,\tau_5(t_5)=\frac1{t_5}\left(\left(\frac{\Theta_{05}}2\right)^2+
\left(\frac{\Theta_{15}}2\right)^2-\left(\frac{\Theta_{\infty5}}2\right)^2\right)+
{\rm tr}\,\Big(\frac1{t_5}A_{05}+\frac{\sigma_3}2\Big)A_{15}.
\end{equation}
In \cite{JM} the function $\sigma_5(t_5)$:
\begin{equation}
 \label{eq:sigma5-def}
\sigma_5(t_5)=\frac{\Theta_{05}+\Theta_{\infty5}}2t_5+t_5\frac{d}{dt_5}\ln\,\tau_5(t_5).
\end{equation}
This function satisfies a second-order ODE which is quadratic with respect to $\sigma_5^{''}$ (see \cite{JM}).
Differentiating (\ref{eq:sigma5-def}) and using (\ref{eq:tau5-matrix-form}), (\ref{eq:system-idm5}), and
parametrization~(\ref{eq:A05-parametrization}), (\ref{eq:A15-parametrization}), we find that
\begin{equation}
 \label{eq:sigma'=-z5}
 \frac{d\sigma_5}{dt_5}=-z_5.
\end{equation}
In the corresponding formula (C.44) of \cite{JM} there is a misprint in the sign: this sign is important
for us to establish the differentiable character of our asymptotic expansion for the function
$\hat\sigma_6(t_6)$ (see Subsection~{\bf II.4} of Section~\ref{sec:3}).
\section{Formal Limit Transitions $P_6\to P_5$}
 \label{sec:3}
As mentioned in the Introduction we consider here two different limits. Our scheme for the derivation of
these limits consists of the following steps:
\begin{enumerate}
\item
A formal limit passage of Equation~(\ref{eq:linear-4-regular}) to Equation~(\ref{eq:linear-5});
\item
Finding conditions which guarantee that simultaneously with the limit passage in item~1 we have
the limit passage of Equation~(\ref{eq:t-6}) to Equation~(\ref{eqP5-linear-t5}). Actually these
conditions are additional to those found in the first step. One proves that if the asymptotic
expansions found in the first step are differentiable with respect to $t_5$, then it is possible
to define asymptotic expansions for the functions $s_{\nu6}$, $\nu=1,t$,
(see Equation~(\ref{eq:A-and-R})) to satisfy additional conditions appearing at this step.
To prove that the asymptotics we found are really differentiable, one has to use the systems of
isomonodromy deformations~(\ref{eq:schlesinger-p6}) (or, in terms of the matrix elements,
Equations (C.51), (C.52), and (C.55) of \cite{JM}) and Equation~(\ref{eq:system-idm5}) (or (C.40) of
\cite{JM}). We leave this proof to the reader and write down only the formulae for $s_{\nu6}$.
\item
Presentation of the formal limit passage in terms of the matrix elements of $A_{\nu6}$ and $A_{\nu5}$.
\item
Presentation of the limit as (formal) asymptotics for the $P_6$- and $\tau_6$-functions.
\end{enumerate}

To eliminate possible confusion, let us agree to supply the parameters $\lambda$ from Equations
(\ref{eq:linear-4-regular}) and (\ref{eq:linear-5}) with the subscripts 6 or 5, respectively.

{\bf I.1.} The first limit passage:
\begin{gather}
 \label{eq:lim1-t-lambda}
\varepsilon\to+0,\quad
\Theta_{16}=-\frac1\varepsilon,\quad
t_6=\varepsilon t_5={\cal O}(\varepsilon),\quad
\lambda_6=\varepsilon t_5\lambda_5=o(\varepsilon),\\
 \label{eq:lim1-Psi}
\underset{\varepsilon\to+0}\lim\,R_{16}^{-1}\Psi_6(\lambda_6,t_6)=\Psi_5(\lambda_5,t_5),\\
 \label{eq:lim1-A}
R_{16}^{-1}A_{06}R_{16}=A_{05}+{\cal O}(\varepsilon),\qquad
R_{16}^{-1}A_{t6}R_{16}=A_{15}+{\cal O}(\varepsilon).
\end{gather}

{\bf I.2.}
\begin{equation}
 \label{eq:lim1-t-R}
 -R_{16}^{-1}\frac{d}{dt_5}R_{16}=\frac{\Theta_{\infty5}}{2t_5}\sigma_3+{\cal O}(\varepsilon),
\end{equation}
where $\Theta_{\infty5}$ appears as a parameter of the limit passage, i.e., an arbitrary complex number,
its notation as one of the formal monodromies related with the fact that in derivation of
Equation~(\ref{eq:lim1-t-R}) we took into account relation~(\ref{eq:diag-A5}).

{\bf I.3.} In terms of the matrix elements $A_{\nu6}$ and  $A_{\nu5}$ we can write the limit {\bf I} in two
possible ways. So for the matrix elements we get the two different formal limits {\bf I.3.a.} and
{\bf I.3.b.}. Both limits are related via the transformation generated by the
Reflection~(\ref{eq:reflection-infinity}); however, the leading terms of one limit do not completely
define the leading terms of the other: some further terms of the expansions are needed. So the formulae
given below for the limit {\bf I.3.a.} do not completely define asymptotics {\bf I.3.b.}

{\bf I.3.a.}
\begin{equation}
 \label{eq:13a-theta}
\Theta_{16}=-\frac1\varepsilon,\quad
\Theta_{\infty6}+\Theta_{16}=\Theta_{\infty5},\quad
\Theta_{06}=\Theta_{05},\quad
\Theta_{t6}=-\Theta_{15}.
\end{equation}
In Section~\ref{sec:2}, in the paragraph following Equation~(\ref{eq:P5coeff}), we explained that the signs
in the last two equations of (\ref{eq:13a-theta}) can be taken arbitrarily; however, from their choice,
the following formulae are strongly depended.
\begin{gather}
\label{eq:13a-z06-zt6}
z_{06}=z_5+{\cal O}(\varepsilon),\quad
z_{t6}=-z_5-\frac{\Theta_{05}-\Theta_{15}+\Theta_{\infty5}}2+{\cal O}(\varepsilon),\\
z_{16}=-\varepsilon\left(z_5-
\frac1{y_5}\left(z_5+\frac{\Theta_{05}+\Theta_{15}+\Theta_{\infty5}}2\right)\!\!\right)\times
\nonumber\\
\left(z_5+\Theta_{05}-y_5\left(z_5+\frac{\Theta_{05}-\Theta_{15}+\Theta_{\infty5}}2\right)\!\!\right)+
{\cal O}(\varepsilon^2),
\label{eq:13a-z16}\\
u_{06}s_{16}^2=\varepsilon^2u_5\frac{z_5+\Theta_{05}}{z_5}+{\cal O}(\varepsilon^3),\qquad
u_{t6}s_{16}^2=\varepsilon^2y_5u_5+{\cal O}(\varepsilon^3),
\label{eq:13a-u06-ut6}\\
u_{16}s_{16}^2=\frac{\varepsilon u_5}{z_5-\frac1{y_5}
\left(z_5+\frac{\Theta_{05}+\Theta_{15}+\Theta_{\infty5}}2\right)}+{\cal O}(\varepsilon^2).
\label{eq:13a-u16}
\end{gather}
As follows from Equation~(\ref{eq:lim1-t-R}) the function $s_{16}$ must satisfy the equation
\begin{equation}
 \label{eq:13a-s16}
 \frac{d}{dt_5}\ln\,s_{16}=\frac{\Theta_{\infty5}}{2t_5}+{\cal O}(\varepsilon).
\end{equation}

{\bf I.3.b.}
\begin{gather}
 \label{eq:13b-Theta}
\Theta_{16}=-\frac1\varepsilon,\quad
\Theta_{16}-\Theta_{\infty6}=\Theta_{\infty5},\quad
\Theta_{06}=\Theta_{05},\quad
\Theta_{t6}=-\Theta_{15},\\
 \label{eq:13b-z06-zt6}
z_{06}=-z_5-\Theta_{05}+{\cal O}(\varepsilon),\quad
z_{t6}=z_5+\frac{\Theta_{05}+\Theta_{15}+\Theta_{\infty5}}2+{\cal O}(\varepsilon),\\
 \label{eq:13b-z16}
z_{16}=\frac1\varepsilon+\varepsilon\left(z_5-
\frac1{y_5}\left(z_5+\frac{\Theta_{05}+\Theta_{15}+\Theta_{\infty5}}2\right)\!\!\right)\times
\nonumber\\
\left(z_5+\Theta_{05}-y_5\left(z_5+\frac{\Theta_{05}-\Theta_{15}+\Theta_{\infty5}}2\right)\!\!\right)+
{\cal O}(\varepsilon^2),\\
 \label{eq:13b-u16/ut6}
\frac{u_{16}}{u_{t6}}=\varepsilon y_5\left(z_5-
\frac1{y_5}\left(z_5+\frac{\Theta_{05}+\Theta_{15}+\Theta_{\infty5}}2\right)\!\!\right)+
{\cal O}(\varepsilon^2),\\
 \label{eq:13b-u16/u06}
\frac{u_{16}}{u_{06}}=\varepsilon\frac{z_5+\Theta_{05}}{z_5}\left(z_5-
\frac1{y_5}\left(z_5+\frac{\Theta_{05}+\Theta_{15}+\Theta_{\infty5}}2\right)\!\!\right)+
{\cal O}(\varepsilon^2),\\
 \label{eq:13b-u16s16}
u_{16}s_{16}^2=-\frac{\varepsilon u_5}{z_5-\frac1{y_5}
\left(z_5+\frac{\Theta_{05}+\Theta_{15}+\Theta_{\infty5}}2\right)}+{\cal O}(\varepsilon^2),\\
 \label{eq:13b-s16}
\frac{d}{dt_5}\ln\,s_{16}=\frac{d}{dt_5}\ln\frac{u_5}{z_5-\frac1{y_5}
\left(z_5+\frac{\Theta_{05}+\Theta_{15}+\Theta_{\infty5}}2\right)}
-\frac{\Theta_{\infty5}}{2t_5}+{\cal O}(\varepsilon).
\end{gather}

{\bf I.4.} Substituting Asymptotics (\ref{eq:lim1-A}) and (\ref{eq:lim1-t-lambda}) into
Equation~(\ref{eq:tau-6}) and using definition~(\ref{eq:tau5-matrix-form}) one finds that
\begin{equation}
 \label{eq:tau-6-asympt}
\frac{d}{dt_5}\ln\tau_6(t_6)=
\frac1{t_5}\left(\left(\frac{\Theta_{\infty5}}2\right)^2-\left(\frac{\Theta_{05}}2\right)^2-
\left(\frac{\Theta_{15}}2\right)^2\right)+\frac{d}{dt_5}\ln\tau_5(t_5)+{\cal O}(\varepsilon).
\end{equation}
Substituting Asymptotics (\ref{eq:13a-z06-zt6}) and (\ref{eq:13a-z16}) into Equation~(\ref{eq:y6})
we find for the limit {\bf I.3.a.}
\begin{equation}
 \label{eq:13a-y6-asympt}
y_6(t_6)=\frac{\varepsilon t_5}{1+y_5\left(1-
\frac{\Theta_{05}+\Theta_{15}-\Theta_{\infty5}}{2(z_5+\Theta_{\infty5})}\right)}+
{\cal O}(\varepsilon^2).
\end{equation}
In case {\bf I.3.b.} Equations~(\ref{eq:13b-z06-zt6}), (\ref{eq:13b-u16/ut6}) and (\ref{eq:13b-u16/u06})
yield
\begin{equation}
 \label{eq:13b-y6-asympt}
y_6(t_6)=\frac{\varepsilon t_5}{1+\frac1{y_5}\left(1+
\frac{\Theta_{05}+\Theta_{15}+\Theta_{\infty5}}{2z_5}\right)}+{\cal O}(\varepsilon^2).
\end{equation}
Now we consider the second formal limit.

{\bf II.1.}
\begin{gather}
 \label{eq:lim2-t-lambda}
\varepsilon\to+0,\quad
\Theta_{t6}=-\frac1\varepsilon,\quad
t_6=\varepsilon t_5={\cal O}(\varepsilon),\quad
\lambda_6=\frac1{\lambda_5},\quad
\frac\varepsilon{\lambda_6^3}=o(1),\\
\label{eq:lim2-Psi}
\underset{\varepsilon\to+0}\lim R_{t6}^{-1}\Psi_6(\lambda_6,t_6)=\Psi_5(\lambda_5,t_5),\\
 \label{eq:lim2-A05-A15}
R_{t6}^{-1}\frac{\Theta_{\infty6}}2\sigma_3R_{t6}=A_{05}+{\cal O}(\varepsilon),\qquad
R_{t6}^{-1}A_{16}R_{t6}=A_{15}+{\cal O}(\varepsilon)
\end{gather}

{\bf II.2.}
\begin{equation}
 \label{eq:lim2-Rt6-derivative}
-R_{t6}^{-1}\frac{d}{dt_5}R_{t6}=\frac1{t_5}\left(\frac{\Theta_{\infty5}}2\sigma_3+A_{05}+A_{15}\right)+
{\cal O}(\varepsilon),
\end{equation}
where $\Theta_{\infty5}$ is, as for the first limit, a parameter of the limit passage satisfying
Equation~(\ref{eq:diag-A5}).

{\bf II.3.}
\begin{equation}
 \label{eq:lim2-theta}
\Theta_{t6}=-\frac1\varepsilon,\quad
\Theta_{t6}+\Theta_{06}=\Theta_{\infty5},\quad
\Theta_{\infty6}=\Theta_{05},\quad
\Theta_{16}=\Theta_{15}.
\end{equation}
Together with the case $\Theta_{t6}+\Theta_{06}=\Theta_{\infty5}$ we can consider another one:
$\Theta_{t6}-\Theta_{06}=\Theta_{\infty5}$. As explained in Section~\ref{sec:2}, transformation
$\Theta_{06}\to-\Theta_{06}$ simply means that the reparametrization of the matrix $A_{06}$, which does no
effect on both the $\tau_6$- and the $P_6$-functions. Thus, contrary to the first limit, it is not
worthwhile to consider separately these possibilities:
\begin{gather}
 \label{eq:II-zt6-z06}
z_{t6}=-\frac{z_5}{\varepsilon\Theta_{05}}+{\cal O}(1),\qquad
z_{06}=\frac{z_5}{\varepsilon\Theta_{05}}+{\cal O}(1),\\
 \label{eq:II-z16}
z_{16}+\Theta_{16}=\left(\frac{z_5(1-y_5)}{\Theta_{05}}+1\right)\!\!
\left(\Theta_{15}+\frac{1-y_5}{y_5}\left(z_5+
\frac{\Theta_{05}+\Theta_{15}+\Theta_{\infty5}}2\right)\!\!\right)+{\cal O}(\varepsilon),\\
 \label{eq:II-ut6-u06}
u_{t6}=-\varepsilon^2\frac{\Theta_{05}u_5}{z_5s_{t6}^2}+{\cal O}(\varepsilon^3),\qquad
u_{06}=u_{t6}+{\cal O}(\varepsilon^3),\\
 \label{eq:II-u16}
u_{16}=u_{t6}\frac{\Theta_{15}+\left(\frac1{y_5}-1\right)\left(z_5+
\frac{\Theta_{05}+\Theta_{15}+\Theta_{\infty5}}2\right)}
{\Theta_{15}+\left(\frac1{y_5}-1\right)\left(z_5+
\frac{\Theta_{05}+\Theta_{15}+\Theta_{\infty5}}2\right)\left(1+\frac{\Theta_{05}}{z_5(1-y_5)}\right)}+
{\cal O}(\varepsilon^3),\\
 \label{eq:II-st6}
\frac{d}{dt_5}\ln\,s_{t6}=\frac{z_5}{\Theta_{05}}\frac{d}{dt_5}\ln\,u_{t6}+{\cal O}(\varepsilon),\\
 \label{eq:II-ut6}
t_5\frac{d}{dt_5}\ln\,u_{t6}=-\frac{\Theta_{05}}{z_5}\left(z_5-
\frac1{y_5}\left(z_5+\frac{\Theta_{05}+\Theta_{15}+\Theta_{\infty5}}2\right)\!\!\right)
+{\cal O}(\varepsilon),\\
 \label{eq:II-ut6zt6st6}
t_5\frac{d}{dt_5}\ln(u_{t6}z_{t6}s_{t6})=-\left(z_5+\Theta_{05}-
y_5\left(z_5+\frac{\Theta_{05}-\Theta_{15}+\Theta_{\infty5}}2\right)\!\!\right)
+{\cal O}(\varepsilon).
\end{gather}
Equations (\ref{eq:II-ut6}) and (\ref{eq:II-ut6zt6st6}) play a more important role in the proof of the
consistency of the second limit with the isomonodromy condition~(\ref{eq:lim2-Rt6-derivative}) than
the analogous formulae for the first limit, that is why they are written here explicitly. Note that
$\Theta_{05}\neq0$ due to the third equation in (\ref{eq:lim2-theta}) and the condition
$\Theta_{\infty6}\neq0$ which is imposed in Section~\ref{sec:2}.

{\bf II.4.} Substituting Equations~(\ref{eq:lim2-t-lambda}) and (\ref{eq:lim2-A05-A15}) into
(\ref{eq:tau-6}) and using normalization~(\ref{eq:niormalization-for-P6}) we find the following
asymptotics for the $\tau_6$-function:
\begin{equation}
 \label{eq:II-tau6-asympt-rough}
\frac{d}{dt_5}\ln\,\tau_6(t_6)=-\frac1{2\varepsilon^2t_5}-\frac{\Theta_{\infty5}}{2\varepsilon t_5}+
{\cal O}(1).
\end{equation}
Thus, the most interesting term of the asymptotic expansion (\ref{eq:II-tau6-asympt-rough}), that is
${\cal O}(1)$, cannot be obtained directly from our result: it requires more precise expansions in
(\ref{eq:II-zt6-z06}) (up to ${\cal O}(\varepsilon)$). Nevertheless, there is the following
trick to overcome this difficulty. Using Equations (\ref{eq:tau-6}), (\ref{eq:sigma-6}),
(\ref{eq:schlesinger-p6}), and (\ref{eq:niormalization-for-P6}), we, following \cite{JM}, find
\begin{equation}
 \label{eq:II-sigma-hat-derivative}
\frac{d\hat\sigma_6(t_6)}{dt_6}=-{\rm tr}\left(\frac{\Theta_{\infty6}}2\sigma_3A_{t6}\right)-
{\rm tr}\,A_{t6}^2.
\end{equation}
Now Equations~(\ref{eq:lim2-A05-A15}) yield
\begin{equation}
 \label{eq:II-sigma-hat-deriv-A05-asympt}
t_6\frac{d\hat\sigma_6(t_6)}{dt_6}=\frac{t_5}2{\rm tr}\,(A_{05}\sigma_3)-\frac{t_5}{2\varepsilon}+
{\cal O}(\varepsilon).
\end{equation}
Applying now (\ref{eq:diag-A5}) one obtains
\begin{equation}
 \label{eq:II-sigma-hat-deriv-A15-asympt}
t_6\frac{d\hat\sigma_6(t_6)}{dt_6}=-\frac{t_5}2\left(\frac1{\varepsilon}+\Theta_{\infty5}+
{\rm tr}(A_{15}\sigma_3)\right)+{\cal O}(\varepsilon).
\end{equation}
Now from (C.59) or directly from Equation~(\ref{eq:II-sigma-hat-derivative}) we find that
\begin{gather*}
t_6\frac{d\hat\sigma_6(t_6)}{dt_6}-\hat\sigma_6(t_6)={\rm tr}\,A_{06}A_{t6}=
\frac12({\rm tr}((A_{06}+A_{t6})^2-{\rm tr}(A_{06})^2-{\rm tr}(A_{t6})^2)\\
=\frac12\big({\rm tr}\left(\frac{\Theta_{\infty6}}2\sigma_3+A_{16}\right)^2-
{\rm tr}(A_{06})^2-{\rm tr}(A_{t6})^2\big)\\
=\frac{\Theta_{\infty6}}2{\rm tr}(A_{16}\sigma_3)+
\frac12\left({\rm tr}\left(\frac{\Theta_{\infty6}}2\sigma_3\right)^2+{\rm tr}\,(A_{16})^2-
{\rm tr}\,(A_{06})^2-{\rm tr}\,(A_{t6})^2\right).
\end{gather*}
Substituting into the latter formula Equations~(\ref{eq:lim2-A05-A15}) and
(\ref{eq:II-sigma-hat-deriv-A15-asympt})
and using definitions (\ref{eq:sigma-6}) and (\ref{eq:tau5-matrix-form}) we obtain
\begin{gather}
\frac{d}{dt_5}\ln\tau_6(t_6)=\left(\frac1{t_5}+\varepsilon+\varepsilon^2t_5\right)\!\!
\left(-\frac1{2\varepsilon^2}-\frac{\Theta_{\infty5}-t_5}{2\varepsilon}+
t_5{\rm tr}\Big(A_{15}\left(\frac1{t_5}A_{05}+\frac{\sigma_3}2\right)\!\!\Big)\right.
+\frac{\Theta_{\infty5}t_5}2\nonumber\\
 \label{eq:II-tau6-asympt}
\left.+\frac{\Theta_{05}^2+\Theta_{15}^2-\Theta_{\infty5}^2}2\right)
+{\cal O}(\varepsilon)=-\frac1{2\varepsilon^2t_5}-\frac{\Theta_{\infty5}}{2\varepsilon t_5}+
\frac{d}{dt_5}\ln\,\tau_5(t_5)+{\cal O}(\varepsilon).
\end{gather}
The second bracket on the r.h.s. of Equation (\ref{eq:II-tau6-asympt}) is nothing but the asymptotic
expansion (up to the order ${\cal O}(\varepsilon)$ for $-\hat\sigma_6(t_6)$. Differentiating it with
$t_5$ and using Equation~(\ref{eq:sigma'=-z5}) we find exactly
Equation ~(\ref{eq:II-sigma-hat-deriv-A15-asympt}). In an analogous manner the differentiable character
of the asymptotic expansions for $u_{\nu6}$ and $z_{\nu6}$ can be proved.

Inserting into the first Equation~(\ref{eq:y6}) asymptotics (\ref{eq:II-z16}) -- (\ref{eq:II-u16}) we find
that
\begin{equation}
 \label{eq:II-y6-asympt}
y_6(t_6)=\frac{1+{\cal O}(\varepsilon)}{1+\frac{y_5-1}{t_5}\left(z_5+
\frac{\Theta_{05}-\Theta_{15}+\Theta_{\infty5}}2-
\frac1{y_5}\left(z_5+\frac{\Theta_{05}+\Theta_{15}+\Theta_{\infty5}}2\right)\!\right)}.
\end{equation}
\section{The First Limit}
 \label{sec:4}
In this section we study the situation when the $\Psi_6$ function in the neighborhood of the
isomonodromic cluster of two regular singularities is described via the $\Psi_5$-function, i.e.,
the function with the irregular singular point. First of all let's define this cluster keeping in
mind the formulae (\ref{eq:lim1-t-lambda}) -- (\ref{eq:13a-s16}) corresponding to the first
limit {\bf I.3.a.} We won't consider the limit {\bf I.3.b.}, as it can be obtained via the
transformation (\ref{eq:reflection-infinity}).

Let us begin with the precise setting. We are considering the $\Psi_6(\lambda_6,t_6)$ function
defined as in Section~\ref{sec:2}. We assume that its manifold of monodromy data
$$
{\cal M}_6(\Theta_{06},\Theta_{16},\Theta_{t6},\Theta_{\infty6})
$$
is given. It is also assumed that
\begin{gather}
 \label{eq:limit1-theta-nu-notZ}
\Theta_{\nu6}\in\mathbb{C}\setminus\mathbb{Z},\qquad\nu=0,1,t,\infty,\\
\Theta_{\infty5}=\Theta_{06}+\Theta_{\infty6}\in\mathbb{C}.
 \label{eq:theta-infty5-definition}
\end{gather}
The function $\Psi_6(\lambda_6,t_6,n)$ is a GSD of $\Psi_6(\lambda_6,t_6)$ defined by the
following recurrence procedure:
\begin{gather*}
\Psi_6(\lambda_6,t_6,0)=\Psi_6(\lambda_6,t_6),\\
\Psi_6(\lambda_6,t_6,n+1)={\cal L}_{1\infty}^{-+}\Psi_6(\lambda_6,t_6,n),\quad
n\in\mathbb{Z}_+=\{0,1,\ldots\}.
\end{gather*}
The monodromy manifold ${\cal M}_6(\Theta_{06},\Theta_{16}^n,\Theta_{t6},\Theta_{\infty6}^n)$ for
$\Psi_6(\lambda_6,t_6,n)$ coincides with that for $\Psi_6(\lambda_6,t_6)$ except for the values of
the two parameters, $\Theta_{16}^n$ and $\Theta_{\infty6}^n$:
\begin{gather}
 \label{eq:theta-16n-epsilon-n}
\Theta_{16}^n=\Theta_{16}-2n\equiv-\frac1{\varepsilon_n},\qquad n\in\mathbb{Z}_+,\\
 \label{eq:theta-infty6n}
\Theta_{16}^n+\Theta_{\infty6}^n=\Theta_{\infty5}
\end{gather}
Equation~(\ref{eq:theta-16n-epsilon-n}) is the definition of the discrete small parameter
$\varepsilon=\varepsilon_n\to+0$, while the parameter $\Theta_{\infty5}$ in
Equation~(\ref{eq:theta-infty6n}) is defined by Equation~(\ref{eq:theta-infty5-definition}).

The main object of our investigation are the GSD's $A_{\nu6}^n(t_6)$, $\nu=0,1,t$. They can be
defined by means of Equation~(\ref{eq:linear-4-regular}):
\begin{equation}
 \label{eq:1-A-GSD-Psi}
\partial_{\lambda_6}\Psi_6(\lambda_6,t_6,n)\Psi_6^{-1}(\lambda_6,t_6,n)=\sum\limits_{\nu=0,1,t}
\frac{A_{\nu6}^n(t_6)}{\lambda-\nu}.
\end{equation}
An alternative (equivalent) definition of $A_{\nu6}^n(t_6)$ (without the usage of the auxiliary
object $\Psi_6(\lambda_6,t_6,n)$) can be given as the following recurrence procedure:
\begin{equation}
 \label{eq:1-A-GSD-rec}
A_{\nu6}^0(t_6)=A_{\nu6}(t_6),\qquad
A_{\nu6}^{n+1}={\tilde A}_{\nu6}^n(t_6),
\end{equation}
where $A_{\nu6}(t_6)$ is some solution of (\ref{eq:schlesinger-p6}), and
${\tilde A}_{\nu6}^n(t_6)$ is obtained via
Equations~(\ref{eq:A-tilde-mu-infty}) -- (\ref{eq:Atilde-infty2}) with $J_{1\infty}^{-+}$
by inserting $A_{\nu6}^n(t_6)$ into the r.-h.s.'s. In fact, in Theorem~\ref{th:1} we assume that
the matrices $\{A_{\nu6}^0(t_6)\}_{\nu=0,1,t}$ correspond to the manifold of monodromy data
${\cal M}_6(\Theta_{06},\Theta_{16},\Theta_{t6},\Theta_{\infty6})$; thus the function
$\Psi_6(\lambda_6,t_6,n)$ is also implicitly presented in the second definition of
$A_{\nu6}^n(t_6)$.

We interpret the formal limit transition {\bf I} as asymptotics (as $n\to+\infty$) of the even
sequences $\{A_{\nu6}^{2n}(\epsilon_{n} t_5\}_{\nu=0,1,t}$, or as asymptotics of the corresponding
sequences of their matrix elements. The asymptotic behavior of the odd sequences
$\{A_{\nu6}^{2n+1}(\epsilon_{n}t_5\}_{\nu=0,1,t}$ is given by exactly the same formulae as for the
even one (see formal limit {\bf I} in Section~\ref{sec:3} with $\varepsilon\to\varepsilon_{n}$.
The monodromy data for ``odd'' limit can be obtained from the monodromy data for the ``even'' limit
by simply changing in the latter formulae $\Theta_6\to\Theta_6-1$.

Sometimes, when it does not cause any confusion, we omit subscripts/superscripts $n$. To make a
difference between the initial values of the parameters $\Theta_{16}$ and $\Theta_{\infty6}$ in
Equations~(\ref{eq:limit1-theta-nu-notZ}), (\ref{eq:theta-infty5-definition}) and the parameters
$\Theta_{16}^n$ and $\Theta_{\infty6}^n$, we agree to denote the pair of initial values as
\begin{equation}
 \label{eq:I-initial-theta}
\Theta_{16}^0\equiv\Theta_6\in\mathbb{C}\setminus\mathbb{Z},\qquad
\Theta_{\infty6}^0=\Theta_{\infty5}-\Theta_6\in\mathbb{C}\setminus\mathbb{Z},
\end{equation}
while, instead of $\Theta_{16}^n$ and $\Theta_{\infty6}^n$, we write $\Theta_{16}$ and
$\Theta_{\infty6}$, respectively. {\it Thus, hereafter in Section~{\rm\ref{sec:4}}, $\Theta_{6}$
and $\Theta_{\infty5}$ are fixed according to {\rm(\ref{eq:I-initial-theta})} and
{\rm(\ref{eq:theta-infty5-definition})}, while $\Theta_{16}$ and $\Theta_{\infty6}$ are dependent
on $n$ such that $\Theta_{16}\underset{n\to+\infty}\to-\infty$ and
$\Theta_{\infty6}\underset{n\to+\infty}\to+\infty$ and the condition~{\rm(\ref{eq:limit1-theta-nu-notZ})}
holds}.

In fact the formal limit {\bf I} contains one more parameter $f_0\in\mathbb{C}\setminus\{0\}$. This
parameter is hidden as the constant of integration in Equation~(\ref{eq:13a-s16}). The asymptotic
expansions of the sequences under investigation
\begin{equation}
 \label{eq:I-even-A}
 A_{\nu6}^{2n}(\epsilon_{n}t_5),\qquad n\in\mathbb{Z}_+,\quad\nu=0,1,t,
\end{equation}
don't depend on $f_0$. On the language of the formulae~(\ref{eq:13a-z06-zt6}) -- (\ref{eq:13a-u16})
this means that if $s_{16}\to s_{16}f_0$, then $u_5\to u_5f_0^2$. Nevertheless we include $f_0$
for completeness. Now, denoting the matrix elements of (\ref{eq:I-even-A}) exactly as that for
$A_{\nu6}$ in (\ref{eq:A-and-R}), and the matrix elements of the monodromy matrices $M_{\nu6}$ as
\begin{equation}
 \label{eq:I-Mnu6-elements-def}
M_{\nu6}=\left(\begin{array}{cc}
m_{11}^{\nu6}&m_{12}^{\nu6}\\
m_{21}^{\nu6}&m_{22}^{\nu6}
\end{array}\right),
\end{equation}
we are ready to formulate our result.
\begin{theorem}
 \label{th:1}
Assume that the complex parameters $\Theta_06$, $\Theta_6$, $\Theta_{t6}$ satisfy Conditions
{\rm(\ref{eq:limit1-theta-nu-notZ})} and {\rm(\ref{eq:I-initial-theta})}. Let a point
$\mu_6\in{\cal M}_6(\Theta_{06},\Theta_6,\Theta_{t6},\Theta_{\infty5}-\Theta_6)$ and
$\{m_{ik}^{\nu6}\}_{i,k=1,2}^{\nu=0,1,t}$ be its coordinates. Suppose that the following
conditions are valid:\\
{\rm1}. $m_{21}^{16}\neq0$;\\
{\rm2}.
\begin{equation}
 \label{eq:I-l-condition}
\cos\,\pi l\equiv=im_{11}^{16}\sin\,\pi(\Theta_{\infty5}-\Theta_6)+
\exp(-\pi i(\Theta_{\infty5}-\Theta_6))\cos\,\pi\Theta_6\neq-1;
\end{equation}\\
{\rm3}. The inverse monodromy problem defined by the pairs $(\mu_6,t_6)$ are solvable for all
$t_6=\varepsilon_{n}t_5$ where $n\geq N\in\mathbb{Z}_+$ and $\varepsilon_{n}$ is defined in
Equation~{\rm(\ref{eq:theta-16n-epsilon-n})}, and $t_5\in\mathbb{C}\setminus\{0\}$ with
$|\arg t_5|<\pi-\delta$ for some $\delta>0$;\\
{\rm4}. The point $\mu_6$ corresponds to the general, i.e., transcendental solution of
$P_6$\footnote[6]{We mean the $P_6$-transcendent, i.e., any solution which cannot be constructed
in terms of the logarithmic derivatives of the hypergeometric functions.}.

Let $l$ be the unique solution of Equation~(\ref{eq:I-l-condition}) under the conditions:
\begin{equation}
 \label{eq:I-l-conditions}
l\neq0,\qquad |{\rm Re}\,l|<1.
\end{equation}
Define the parameters $\alpha$, $\beta$, $d_0^2$:
\begin{equation}
 \label{eq:alpha-beta-d0}
\alpha=\frac{l-\Theta_{\infty5}}2,\qquad
\beta=-\frac{l+\Theta_{\infty5}}2,\qquad
d_0^2=\frac{m_{21}^{16}}{2\pi i}\Gamma(1-\alpha)\Gamma(1-\beta),
\end{equation}
where $\Gamma(\cdot)$ is the gamma function~{\rm\cite{BE}}.

For arbitrary $f_0\in\mathbb{C}\setminus\{0\}$ define the matrix
\begin{equation}
 \label{eq:I-K}
K=-f_0^{\sigma_3}e^{-\frac{\pi i}2\Theta_6\sigma_3}
\left(\begin{array}{cc}
1&\frac{\pi}{\Gamma(\alpha)\Gamma(\beta)\sin\,\pi(\Theta_{\infty5}-\Theta_6)}\\
0&1
\end{array}\right)
d_0^{\sigma_3},
\end{equation}
and consider the following equations:
\begin{gather}
 \label{eq:I-def-theta's}
\Theta_{05}=\Theta_{06},\qquad
\Theta_{15}=-\Theta_{t6},\qquad
\Theta_{\infty5}=\Theta_{16}+\Theta_{\infty6}=-(\alpha+\beta),\\
 \label{eq:I-def-monodromy5}
\tilde M_{05}=KM_{06}K^{-1},\qquad
\tilde M_{15}=KM_{t6}K^{-1},\qquad
\tilde M_{\infty5}=KM_{\infty6}M_{16}K^{-1},
\end{gather}
as defining the point $\tilde\mu_5\in\tilde{\cal M}_5(\Theta_{05},\Theta_{15},\Theta_{\infty5})$.
Suppose that the inverse monodromy problem for $(\tilde\mu_5,t_5)$ is solvable and the functions
$u_5=u_5(\tilde\mu_5,t_5)$, $y_5=y_5(\tilde\mu_5,t_5)$, and $z_5=z_5(\tilde\mu_5,t_5)$ represent
this solution.

Then it is possible to construct the sequences~{\rm(\ref{eq:I-even-A})}, corresponding to the given
manifold ${\cal M}_6(\Theta_{06},\Theta_6,\Theta_{t6},\Theta_{\infty5}-\Theta_6)$.
The formulae~{\rm(\ref{eq:13a-theta})--(\ref{eq:13a-s16})} are asymptotic expansions as
$n\to+\infty$ of the matrix elements of {\rm(\ref{eq:I-even-A})} if:\\
{\rm1}. The functions $u_5$, $y_5$, and $z_5$ are identified with the solution of the inverse
monodromy problem for $(\tilde\mu_5,t_5)$;\\
{\rm2}. Equation~{\rm(\ref{eq:13a-s16})} is supplemented with
\begin{equation}
 \label{eq:I-asympt-s16}
s_{16}=\varepsilon f_0d_0(\varepsilon t_5)^{\frac{\Theta_{\infty5}}2}(1+{\cal O}(\varepsilon));
\end{equation}
{\rm3}. $\varepsilon=\varepsilon_{n}$ is substituted into
{\rm(\ref{eq:13a-theta})--(\ref{eq:13a-s16})} and {\rm(\ref{eq:I-asympt-s16})}.
\end{theorem}
\begin{remark}{\rm
The complex number $l$ is defined by Conditions (\ref{eq:I-l-condition}) and
(\ref{eq:I-l-conditions}) up to a sign, which means simply the permutation
$\alpha\leftrightarrow\beta$. This permutation doesn't
influence the asymptotics that we study, as well as the indefiniteness of the sign $d_0$ in
(\ref{eq:alpha-beta-d0}).
}\end{remark}
\begin{remark}{\rm
Some of the conditions imposed on $\mu_6$ (see Conditions $1$, $2$, $4$) are not necessary for the
possibility to interpret the formal limit as an asymptotic expansion; however, their violation
requires special investigation. In particular, we can apply this theorem not only to the $P_6$
transcendent, but also for most of the solutions that can be constructed via the classical special
functions.
}\end{remark}
\begin{remark}{\rm
To find the Stokes multipliers for the limiting $P_5$ equation one can use
Equation~(\ref{eq:M-infty-stokes-k}).
}\end{remark}
\begin{derivation}
Here I outline only the calculational scheme for the solution of the direct and inverse
monodromy problem, while an explanation of how such calculations work for the justification of
the asymptotic expansions, the reader will find in Section~\ref{sec:5} in the derivation of
Theorem~\ref{th:2}. Although the derivations are different, the scheme of justification is the same
and can be based on the work~\cite{K4}.

For some $\delta_{in}$, $\delta_{out}$ such that $\frac13<\delta_{in}<\delta_{out}<1$ and for all
rather small $\varepsilon>0$ we can present the $\lambda_6$-complex plane as
$\mathbb{C}=\Omega_{in}(\varepsilon)\cup\Omega_{out}(\varepsilon)$, where
$$
\lambda_6\in\Omega_{in}(\varepsilon)\Leftrightarrow
|\lambda_6|\leq{\cal O}(\varepsilon^{\frac12+\frac{\delta_{in}}2}),\qquad
\lambda_6\in\Omega_{out}(\varepsilon)\Leftrightarrow
|\lambda_6|\geq{\cal O}(\varepsilon^{\frac12+\frac{\delta_{out}}2}).
$$
Thus the matching domain,
$\Omega_{mat}(\varepsilon)=\Omega_{in}(\varepsilon)\cap\Omega_{out}(\varepsilon)$, is nonempty for
all rather small $\varepsilon$.

In the cluster domain we approximate the function $\Psi_6=\Psi_6(\lambda_6,t_6,n)$ as follows
\begin{equation}
 \label{eq:I-psi-approx-in}
\Psi_6\underset{\varepsilon\to0}=R_{16}(1+{\cal O}(\varepsilon^\delta))\Psi_5^0K,\qquad
\lambda_6\in\Omega_{in}(\varepsilon),\quad\delta>0,
\end{equation}
where $\Psi_5^0=\Psi_5^0(\lambda_5,t_5)$ is the canonical solution of Equation~(\ref{eq:linear-5})
with monodromy data to be determined, as well as the constant matrix $K\in{\rm SL}(2,\mathbb{C})$.
Using Approximation~(\ref{eq:I-psi-approx-in}) one finds
Equations~(\ref{eq:13a-theta}) -- (\ref{eq:13a-s16}), (\ref{eq:I-def-monodromy5}), and the first two
equations in (\ref{eq:I-def-theta's}). We have to prove (\ref{eq:I-K}), the last equation in
(\ref{eq:I-def-theta's}) and (\ref{eq:I-asympt-s16}), as well as for the justification, it is
important to establish the following
\begin{equation}
 \label{eq:K-matching}
K=\underset{\begin{array}{c}\varepsilon\to0\\\lambda_0\in\Omega_{mat}(\varepsilon)\end{array}}\lim
(R_{16})(\Psi_5^0)^{-1}\Psi_{out},
\end{equation}
where $\Psi_{out}$ is an approximation for $\Psi_6$ in the domain $\Omega_{out}(\varepsilon)$:
\begin{equation}
 \label{eq:I-psi-out-approx}
\Psi_6=(I+{\cal O}(\varepsilon^\delta))\Psi_{out},\qquad\lambda_6\in\Omega_{out}(\varepsilon),
\;\;\delta>0.
\end{equation}

We construct $\Psi_{out}$ by making use of the function $Y(x;\alpha,\beta,\gamma)$ constructed by
Jimbo~\cite{J}. While constructing $\Psi_{out}$ we find the conditions
(\ref{eq:I-l-condition}) -- (\ref{eq:alpha-beta-d0}) and $m_{21}^{16}\neq0$. Since the function
$Y(x;\alpha,\beta,\gamma)$ plays an important role not only in the present derivation but in the one
in Section~\ref{sec:5}, we, following~\cite{J}, recall its basic properties:
\begin{gather}
Y(x)\equiv Y(x;\alpha,\beta,\gamma)x^{\frac{\gamma-1}2}(x-1)^{\frac{\alpha+\beta+1-\gamma}2}\nonumber\\
=\left(\begin{array}{c}
F(\alpha,\alpha+1-\gamma,\alpha-\beta;\frac1x,\frac{\beta(\beta+1-\gamma)}{(\beta-\alpha)(\beta-\alpha+1)x}
F(\beta,\beta+2-\gamma,\beta-\alpha+2;\frac1x)\\
\frac{\alpha(\alpha+1-\gamma)}{(\alpha-\beta)(\alpha-\beta+1)x}
F(\alpha+1,\alpha+2-\gamma,\alpha-\beta+2;\frac1x),F(\beta,\beta+1-\gamma,\beta-\alpha;\frac1x)
\end{array}\right)
 \label{eq:Y(x)-def}
\\
\times x^{\frac{\beta-\alpha}2\sigma_3}x^{\frac{\gamma-1-\alpha-\beta}2}
(x-1)^{\frac{\alpha+\beta+1-\gamma}2},\nonumber\\
 \label{eq:alpha,beta,gamma-conditions-Y}
\Theta_{\infty Y}\equiv\alpha-\beta,\;\Theta_{0Y}\equiv1-\gamma,\;
\Theta_{1Y}\equiv-\alpha-\beta-1+\gamma\in\mathbb{C\setminus Z},
\end{gather}
where $F(\cdot,\cdot,\cdot;\cdot)$ denotes the Gauss hypergeometric function~\cite{BE}.
\begin{equation}
 \label{eq: Y(x)-asympt}
Y(x)\underset{x\to\nu=0,1}=G_{\alpha\beta\gamma}^\nu(I+{\cal O}(x-\nu))
(x-\nu)^{\frac{\Theta_{\nu Y}}2\sigma_3}C_{\alpha\beta\gamma}^\nu\underset{x\to\infty}
=\left(I+{\cal O}(x^{-1})\right)x^{\frac{\beta-\alpha}2\sigma_3},
\end{equation}
where
\begin{align}
 \label{eq:I-G0-G1-alpha-beta-gamma}
G_{\alpha\beta\gamma}^0&=\frac1{\beta-\alpha}\left(\begin{array}{cc}
\beta+1-\gamma&\beta\\
\alpha+1-\gamma&\alpha
\end{array}\right),\quad
G_{\alpha\beta\gamma}^1=\frac1{\beta-\alpha}\left(\begin{array}{cc}
1&\beta(\beta+1-\gamma)\\
1&\alpha(\alpha+1-\gamma)
\end{array}\right),\\
 \label{eq:C0-alpha-beta-gamma}
C_{\alpha\beta\gamma}^0&=\left(\begin{array}{cc}
e^{-\pi i(\alpha+1-\gamma)}
\frac{\Gamma(\gamma-1)\Gamma(\alpha-\beta+1)}{\Gamma(\gamma-\beta)\Gamma(\alpha)}&
-e^{-\pi i(\beta+1-\gamma)}
\frac{\Gamma(\gamma-1)\Gamma(\beta-\alpha+1)}{\Gamma(\gamma-\alpha)\Gamma(\beta)}\\
e^{-\pi i\alpha}\frac{\Gamma(1-\gamma)\Gamma(\alpha-\beta+1)}{\Gamma(1-\beta)\Gamma(\alpha+1-\gamma)}&
-e^{-\pi i\beta}\frac{\Gamma(1-\gamma)\Gamma(\beta-\alpha+1)}{\Gamma(1-\alpha)\Gamma(\beta+1-\gamma)}
\end{array}\right),\\
 \label{eq:C1-alpha-beta-gamma}
C_{\alpha\beta\gamma}^1&=\left(\begin{array}{cc}
-\frac{\Gamma(\alpha+\beta+1-\gamma)\Gamma(\alpha-\beta+1)}{\Gamma(\alpha+1-\gamma)\Gamma(\alpha)}&
\frac{\Gamma(\alpha+\beta+1-\gamma)\Gamma(\beta-\alpha+1)}{\Gamma(\beta+1-\gamma)\Gamma(\beta)}\\
-e^{-\pi i(\gamma-\alpha-\beta-1}\frac{\Gamma(\gamma-\alpha-\beta-1)\Gamma(\alpha-\beta+1)}{\Gamma(1-\beta)\Gamma(\gamma-\beta)}&
e^{-\pi i(\gamma-\alpha-\beta-1)}\frac{\Gamma(\gamma-\alpha-\beta-1)\Gamma(\beta-\alpha+1)}{\Gamma(1-\alpha)\Gamma(\gamma-\alpha)}
\end{array}\right),
\end{align}
\begin{align}
 \label{eq:hypergeo-Y(x)}
\frac{dY(x)}{dx}&=\left(\frac{A_{0Y}}x+\frac{A_{1Y}}{x-1}\right)Y(x),\qquad
A_{0Y}+A_{1Y}=\frac{\beta-\alpha}2\sigma_3,\\
 \label{eq:A0Y}
A_{0Y}&=\frac1{\beta-\alpha}\left(\begin{array}{cc}
-\frac{(\alpha+\beta)(1-\gamma)+2\alpha\beta}2&\beta(\beta+1-\gamma)\\
-\alpha(\alpha+1-\gamma)&\frac{(\alpha+\beta)(1-\gamma)+2\alpha\beta}2
\end{array}\right),\\
 \label{eq:A1Y}
A_{1Y}&=\frac1{\beta-\alpha}\left(\begin{array}{cc}
-\frac{(\alpha+\beta)(1-\gamma)+\alpha^2+\beta^2}2&-\beta(\beta+1-\gamma)\\
-\alpha(\alpha+1-\gamma)&-\frac{(\alpha+\beta)(1-\gamma)+\alpha^2+\beta^2}2
\end{array}\right).
\end{align}
To the above-mentioned properties of $Y(x)$ pointed out by Jimbo, we add up
the following one, which is important in the derivation of our main results
in this and the next section:
\begin{gather}
 \label{eq:Y(x)-spec-asympt}
Y(x)\underset{x\in(\ref{eq:x-spec-asympt-cond})}=\frac1{\sqrt{\alpha-\beta}}
\left(\!\begin{array}{cc}
1&\beta\\
1&\alpha
\end{array}\!\right)\!\!
\left(I+{\cal O}\Big(x^{\tilde\delta}\Big)\!\right)
\exp\left(\!\!\Big(\frac{1-\gamma}{2x}+\frac{\alpha+\beta}2\ln\frac{1-\gamma}x\Big)\sigma_3\!\right)\!
\hat K(\varkappa),\\
 \label{eq:K(varkappa)}
\hat K(\varkappa)=\left(I+{\cal O}\Big(x^{\tilde\delta}\Big)\!\right)\sqrt{\alpha-\beta}\sigma_3
\left(\!\begin{array}{cc}
\frac{\Gamma(\alpha-\beta)}{\Gamma(\alpha)}&\frac{\Gamma(\beta-\alpha)}{\Gamma(\beta)}\\
-e^{i\pi\alpha\varkappa}\frac{\Gamma(1+\alpha-\beta)}{\Gamma(1-\beta)}&e^{i\pi\beta\varkappa}
\frac{\Gamma(1+\beta-\alpha)}{\Gamma(1-\alpha)}
\end{array}\!\right)(1-\gamma)^{\frac{\beta-\alpha}2\sigma_3},
\end{gather}
\begin{equation}
 \label{eq:x-spec-asympt-cond}
\left|\frac{1-\gamma}{x^2}\right|\underset{\varepsilon\to+0}={\cal O}\big(x^{\tilde\delta}\big),\quad
\left|\frac{1-\gamma}x\right|\underset{\varepsilon\to+0}=+\infty,\quad
-\pi+\frac{\varkappa\pi}2<\arg\frac{1-\gamma}x<\pi+\frac{\varkappa\pi}2,\;\;\varkappa=\pm1.
\end{equation}
In Equation~(\ref{eq:Y(x)-spec-asympt}), and thereafter, the notation
${\cal O}\big(x^{\tilde\delta}\big)$ means
that the corresponding estimate holds for some $\tilde\delta>0$. The precise (the largest possible)
value of $\tilde\delta$ in such estimates are not important in our scheme of derivation.

For a proof of Asymptotics~(\ref{eq:Y(x)-spec-asympt}) we need a special asymptotic expansion for the
Gauss hypergeometric function $F(a,b,c;x)$ which can be found in \cite{BE}:
\begin{equation}
 \label{eq:hyper-asympt-spec}
\begin{gathered}
F(a,b,c;z)=e^{i\pi a\varkappa}\frac{\Gamma(c)}{\Gamma(c-a)}(bz)^{-a}(1+{\cal O}(|bz|^{-1}))
+\frac{\Gamma(c)}{\Gamma(a)}e^{bz}(bz)^{a-c}(1+{\cal O}(|bz|^{-1})),\\
0<|z|<1,\quad|bz|\to\infty,\quad-\pi+\frac{\varkappa\pi}2<\arg(bz)<\pi+\frac{\varkappa\pi}2,\;\;
\varkappa=\pm1.
\end{gathered}
\end{equation}
In the domain~(\ref{eq:x-spec-asympt-cond}) one finds:
\begin{equation}
 \label{eq:x-power-aympt-simp}
x^{\frac{\beta-\alpha}2\sigma_3}x^{\frac{\gamma-1-\alpha-\beta}2}(x-1)^{\frac{\alpha+\beta+1-\gamma}2}=
x^{\frac{\beta-\alpha}2\sigma_3}e^{-\frac{\alpha+\beta+1-\gamma}{2x}}
\left(1+{\cal O}\big(x^{\tilde\delta}\big)\!\right).
\end{equation}
Substituting Expansions (\ref{eq:hyper-asympt-spec}) and (\ref{eq:x-power-aympt-simp}) into
Equation~(\ref{eq:Y(x)-def}) one arrives at
Equations~(\ref{eq:Y(x)-spec-asympt})--(\ref{eq:x-spec-asympt-cond}).

We construct the function $\Psi_{out}$ as follows:
\begin{gather}
\Psi_{out}(\lambda_6)=G^{-1}Y(x)C_0^{-1},\qquad x=\frac1{\lambda_6},\nonumber\\
G=G_{\alpha\beta\gamma}^0(\beta-\alpha)(1-\gamma)^{-\frac12}d_0^{\sigma_3}
(1-\gamma)^{\frac{\alpha+\beta-1}2\sigma_3},
 \label{eq:I-G-Psi-out}\\
C_0=e^{\frac{\pi i}2(1-\gamma+\alpha+\beta)}(1-\gamma)^{\frac12}(\beta-\alpha)^{-1}
(1-\gamma)^{\frac{1-\alpha-\beta}2\sigma_3}d_0^{-\sigma_3}C_{\alpha\beta\gamma}^0,
 \label{eq:I-C-Psi-out}
\end{gather}
where $\alpha$ and $\beta$, satisfying the conditions~(\ref{eq:alpha,beta,gamma-conditions-Y}),
and $d_0\in\mathbb{C}\setminus\{0\}$ are the parameters to be determined, and
\begin{equation}
 \label{eq:I-gamma-cond}
1-\gamma=\Theta_{\infty5}+\frac1{\varepsilon_{n}}\to+\infty,\quad
1-\gamma=\Theta_{\infty5}-\Theta_6+2n,\quad
\gamma-\alpha-\beta-1=\Theta_6-2n,
\end{equation}
as it follows from Equations~(\ref{eq:theta-16n-epsilon-n}), (\ref{eq:I-initial-theta}),
(\ref{eq:I-def-theta's}), and (\ref{eq: Y(x)-asympt}).

By using (\ref{eq: Y(x)-asympt}) and (\ref{eq:Y(x)-spec-asympt}) - (\ref{eq:x-spec-asympt-cond})
we find that the function $\Psi_{out}$ has the following asymptotic behavior
\begin{align}
 \label{eq:Psi-out-asympt-infty}
\lambda_6&\to\infty:&\Psi_{out}&\sim\left(\frac1{\lambda_6}\right)^{\frac{\Theta_{\infty6}}2\sigma_3},\\
 \label{eq:Psi-out-asympt-1}
\lambda_6&\to1:&\Psi_{out}&\sim G^{-1}G_{\alpha\beta\gamma}^1
\left(\frac1{\lambda_6}-1\right)^{\frac{\Theta_{16}}2\sigma_3}C_{\alpha\beta\gamma}^1C_0^{-1},\\
\lambda_6&\to0,&\lambda_6&\in\Omega_{mat}(\varepsilon),\qquad-\pi+\frac{\varkappa\pi}2<\arg\lambda_6<
\pi+\frac{\kappa\pi}2,\nonumber
\end{align}
\begin{equation}
 \label{eq:Psi-out-asympt-0}
\Psi_{out}\sim\frac{G^{-1}}{\sqrt{\alpha-\beta}}
\left(\begin{array}{cc}
1&\beta\\
1&\alpha
\end{array}\right)
\exp\left(\!\Big(\frac{(1-\gamma)\lambda_6}2-
\frac{\alpha+\beta}2\ln\frac1{(1-\gamma)\lambda_6}\Big)\sigma_3\!\right)\hat K(\varkappa)C_0^{-1}.
\end{equation}
We see from (\ref{eq:Psi-out-asympt-infty}) that $\Psi_{out}$ satisfies the same normalization
condition as the function $\Psi_6$. Thus, the monodromy matrices at $\lambda_6=\infty$ for both
functions coincide, $M_{\infty\,out}=M_{\infty6}$. Our next step will be to set the parameters $\alpha-\beta$ and
$d_0$ such that the following equation
\begin{equation}
 \label{eq:I-mon-matrix-cond-1}
M_{1\,out}=M_{16}
\end{equation}
holds. To calculate asymptotically $M_{1\,out}$ we use Expansion~(\ref{eq:Psi-out-asympt-1}),
Equations (\ref{eq:I-C-Psi-out}), (\ref{eq:C0-alpha-beta-gamma}) and the following formulae for
the $\Gamma$-function~\cite{BE}:
\begin{gather}
 \label{eq:Gamma-ratio}
\frac{\Gamma(z+\mu)}{\Gamma(z+\nu)}=z^{\mu-\nu}(1+{\cal O}(z^{-1})),
\qquad |z|\to\infty,\qquad|\arg\,z|<\pi,\\
 \label{eq:Gamma-product}
\Gamma(z)\Gamma(1-z)=\frac\pi{\sin(\pi z)}.
\end{gather}
One notices that
$$
\det C_{\alpha\beta\gamma}^0=\frac{\gamma-1}{\beta-\alpha}e^{-\pi i(\alpha+\beta+1-\gamma)},
$$
and finds
\begin{equation}
 \label{eq:I-C0-C}
 \begin{gathered}
C_0^{-1}=(1-\gamma)^{\frac{\alpha-\beta}2\sigma_3}e^{\frac{\pi i}2(\alpha+\beta+1-\gamma)}
\left(C+{\cal O}\big(x^{\tilde\delta}\big)\!\right)d_0^{\sigma_3},\\
C=\left(\begin{array}{cc}
\frac{\Gamma(\beta-\alpha+1)}{\Gamma(1-\alpha)}e^{\pi i(1-\beta)}&
\frac{\Gamma(\beta-\alpha+1)}{\Gamma(\beta)}\frac{\sin\,\pi(\alpha+1-\gamma)}{\sin\,\pi(1-\gamma)}
e^{-\pi i(\beta-\gamma)}\\
\frac{\Gamma(\alpha-\beta+1)}{\Gamma(1-\beta)}e^{\pi i(1-\alpha)}&
\frac{\Gamma(\alpha-\beta+1)}{\Gamma(\alpha)}\frac{\sin\,\pi(\beta+1-\gamma)}{\sin\,\pi(1-\gamma)}
e^{-\pi i(\alpha-\gamma)}
\end{array}\right).
 \end{gathered}
\end{equation}
Note that according to the last two equations (\ref{eq:I-gamma-cond}) $C$ is the constant matrix.
Using Equations~(\ref{eq:Gamma-ratio}) and (\ref{eq:Gamma-product}) we obtain the leading term of
asymptotics for $C_{\alpha\beta\gamma}^1$ (\ref{eq:C1-alpha-beta-gamma}):
\begin{equation}
 \label{eq:I-C1-alpha-beta-gamma-asympt}
\begin{gathered}
C_{\alpha\beta\gamma}^1\underset{\varepsilon\to0}\sim(1-\gamma)^{\frac{\alpha+\beta}2\sigma_3}
\left(\begin{array}{cc}
1&0\\
0&\frac{e^{\pi i(\alpha+\beta+1-\gamma)}}{1-\gamma}
\end{array}\right)\\
\times\left(\begin{array}{cc}
-\frac{\Gamma(\alpha-\beta+1)}{\Gamma(\alpha)}&\frac{\Gamma(\beta-\alpha+1)}{\Gamma(\beta)}\\
\frac{\sin\,\pi(\beta+1-\gamma)}{\sin\,\pi(\alpha+\beta+1-\gamma)}
\frac{\Gamma(\alpha-\beta+1)}{\Gamma(1-\beta)}&
-\frac{\sin\,\pi(\alpha+1-\gamma)}{\sin\,\pi(\alpha+\beta+1-\gamma)}
\frac{\Gamma(\beta-\alpha+1)}{\Gamma(1-\alpha)}
\end{array}\right)
(1-\gamma)^{\frac{\beta-\alpha}2\sigma_3}.
\end{gathered}
\end{equation}
Now Equations (\ref{eq:I-C0-C}) and (\ref{eq:I-C1-alpha-beta-gamma-asympt}) yield
\begin{equation}
 \label{eq:M1out}
M_{1\,out}=C_1^{-1}e^{\pi i(\gamma-1-\alpha-\beta)}C_1.
\end{equation}
Hence Equation~(\ref{eq:I-mon-matrix-cond-1}) in the notation~(\ref{eq:I-Mnu6-elements-def}) reads as
\begin{align}
 \label{eq:I-m-11-16}
m_{11}^{16}&=\frac{i}{\sin\,\pi(1-\gamma)}(\cos\,\pi(\alpha+\beta+1-\gamma)e^{-\pi i(1-\gamma)}-
\cos\,\pi(\alpha-\beta)),\\
m_{22}^{16}&=\frac{i}{\sin\,\pi(1-\gamma)}(\cos\,\pi(\alpha-\beta)-
\cos\,\pi(\alpha+\beta+1-\gamma)e^{\pi i(1-\gamma)}),\nonumber\\
 \label{eq:I-m-21-16=m12-16}
m_{21}^{16}&=\frac{2\pi id_0^2}{\Gamma(1-\alpha)\Gamma(1-\beta)},\quad
m_{12}^{16}=-\frac{2\pi i\sin\,\pi(\alpha+1-\gamma)\sin\,\pi(\beta+1-\gamma)}
{d_0^2\Gamma(\alpha)\Gamma(\beta)\sin^2\pi(1-\gamma)}.
\end{align}
Now define $l=\alpha-\beta$ and recall (\ref{eq:I-gamma-cond}) to derive from
Equations~(\ref{eq:I-m-11-16}) and (\ref{eq:I-m-21-16=m12-16}) Formulae (\ref{eq:I-l-condition})
and (\ref{eq:alpha-beta-d0}), respectively.
It was supposed that $l\neq0,\pm1$ (see Equation~(\ref{eq:alpha,beta,gamma-conditions-Y})).
The natural requirement $|{\rm Re}|\,l<1$ means no additional restrictions, since, thanks to
the equation $\alpha+\beta=-\Theta_{\infty5}$, the shift $l\to l+2$ leads, simply, to a
redefinition of $d_0$; thus we get (\ref{eq:I-l-conditions}).

Turning to the matching (\ref{eq:K-matching}), we notice that $(1-\gamma)/x=\lambda_5t_5$; and by
recalling asymptotics of $\Psi_5^0$ (\ref{eq:Psi5-asymptotics-infty}) we confirm not only the
last equation (\ref{eq:I-def-theta's}), but also obtain that
\begin{equation}
 \label{eq:R16-asympt}
R_{16}=G^{-1}\left(\begin{array}{cc}
1&\beta\\
1&\alpha
\end{array}\right)
\left(1+{\cal O}\big(x^{\tilde\delta}\big)\!\right)
t_5^{-\frac{\Theta_{\infty5}}2\sigma_3}f_0^{-\sigma_3},\qquad\in\mathbb{C}\setminus\{0\},
\end{equation}
\begin{equation}
 \label{eq:K-limit}
K=\underset{\varepsilon\to0}\lim\,f_0^{\sigma_3}\hat K(-1)C_0^{-1}/\sqrt{\alpha-\beta}.
\end{equation}
We can further simplify Equation~(\ref{eq:R16-asympt}) with the help of
Equations~(\ref{eq:I-G-Psi-out}) and (\ref{eq:I-G0-G1-alpha-beta-gamma}):
\begin{equation}
 \label{eq:r16-asympt-final}
R_{16}=d_0^{-\sigma_3}(1-\gamma)^{\frac{\Theta_{\infty5}}2\sigma_3}
\left(1+{\cal O}\big(x^{\tilde\delta}\big)\!\right)
t_5^{-\frac{\Theta_{\infty5}}2\sigma_3}f_0^{-\sigma_3}.
\end{equation}
Equation~(\ref{eq:r16-asympt-final}) is equivalent, up to the leading term, to
Equation~(\ref{eq:I-asympt-s16}).
Finally, substituting into (\ref{eq:K-limit}) formulae (\ref{eq:K(varkappa)}), (\ref{eq:I-gamma-cond}),
and (\ref{eq:I-C0-C}), we arrive at Equation~(\ref{eq:I-K}).
\end{derivation}
\section{The Second Limit}
 \label{sec:5}
Comparing the formulae for $s_{16}$ (\ref{eq:13a-s16}) and (\ref{eq:13b-s16}) for the first limit
passage with the analogous formula for $s_{t6}$ (\ref{eq:II-st6}) for the second limit, one finds
the ``principle'' distinction between the limits: while Equations~(\ref{eq:13a-s16}) and
(\ref{eq:13b-s16}) are ``integrable'' Equation \ref{eq:II-st6}) is not. Thus it is not clear how
to set the constant of integration in the definition of $s_{t6}$. As the result asymptotics of
$u_{\nu6}$, $\nu=0,1,t$ is found here up to the factor $c$: $s_{t6}\to cs_{t6}$ $\Leftarrow$
$u_{\nu6}\to c^{-2}u_{\nu6}$. This fact does not influence the functions $y_6$ and $\tau_6$,
whose asymptotics are properly defined in the case of the second limit passage (see below
Theorem~\ref{th:2}). The problem of how to cope with the ambiguity of $c$, i.e., to set $s_{t6}$ in
terms of ${\cal M}_6(\Theta_{06},\Theta_{16},\Theta_{t6},\Theta_{\infty6})$, is left for further
investigation. This explains some differences appearing in the formulations of Theorems~\ref{th:1} and
~\ref{th:2}; here we also omit the nonessential $f_0$-like parameter (see (\ref{eq:I-K})).

As it is mentioned in Introduction on the level of the Painlev\'e and the corresponding $\tau$-functions
both limits are equivalent. We show that at the end of this section in Proposition~\ref{prop:1}.
Nevertheless we present also the direct derivation because:\\
1. This result is easy to extend for the cluster on the arbitrary regular background, i.e., for
the Garnier systems (see Corollary~\ref{cor:1}).\\
2. It is interesting to see how the ``nonintegrability'' of Equation~(\ref{eq:II-st6}) manifests
itself in our asymptotic calculations which do not contain any integration in the usual sense.

Actually I started my studies with the second cluster as it naturally appeared when one inserts
the expansion $(\lambda_6-t_6)^{-1}=\lambda_6^{-1}+t_6\lambda_6^{-2}+\ldots$ into
Equation~(\ref{eq:linear-4-regular}). After I realized that the above ``$c$-problem'' is not the
intrinsic ``cluster problem'', I found the first limit, where such a ``$c$-problem'' does not
appear. Note that applying Proposition~\ref{prop:1} to Theorem~\ref{th:2} we cannot
obtain the complete result for the first limit, as it stated in Theorem~\ref{th:1}.

Let us begin with the regularization of the second limit. Suppose that\\
${\cal M}_6(\Theta_{06},\Theta_{16},\Theta_{t6},\Theta_{\infty6})$ is given and define the
parameters $\varepsilon_n$, $\Theta_{\infty5}$, $\Theta_6$ as follows:
\begin{gather}
 \label{eq:II-theta-t6}
\Theta_{t6}=\Theta_6-2n\equiv-\frac1{\varepsilon_n},\quad n\in\mathbb{Z}_+,\quad
\Theta_6\in\mathbb{C}\setminus\mathbb{Z},\\
 \label{eq:II-def-theta-infty5}
\Theta_{t6}+\Theta_{06}=\Theta_{\infty5},\qquad\Theta_{\infty5}-\Theta_6\in\mathbb{C}\setminus\mathbb{Z}.
\end{gather}
Here again our parameter $\varepsilon=\varepsilon_n\to+0$ is discrete. The main object of our
investigation is the GSD $A_{\nu6}^k(t_6)$ $\nu=0,1,t$ and $k\in\mathbb{Z}_+$, which is defined via
the recurrence procedure:
\begin{equation}
 \label{eq:II-A-nu6-k-rec}
A_{\nu6}^0(t_6)=A_{\nu6}(t_6),\qquad A_{\nu6}^{k+1}(t_6)={\tilde A}_{\nu6}^k(t_6),
\end{equation}
where $A_{\nu6}(t_6)$ is the solution of System~(\ref{eq:schlesinger-p6}) corresponding to
${\cal M}_6(\Theta_{06},\Theta_{16},\Theta_{t6},\Theta_{\infty6})$, and ${\tilde A}_{\nu6}^k(t_6)$
is obtained via (\ref{eq:Atilde-nu})--(\ref{eq:Atilde-nu-result}) (with $J_{0t}^{+-}$) by
inserting $A_{\nu6}^k(t_6)$ into these formulae instead of $A_{\nu6}$.

We interpret the formal limit {\bf II} as the asymptotics as $n\to+\infty$  of the sequence
$A_{\nu6}^{2n}(\varepsilon_nt_5)$, or, equivalently, as the asymptotics of the corresponding
sequences of their matrix elements. Here $t_5$ is a real positive number, which is further assumed
as a given parameter. Asymptotics of the odd sequence, $A_{\nu6}^{2n+1}(\varepsilon_nt_5)$, is
given by exactly the same formulae as for the even one, $A_{\nu6}^{2n}(\varepsilon_nt_5)$, (see
Equations~(\ref{eq:13b-Theta})-(\ref{eq:13b-s16}) with $\varepsilon\to\varepsilon_n$) but in the
corresponding formulae for the monodromy data one has to change $\Theta_6\to\Theta_6-1$. In the
derivation of the results stated in Theorem~\ref{th:2} below, we omit the super/subscript $n$,
if it does not cause a confusion.

To formulate our result we need to introduce some preliminary notation. It follows from Equations
(\ref{eq:monodromy-local}), (\ref{eq:II-theta-t6}), and (\ref{eq:II-def-theta-infty5}) that
$e^{\pm\pi i\Theta_6}$ are the eigenvalues of $M_{t6}$ and
$e^{\pm\pi i(\Theta_{\infty5}-\Theta_6)}$ are the eigenvalues of $M_{06}$. Define
\begin{gather}
 \label{eq:II-T}
{\bf T}={\rm tr}(M_{t6}-e^{+\pi i\Theta_6})(M_{06}-e^{-\pi i(\Theta_{\infty5}-\Theta_6)}),\\
 \label{eq:II-l}
l=\frac1{\pi i}\ln\left(\cos(\pi\Theta_{\infty5})-\frac12{\bf T}
\pm\sqrt{\Big(\cos(\pi\Theta_{\infty5})-\frac12{\bf T}\Big)^2-1}\right),\qquad-1<{\rm Re}\,l<1,\\
 \label{eq:II-alpha-beta}
\alpha=-\frac12(\Theta_{\infty5}-l),\qquad\beta=-\frac12(\Theta_{\infty5}+l).
\end{gather}
The sign before the square root in Equation~(\ref{eq:II-l}) can be chosen arbitrary: the change of
the sign means simply the change $\alpha\leftrightarrow\beta$ in our construction, which is invariant
under this transformation.
\begin{theorem}
 \label{th:2}
Let $\mu_6\in{\cal M}_6(\Theta_{\infty5}-\Theta_6,\Theta_{16},\Theta_6,\Theta_{\infty6})$. Suppose
the following conditions are valid:\\
{\rm1. (\ref{eq:II-theta-t6}), (\ref{eq:II-def-theta-infty5})}, and
$\Theta_{\nu6}\in\mathbb{C\setminus Z}$, $\nu=0,1,\infty$;\\
{\rm2}. $\alpha,\beta,l\in\mathbb{C\setminus Z}$;\\
{\rm3}. $(M_{t6}-e^{-\pi i\Theta_6})(M_{06}-e^{-\pi i(\Theta_{\infty5}-\Theta_6)})\neq0$.\\
{\rm4}. The inverse monodromy problem for Equation~{\rm(\ref{eq:linear-4-regular})} is solvable
for all pairs $(\mu_6,t_6)$ such that $t_6=\varepsilon_nt_5>0$, where $\varepsilon_n$ is given by
Equation~{\rm(\ref{eq:II-theta-t6})}.\\
{\rm5}. It is possible to define the sequence $A_{\nu6}^{2n}(\varepsilon_nt_5)$: it is true in
particular, if $\mu_6$ corresponds to the non-classical solution of Equation~{\rm(\ref{eq:p6})}
\footnote[8]{The last condition is not necessary: such sequence is possible to organize for the
classical solutions too, just for some of them $n\in\mathbb{Z}_-$ instead of $\mathbb{Z}_+$;
we do not discuss here the corresponding modifications.}.

Then define $\mu_5\in{\cal M}_5(\Theta_{05},\Theta_{15},\Theta_{\infty5})$:
\begin{gather}
 \label{eq:II-P5-theta}
\Theta_{05}=\Theta_{\infty6},\qquad\Theta_{15}=\Theta_{16},\\
 \label{eq:II-P5-Stokes}
S_0=\left(\begin{array}{cc}
1&0\\
-\frac{2\pi i}{\Gamma(1-\alpha)\Gamma(1-\beta)}&1
\end{array}\right),\qquad
S_1=\left(\begin{array}{cc}
1&-\frac{2\pi i}{\Gamma(\alpha)\Gamma(\beta)}e^{\pi i\Theta_{\infty5}}\\
0&1
\end{array}\right),\\
 \label{eq:P5-M05-M15}
KM_{\infty6}K^{-1}=M_{05},\qquad
KM_{16}K^{-1}=M_{15},
\end{gather}
where $K$ is the unique (up to the sign) solution of the system:
\begin{equation}
 \label{eq:II-P5-K-system}
KM_{t6}K^{-1}=S_0e^{\pi i\Theta_6\sigma_3},\qquad
KM_{06}K^{-1}=e^{-\pi i\Theta_6\sigma_3}S_1e^{\pi i\Theta_{\infty5}\sigma_3}.
\end{equation}
See Theorem~{\rm\ref{th:A}} in the Appendix and Equations~{\rm((\ref{eq:A.19})--(\ref{eq:A.21})}.

If for a given pair $(\mu_5,t_5)$ the inverse monodromy problem is solvable, then an
asymptotic expansion of the GSD $A_{\nu6}^{2n}(\varepsilon_nt_5)$ as $\varepsilon_n\to+0$ is
given by Equations {\rm(\ref{eq:lim2-t-lambda})--(\ref{eq:lim2-Rt6-derivative})} with
$\varepsilon=\varepsilon_n$.
\end{theorem}
\begin{derivation}
Suppose that the matrices $A_{\nu6}$ in Equation~(\ref{eq:linear-4-regular})
satisfy Equations~(\ref{eq:lim2-t-lambda})--(\ref{eq:lim2-Rt6-derivative}) corresponding to the
formal limit {\bf II}. Then one proves that in the domain
$\lambda_6\in\Omega_{out}:=\{\lambda_6\in\mathbb{C}, \varepsilon/|\lambda_6|^3
\leq{\cal O}(\varepsilon^{3\delta_0-2}),\;
|\arg\,\lambda_6+\frac{\pi}2|<\delta_1,\;\delta_0>\frac23,\;0\leq\delta_1<\pi\}$ the $\Psi_6$
function has the following asymptotics
\begin{equation}
 \label{eq:psi6}
\Psi_6(\lambda_6,t_6)K^{-1}=\varepsilon^{\sigma_3}R\Psi_5^0(\lambda_5,t_5)(I+o(1)),
\end{equation}
where $K,R\in{\rm SL}(2,\mathbb{C})$, $R=\varepsilon^{-\sigma_3}R_{t6}+o(1)$,
$\varepsilon^{-\sigma_3}R_{t6}={\cal O}(1)$. The matrix $K$ in Equation~(\ref{eq:psi6}) is
independent of $t_\nu$, $\lambda_\nu$, $\nu=5,6$ and $\varepsilon$. Equations~(\ref{eq:P5-M05-M15})
are an immediate consequence of (\ref{eq:psi6}) and (\ref{eq:lim2-t-lambda}). To this end our
problem is to find $K$. To solve it we consider the function $\Phi$ which solves the hypergeometric
equation:
\begin{equation}
 \label{eq:II-Phi-hyper}
\frac{d\Phi}{d\lambda_6}=\left(\frac{A_{06}}{\lambda_6}+\frac{A_{16}}{\lambda_6-t_6}\right)\Phi,
\end{equation}
and has the same monodromy matrices at regular singularities $\lambda_6=0$ and $t$ as
the function $\Psi_6(\lambda_6,t_6)K^{-1}$. We call the latter property of $\Phi$ as the
condition $M$. The fundamental solution of Equation ~(\ref{eq:II-Phi-hyper}) can be
presented as
\begin{equation}
 \label{eq:II-Phi-Y}
\Phi=PY(x),\qquad x=\lambda_6/t_6,
\end{equation}
where $Y(x)$ is given by (\ref{eq:Y(x)-def}) with the parameters $\alpha$, $\beta$, $\gamma$
satisfying the equations
\begin{equation}
 \label{eq:II-alpha-beta-gamma}
1-\gamma=\Theta_{06},\quad\gamma-1-\alpha-\beta=\Theta_{t6}, 0<|\beta-\alpha||<1,
\end{equation}
and $P^{-1}A_{\nu6}P\approx A_{\frac\nu{t}Y}$, where $A_{\frac\nu{t}Y}$ is defined by
(\ref{eq:A0Y}) and (\ref{eq:A1Y}). The last condition in (\ref{eq:II-alpha-beta-gamma}) is assumed
by taking into account the theorem formulated in \cite{J}[\S2, pp.1145-1146].
Now we have to determine $l=\alpha-\beta$. To do this we use condition $M$ in the following way.
First, for $m=0,1$ we calculate the matrix $\hat K(2m-1)$ by help of the equation
\begin{equation}
 \label{eq:hat-K-calc}
\Phi\hat K^{-1}(2m-1)=\varepsilon^{\sigma_3}R\Psi_5^m,
\qquad\lambda_6\sim{\cal O}(\varepsilon^{1-\delta_0}).
\end{equation}
Thus, using definitions of Section~\ref{sec:2} we find the Stokes multiplier
\begin{equation}
 \label{eq:II-S0}
S_0=\hat K(-1)\hat K^{-1}(1),
\end{equation}
and the monodromy matrix at the infinity point
\begin{equation}
 \label{eq;II-M-infty}
M_{\infty5}=\hat K(-1)M_{1Y}M_{0Y}\hat K^{-1}(-1).
\end{equation}
Now the direct calculation shows
\begin{equation}
 \label{eq:hat-K-M1Y-S0}
\hat K(-1)M_{1Y}\hat K^{-1}(-1)=S_0e^{\pi i\Theta_{t6}\sigma_3}.
\end{equation}
Comparing Equation~(\ref{eq:M-infty-stokes-k}) for $k=0$ with Equations~(\ref{eq;II-M-infty}) and
(\ref{eq:hat-K-M1Y-S0}) we obtain the other Stokes multiplier:
\begin{equation}
 \label{eq:II-S1-calc}
\hat K(-1)M_{0Y}\hat K^{-1}(-1)=e^{\pi i\Theta_{06}\sigma_3}e^{-\pi i\Theta_{\infty5}\sigma_3}
S_1e^{\pi i\Theta_{\infty5}\sigma_3}.
\end{equation}
Equations (\ref{eq:II-S0}) and (\ref{eq:II-S1-calc}) are equivalent to the ones presented in
(\ref{eq:II-P5-Stokes}), but the values of the parameters $\alpha$ and $\beta$ still remain
undetermined. Comparing Equations (\ref{eq:hat-K-calc}), (\ref{eq:II-Phi-Y}), and (\ref{eq:psi6}),
one finds
\begin{equation}
 \label{eq:II-M-nu6-M-nu-Y}
KM_{\nu6}K^{-1}=\hat K(-1)M_{\frac\nu{t}Y}\hat K^{-1}(-1),\qquad\nu=0,t.
\end{equation}
Now Equations (\ref{eq:II-M-nu6-M-nu-Y}), (\ref{eq:II-S1-calc}), and (\ref{eq:hat-K-M1Y-S0}) yield
(\ref{eq:II-P5-K-system}), which are the {\it exact} formulae, since both sides of these equations
are independent of $\varepsilon$. Finally we use Theorem~\ref{th:A} of the Appendix to find
formulae (\ref{eq:II-T}) -- (\ref{eq:II-alpha-beta}) and the matrix $K$.

The rigorous aspect of the above calculation is based on the justification scheme suggested in
\cite{K4}. The scheme can be explained as follows. Suppose that the inverse monodromy problem
for Equation ~(\ref{eq:linear-5}) is solvable for a given value of the parameter $t_5$ in some
neighborhood ${\cal O}(\mu_5)$ of the given point
$\mu_5\in{\cal M}_5(\Theta_{05},\Theta_{15},\Theta_{\infty5})$. Thus for all
$\tilde\mu_5\in{\cal O}(\mu_5)$ we have functions $\tilde u_5$, $\tilde z_5$, and $\tilde y_5$ whose
monodromy data coincide with $\tilde\mu_5$. Then using Equations
(\ref{eq:II-zt6-z06}) --(\ref{eq:II-ut6zt6st6}), and (\ref{eq:lim2-t-lambda}) with
$\tilde u_5$, $\tilde z_5$, and $\tilde y_5$ instead of $u_5$, $z_5$, and $y_5$ and the discrete
$\varepsilon=\varepsilon_n$, substitute them into the matrix elements of the residue matrices of
Equation~(\ref{eq:linear-4-regular}). Then the derivation presented above can be interpreted as
as the proof that the monodromy data $\hat\mu_5$, of thus obtained
Equation~(\ref{eq:linear-4-regular}), differ from the data $\tilde\mu_5$, obtained by the inversion of Equations~(\ref{eq:P5-M05-M15})
and (\ref{eq:II-P5-K-system}), on quantities estimated as $\parallel\tilde\mu-\hat\mu\parallel<o(1)$.
To be able to apply the scheme of justification suggested in \cite{K4} the last estimate must possess
an important property of ``local uniformness'': this means that there exists a neighborhood of
$\tilde\mu_5$ such that the above estimate is uniform for all points $\hat\mu_5$ in this neighborhood.
In the terminology of the work \cite{K4}: the inverse monodromy problem for
Equation~(\ref{eq:linear-4-regular}) for the monodromy data $\tilde\mu_5$ is asymptotically solvable.
The main result of \cite{K4} says that if the monodromy problem is asymptotically solvable, then
it is exactly solvable. Since the solution of the inverse monodromy problem is unique, this solution
coincides with the GSD introduced in the beginning of this section.

Now I am going to give some details to the calculation outlined above. Let us begin with the matching
(\ref{eq:hat-K-calc}) as it is the most crucial point. Consider Equation~(\ref{eq:II-Phi-Y}).
The matrix $P$ which maps the function $Y(x)$ into a fundamental solution of the hypergeometric
equation (\ref{eq:II-Phi-hyper}) is defined as follows
\begin{equation}
 \label{eq:II-P-details}
P^{-1}A_{\nu6}P=A_{\frac\nu{t}Y}+{\cal O}\big(\varepsilon^{\delta_0-1}\big),\qquad\nu=0,1,
\qquad P\in{\rm SL}(2,\mathbb{C}).
\end{equation}
Let's explain the estimate ${\cal O}\big(\varepsilon^{\delta_0-1}\big)$ in
Equation (\ref{eq:II-P-details}): it is an important and rather subtle moment. First we show that
the estimate cannot be less than ${\cal O}(1)$. Actually, summing up equations in (\ref{eq:II-P-details})
(with the error estimate ${\cal O}(1)$) and using the second equality in (\ref{eq:hypergeo-Y(x)})
and the first equation in (\ref{eq:niormalization-for-P6}), one arrives at
\begin{equation}
 \label{eq:II-P-details-2}
P^{-1}\left(\frac{\Theta_{\infty6}}2\sigma_3+A_{16}\right)P=-\frac{\beta-\alpha}2\sigma_3+{\cal O}(1).
\end{equation}
Suppose that in Equation~(\ref{eq:II-P-details-2}) one have a better estimate, $o(1)$, instead of
${\cal O}(1)$. Then it implies
\begin{equation}
 \label{eq:II-det-details}
\det\left(\frac{\Theta_{\infty6}}2\sigma_3+A_{16}\pm\frac{\beta-\alpha}2\right)=o(1).
\end{equation}
The last equation leads to the nontrivial dependence of $\alpha$ and $\beta$ from $t_5$ (see
Equations~(\ref{eq:II-z16}) and the first equation (\ref{eq:A-and-R}) with $\nu=1$.
This dependence contradicts the isomonodromy condition, since the monodromy data in particular
Stokes multipliers $S_k$ must be independent of $t_5$ (see (\ref{eq:II-P5-Stokes})). Thus the
minimal possible error estimate in Equation~(\ref{eq:II-P-details}) is ${\cal O}(1)$,
i.e., $\delta_0\leq1$. In the last case Equation~(\ref{eq:II-P-details-2}) becomes uninformative.

One of the most important moments of our derivation is the matching procedure of $\Phi$ with
$\varepsilon^{\sigma_3}R_{t6}\Psi_5^0$. This matching occurs at the domain
$\lambda_6\sim\varepsilon^{1-\delta_0}$, $\frac23<\delta_0<1$, where the lower limit for $\delta_0$
is defined via the matching of $\Psi_6$ and $R_{t6}\Psi_5^0$ (see Equation~(\ref{eq:psi6})):
The term $t_6^2A_{t6}/\lambda_6^3\sim{\cal O}\big(\varepsilon^{3\delta_0-2}\big)$, which was neglected
in Equation~(\ref{eq:linear-4-regular}) to obtain the equation for $R_{t6}\Psi_5^0$, should decrease.
The upper limit for $\delta_0:\;\delta_0<1$ is originated from the matching of $\Psi_6$ and $\Phi$:
the corresponding equations differ by the term
${\cal O}(t_6A_{t6}/\lambda_6^2)\sim{\cal O}(\varepsilon^{2(\delta_0-1)})$. Thus to guarantee
the matching we must define $P$ in Equation~(\ref{eq:II-P-details}) with the error estimate not
worse than $\lambda_6{\cal O}(\varepsilon^{2(\delta_0-1)})={\cal O}(\varepsilon^{\delta_0-1})$.
To summarize: the primary qualities of $\Phi$ (according to \cite{K4}) are its monodromy data and
the matching with $\varepsilon^{\sigma_3}R_{t6}\Psi_5^0$. It is because of these properties the
function $\Phi$, in fact, is not the exact solution of Equation~(\ref{eq:II-Phi-hyper}), but solves
it up to the leading order with the rather large error ${\cal O}(\varepsilon^{\delta_0-1})$.
It seems that the function $\Phi$ with the required properties can be defined to satisfy
Equation~(\ref{eq:II-Phi-hyper}) more precisely, with the error ${\cal O}(1)$, but for this goal
we need to substantially increase complexity of our calculations. In any case we cannot define the
function $\Phi$ to satisfy Equation~(\ref{eq:II-Phi-hyper}) with the error less than ${\cal O}(1)$,
because it contradicts the matching procedure. Actually, the asymptotic expansion for $A_{16}$,
which we substitute into Equation~(\ref{eq:II-det-details}) to get the contradiction, is obtained
from the differential equation for $R_{t6}\Psi_5^0$ and, hence, it is inexplicit consequence of the
matching.

Taking into account the discussion above, we rewrite System~(\ref{eq:II-P-details}) with
the error ${\cal O}(\varepsilon^{\delta_0-1})$ changed by ${\cal O}(1)$ as follows:
\begin{equation}
 \label{eq:II-P-calc-elements}
\begin{aligned}
(-1)^{\frac{\nu}t}\frac{\beta(1-\gamma)}{\beta-\alpha}&=
\left(d-\frac{b\varepsilon^2}{u_{\nu6}}\right)
\left(b\Theta_{\nu6}-\frac{u_{\nu6}z_{\nu6}}{\varepsilon^2}
\left(d-\frac{b\varepsilon^2}{u_{\nu6}}\right)\right)+{\cal O}(1),\\
(-1)^{1+\frac{\nu}t}\frac{\alpha(1-\gamma)}{\beta-\alpha}&=
\left(c-\frac{a\varepsilon^2}{u_{\nu6}}\right)
\left(\frac{u_{\nu6}z_{\nu6}}{\varepsilon^2}
\left(c-\frac{a\varepsilon^2}{u_{\nu6}}\right)
-a\Theta_{\nu6}\right)+{\cal O}(1),\\
(-1)^{\frac{\nu}t}\frac{(\alpha+\beta)(1-\gamma)}{2(\beta-\alpha)}&=
\left(d-\frac{b\varepsilon^2}{u_{\nu6}}\right)
\left(c-\frac{a\varepsilon^2}{u_{\nu6}}\right)
\frac{u_{\nu6}z_{\nu6}}{\varepsilon^2}\\
&-\frac{\Theta_{\nu6}}2\left(a\Big(d-\frac{b\varepsilon^2}{u_{\nu6}}\Big)
+b\Big(c-\frac{a\varepsilon^2}{u_{\nu6}}\Big)\right)+{\cal O}(1),
\end{aligned}
\end{equation}
where $\nu=0,t$ and $a$, $b$, $c$, and $d$ are the matrix elements of $P$:
\begin{gather}
 \label{eq:P-abcd}
P=\varepsilon^{\sigma_3}\left(\begin{array}{cc}
a&b\\
c&d
\end{array}\right),\qquad\text{with}\;\;a,b,c,d={\cal O}(1),\\
 \label{eq:ad-bc=1}
ad-bc=1.
\end{gather}
Using now Formulae~(\ref{eq:II-zt6-z06}), (\ref{eq:II-ut6-u06}), (\ref{eq:II-alpha-beta-gamma}), and
(\ref{eq:P-abcd}), one rewrites Equations~(\ref{eq:II-P-calc-elements}) as follows:
\begin{equation}
 \label{eq:u-nu6-z-nu6}
u_{\nu6}=\varepsilon^2\frac{a+b}{c+d}+{\cal O}(\varepsilon^3),\quad
\frac{z_{\nu6}}{\Theta_{\nu6}}=\frac{(a+b)(c+d)}2
\left(\frac{b-a}{b+a}-\frac{\beta-\alpha}{\beta+\alpha}\right)\!+{\cal O}(\varepsilon),\;\;\nu=0,t.
\end{equation}
To determine the matrix $P$ uniquely we must add to (\ref{eq:ad-bc=1}) and (\ref{eq:u-nu6-z-nu6})
one more equation. We get it via the matching (\ref{eq:hat-K-calc}): the matching domain now is
\begin{equation}
 \label{eq:II-matching-domain}
x=\frac{\lambda_6}{t_6}=(\varepsilon\lambda_5t_5)^{-1}={\cal O}\big(\varepsilon^{-\delta_0}\big),\qquad
\frac23<\delta_0<1,
\end{equation}
so that Asymptotics~(\ref{eq:Y(x)-spec-asympt}), (\ref{eq:x-spec-asympt-cond}) imply the following
equations:
\begin{gather}
 \label{eq:varkappa}
\varkappa=2m-1,\qquad\text{for}\quad m=0,1\quad\alpha+\beta=-\Theta_{\infty5},\\
 \label{eq:R-t6-abcd}
\frac{\varepsilon^{\sigma_3}}{\sqrt{\alpha-\beta}}\left(\begin{array}{cc}
a&b\\
c&d
\end{array}\right)
\left(\begin{array}{cc}
1&\beta\\
1&\alpha
\end{array}\right)
t^{-\frac{\Theta_{\infty5}}2\sigma_3}=R_{t6}+\varepsilon^{\sigma_3}
{\cal O}\big(\varepsilon^{\tilde\delta}\big).
\end{gather}
Recall the definition of ${\cal O}\big(\varepsilon^{\tilde\delta}\big)$ below
Equation~(\ref{eq:x-spec-asympt-cond}). The first equation in (\ref{eq:varkappa}) and
(\ref{eq:hat-K-calc}) show that (\ref{eq:II-S0}) and (\ref{eq;II-M-infty}) hold up to the
error bound ${\cal O}\big(\varepsilon^{\tilde\delta}\big)$. From (\ref{eq:R-t6-abcd}) one finds
that Equations~(\ref{eq:ad-bc=1}) and (\ref{eq:u-nu6-z-nu6}) are valid up to the estimate
${\cal O}\big(\varepsilon^{\tilde\delta}\big)$. Equation (\ref{eq:R-t6-abcd}) together with
the first equation in (\ref{eq:II-ut6-u06}) and the equality $\varepsilon_n\Theta_{t6}=-1$ we
arrive at:
\begin{gather*}
a=-\frac1{s_{t6}\sqrt{\alpha-\beta}}\Big(\alpha t_5^{\frac{\Theta_{\infty5}}2}
+u_5t_5^{-\frac{\Theta_{\infty5}}2}\Big),\qquad
b=\frac1{s_{t6}\sqrt{\alpha-\beta}}\Big(\beta t_5^{\frac{\Theta_{\infty5}}2}
+u_5t_5^{-\frac{\Theta_{\infty5}}2}\Big),\\
c=\frac{s_{t6}}{\sqrt{\alpha-\beta}}\left(\frac\alpha{u_5}
\frac{z_5}{\Theta_{05}}t_5^{\frac{\Theta_{\infty5}}2}+\Big(\frac{z_5}{\Theta_{05}}+1\Big)
t_5^{-\frac{\Theta_{\infty5}}2}\right),\\
d=-\frac{s_{t6}}{\sqrt{\alpha-\beta}}\left(\frac\beta{u_5}
\frac{z_5}{\Theta_{05}}t_5^{\frac{\Theta_{\infty5}}2}+\Big(\frac{z_5}{\Theta_{05}}+1\Big)
t_5^{-\frac{\Theta_{\infty5}}2}\right).
\end{gather*}
Now one uses (\ref{eq:II-S0}) and (\ref{eq:varkappa}) to find, again up to
${\cal O}\big(\varepsilon^{\tilde\delta}\big)$, the first equation in
(\ref{eq:II-P5-Stokes}). To finish the proof we have to confirm the second equation in
(\ref{eq:II-P5-Stokes}). To do it, we have to prove (\ref{eq:hat-K-M1Y-S0}) and to find
$S_1$ from (\ref{eq:II-S1-calc}). Details for these calculations are as follows: apply
the formulae for the $\Gamma$-function (\ref{eq:Gamma-ratio}) and (\ref{eq:Gamma-product})
to Equations (\ref{eq:C0-alpha-beta-gamma}) and (\ref{eq:C1-alpha-beta-gamma}) to find:
\begin{gather*}
C_{\alpha\beta\gamma}^0{\hat K}^{-1}(-1)=D_0\left(1+{\cal O}\big(\varepsilon^{\tilde\delta}\big)\!\right)\!
\left(\begin{array}{cc}
1&-\frac{\pi\exp(\pi(\alpha+\beta+1-\gamma))}{(\alpha-\beta)\sin(\pi\gamma)\Gamma(\alpha)\Gamma(\beta)}\\
0&1
\end{array}\right)\!
\frac1{\sqrt{\alpha-\beta}^{\sigma_3}}\left(1+{\cal O}\big(\varepsilon^{\tilde\delta}\big)\!\right),\\
C_{\alpha\beta\gamma}^1{\hat K}^{-1}(-1)=D_1\left(1+{\cal O}\big(\varepsilon^{\tilde\delta}\big)\!\right)\!
\left(\begin{array}{cc}
1&0\\
\frac{\pi\exp(\pi(\gamma-\alpha-\beta))(\beta-\alpha)}{\sin(\pi(\alpha+\beta-\gamma))\Gamma(1-\alpha)
\Gamma(1-\beta)}&1
\end{array}\right)\!
\left(1+{\cal O}\big(\varepsilon^{\tilde\delta}\big)\!\right),
\end{gather*}
where $D_k$ are diagonal matrices. Finally, one uses the definition of $M_{\nu Y}$:
$$
M_{\nu Y}=(C_{\alpha\beta\gamma}^\nu)^{-1}e^{\pi i\Theta_{t\cdot\nu6}\sigma_3}C_{\alpha\beta\gamma}^\nu,
\qquad\nu=0,1
$$
to prove (\ref{eq:hat-K-M1Y-S0}) and (\ref{eq:II-S1-calc}).
\end{derivation}

The result of Theorem~\ref{th:2} can be immediately generalized as follows. For $k=1,\ldots, n$ define
$a_k\in\mathbb{C}\setminus\{0\}$, $a_k\neq a_l$, for $k\neq l$, $\frac{\partial a_k}{\partial t_6}=
\frac{\partial a_k}{\partial\lambda_6}=0$. Consider the functions
$\Psi_\nu=\Psi_\nu(\lambda_\nu,t_\nu;a_1,\ldots,a_n)\in{\rm SL}(2,\mathbb{C})$ $\nu=5,6$, which
generalize the functions $\Psi_\nu$ defined in Section~\ref{sec:2}, as the normalized
(\ref{eq:psi6-infty-normalization}) and (\ref{eq:Psi5-asymptotics-infty}) fundamental solutions of the
systems:
\begin{align*}
\frac{d}{d\lambda_6}\Psi_6&=\left(\frac{A_{06}}{\lambda_6}+\frac{A_{t6}}{\lambda_6-t_6}+
\sum\limits_{k=1}^{m}\frac{A_{k6}}{\lambda_6-a_k}\right)\Psi_6,\\
\frac{d}{dt_6}\Psi_6&=-\frac{A_{t6}}{\lambda_6-t_6}\Psi_6,
\end{align*}
\begin{align*}
\frac{d}{d\lambda_5}\Psi_5&=\left(\frac{t_5}2\sigma_3+\frac{A_{05}}{\lambda_5}+
\sum\limits_{k=1}^{m}\frac{A_{k5}}{\lambda_5-1/a_k}\right)\Psi_5,\\
\frac{d}{dt_5}\Psi_5&=\left(\frac{\lambda_5}2\sigma_3+\frac1{t_5}\left(\frac{\Theta_{\infty5}}2\sigma_3+
\sum\limits_{k=0}^{m}A_{k5}\right)\right)\Psi_5,
\end{align*}
$$
A_{t6}+\sum\limits_{k=0}^{m}A_{k6}\equiv-\frac{\Theta_{\infty6}}2\sigma_3,\qquad
{\rm diag}\sum\limits_{k=0}^{m}A_{k5}\equiv-\frac{\Theta_{\infty5}}2\sigma_3,\qquad
A_{k\nu}\in sl_2(\mathbb{C}).
$$
The local objects like $\Theta_{k\nu}$, $R_{k\nu}$, and $M_{k\nu}$, $k=1,\ldots,m$, are defined in the
same manner as that in Section~\ref{sec:2} for $k=1$ and $a_1=1$.
Systems of isomonodromy deformations are the compatibility conditions for the systems of the linear
ODEs written above, and the GSDs $A_{k6}^{2n}(\varepsilon_nt_5)$ can be defined as in the beginning of
this section.
\begin{corollary}
 \label{cor:1}
Let conditions of Theorem~{\rm\ref{th:2}} be valid for
$$
{\cal M}_6={\cal M}_6(\Theta_{06},\Theta_{t6},\Theta_{16},\ldots,\Theta_{m6},\Theta_{\infty6}),
$$
then asymptotics ($n\to+\infty$) of the GSDs $A_{k6}^{2n}(\varepsilon_nt_5)$, $(k=1,\ldots,m)$, are
given by Equations~{\rm(\ref{eq:lim2-t-lambda}) -- (\ref{eq:lim2-A05-A15})} which are supplemented
by
\begin{equation}
 \label{eq:II-A-k6-Ak5}
 R_{t6}^{-1}A_{k6}R_{t6}=A_{k5}+{\cal O}(\varepsilon),
\end{equation}
Equation~{\rm(\ref{eq:lim2-Rt6-derivative})} should be generalized as follows
\begin{equation}
 \label{eq:II-R-t6-derivative-generalized}
-R_{t6}^{-1}\frac{d}{dt_5}R_{t6}=\frac1{t_5}\left(\frac{\Theta_{\infty5}}2\sigma_3+
\sum\limits_{k=0}^mA_{k5}\right)+{\cal O}(\varepsilon),
\end{equation}
and with the change $\varepsilon\to\varepsilon_n$ in all formulae. The monodromy data for the
limiting $\Psi_5$ function is given by exactly the same formulae as that in Theorem~{\rm\ref{th:2}}
supplemented with the equations:
\begin{equation}
 \label{eq:II-theta-k6-M-k6}
\Theta_{k6}=\Theta_{k5},\qquad
KM_{k6}K^{-1}=M_{k5},\qquad
k=1,\ldots,m.
\end{equation}
\end{corollary}
\begin{proposition}
 \label{prop:1}
For the sixth Painlev\'e function (and the corresponding $\tau$-function) both limits considered
in Section~{\rm\ref{sec:3}} are equivalent under the definition given in the Introduction.
\end{proposition}
\begin{proof}
Let us supply with the superscripts $I$ and $II$ the objects corresponding to the first/second
limits. Then for $\Psi_6$ function one finds:
\begin{equation}
 \label{eq:Psi-I-II-equivalence}
\Psi_6^{II}(\lambda_6^{II},t_6)=t_6^{-\frac{\Theta_{05}}2\sigma_3}(R_{06}^I)^{-1}
\Psi_6^I(\lambda_6^I,t_6)(C_{06}^I)^{-1},\qquad\lambda_6^{II}=t_6/\lambda_6^I.
\end{equation}
This transformation is controlling the limiting procedure. Thus it is not necessary to consider
especially what is happening on the lower ${\cal M}$-plane of the diagram on Figure~3.
Formula~(\ref{eq:Psi-I-II-equivalence}) generates the following transformation for the matrices
$A_{\nu6}$:
\begin{equation}
 \label{eq:A-I-II-equivalence}
A_{\nu6}^{II}=t_6^{-\frac{\Theta_{06}^I}2\sigma_3}(R_{06}^I)^{-1}A_{\frac{t}\nu6}^IR_{06}^I
t_6^{\frac{\Theta_{06}^I}2\sigma_3},\qquad\nu=1,t.
\end{equation}
For the formal monodromy one finds:
\begin{gather*}
\Theta_{05}^{II}=\Theta_{\infty6}^{II}=\Theta_{06}^I=\Theta_{05}^I,\\
\Theta_{15}^{II}=\Theta_{16}^{II}=\Theta_{t6}^I=-\Theta_{15}^I,\\
\Theta_{t6}^{II}=\Theta_{16}^I=-\frac1{\varepsilon},\\
\Theta_{06}^{II}=\Theta_{\infty6}^I\qquad\Theta_{\infty6}^I+\Theta_{16}^I=
\Theta_{06}^{II}+\Theta_{t6}^{II}=\Theta_{\infty5}.
\end{gather*}
Substituting~(\ref{eq:lim1-Psi}) and (\ref{eq:lim2-Psi}) into Equation~(\ref{eq:Psi-I-II-equivalence})
one obtains
$$
R_{t6}^{II}\Psi_5^I(\lambda_5,t_5)K^{II}=t_6^{-\frac{\Theta_{06}^I}2\sigma_3}(R_{06}^I)^{-1}
R_{16}^I\Psi_5^I(\lambda_5,t_5)K^I(C_{06}^I)^{-1}.
$$
Where $K^p$ is the matrix introduced in Theorem $p$ ($p=I,II$). Thus we have
\begin{gather}
 \label{eq:R-relation}
R_{t6}^{II}f_0^{-\sigma_3}=t^{-\frac{\Theta_{06}^I}2\sigma_3}(R_{06}^I)^{-1}R_{16}^I,\\
 \label{eq:K-relation}
K^{II}=f_0^{-\sigma_3}K^I(C_{06}^I)^{-1}.
\end{gather}
Using Relations (\ref{eq:R-relation}) and (\ref{eq:K-relation}) one proves that
Equations~(\ref{eq:lim1-A}) and (\ref{eq:lim2-A05-A15}) are equivalent. For the matrix elements  one
has to take into account that $\Theta_{15}^{II}=-\Theta_{15}^I$ and, hence,
Equations~(\ref{eq:A05-parametrization}) and (\ref{eq:A15-parametrization}) yield
$$
u_5^{II}=u_5^I,\;\;
z_5^{II}=z_5^I,\;\;
y_5^{II}\!\left(z_5+\frac{\Theta_{05}-\Theta{15}^{II}+\Theta_{\infty5}}2\right)\!=
y_5^I\!\left(z_5+\frac{\Theta_{05}-\Theta{15}^I+\Theta_{\infty5}}2\right)\!.
$$
The so-called $c$-problem discussed in the beginning of this Section reveals itself in
Equations~(\ref{eq:R-relation}) and (\ref{eq:K-relation}) as the ambiguity in the definitions of
$R_{06}$, $C_{06}$:
$$
R_{06}\to R_{06}c^{\sigma_3},\qquad
C_{06}\to c^{-\sigma_3}C_{06}.
$$
\end{proof}
\appendix
\section{Appendix}
For $\nu=0,1$ consider the matrices
$$
M_{\nu6}=\left(\begin{array}{cc}
a_\nu&b_\nu\\
c_\nu&d_\nu
\end{array}\right)
\in{\rm SL}(2,\mathbb{C})
$$
with the eigenvalues $r_\nu\in\mathbb{C}\setminus\{0,\pm1\}$. The numbers $1/r_\nu$ have the same
properties. For the matrix elements one finds
\begin{equation}
 \label{eq:A.1}
b_\nu c_\nu=-(a_\nu-r_\nu)(a_\nu-1/r_\nu).
\end{equation}
Define the upper- and lower-triangular matrices:
$$
M_{0\Delta}=\left(\begin{array}{cc}
r_0&f_0\\
0&1/r_0
\end{array}\right),\qquad
M_{1\Delta}=\left(\begin{array}{cc}
r_1&0\\
f_1&1/r_1
\end{array}\right),
$$
where $f_\nu\in\mathbb{C}$. Here we study the following problem.\\
{\it Problem}: For given matrices $M_{\nu6}$ find the parameters $f_\nu$
such that there exists the matrix $K\in{\rm SL}(2,\mathbb{C})$ which solves
\begin{equation}
 \label{eq:A.2}
KM_{\nu6}=M_{\nu\Delta}K,\qquad\nu=0,1.
\end{equation}
The solution of the problem is formulated as Theorem~\ref{th:A} below.
\begin{proposition}
 \label{prop:A1}
The problem is unsolvable if the matrices $M_{\nu6}$ satisfy one of the following equations:
\begin{gather}
 \label{eq:A.3}
(M_{06}-r_0)(M_{16}-r_1)=0,\\
 \label{eq:A.4}
(M_{16}-1/r_1)(M_{06}-1/r_0)=0.
\end{gather}
\end{proposition}
\begin{proof}
Suppose the problem is solvable and one of the equations, say, (\ref{eq:A.3}) holds.
Then there exists the matrix $K$ satisfying Equation~(\ref{eq:A.2}). Using it one finds;
\begin{gather}
0=K0K^{-1}=K(M_{06}-r_0)(M_{16}-r_1)K^{-1}=(M_{0\Delta}-r_0)(M_{1\Delta}-r_1)\\
=\left(\begin{array}{cc}
0&f_0\\
0&1/r_0-r_0
\end{array}\right)
\left(\begin{array}{cc}
0&0\\
f_1&1/r_1-r_1
\end{array}\right)\\
=\left(\begin{array}{cc}
f_0f_1&-f_0(r_1-1/r_1)\\
-f_1(r_0-1/r_0)&(r_0-1/r_0)(r_1-1/r_1)
\end{array}\right).
\end{gather}
This equation contradicts to the condition $r_\nu-1/r_\nu\neq0$. One arrives at the same contradiction
by making analogous calculations for Equation~(\ref{eq:A.4}).
\end{proof}
\begin{remark}
 \label{rem:A1}{\rm
Equations (\ref{eq:A.3}) and (\ref{eq:A.4}) are equivalent. To prove it define the eigenvectors,
$e_{\nu\epsilon}$, $\nu=0,1$, $\epsilon=-1,1$:
$$
(M_{\nu6}-r_\nu^\epsilon)e_{\nu\epsilon}=0,\qquad\parallel e_{\nu\epsilon}\parallel_{\mathbb{C}^2}=1.
$$
Each pair $(e_{01},e_{0-1})$ and $(e_{11},e_{1-1})$ is a basis in $\mathbb{C}^2$. If
Equation~(\ref{eq:A.3}) holds, then we can choose $e_{01}=e_{1-1}$. Thus, if $x\in\mathbb{C}^2$, then
$x=\alpha e_{01}+\beta e_{0-1}$, where $\alpha,\beta\in\mathbb{C}$, and
$$
(M_{16}-1/r_1)(M_{06}-1/r_0)x=\alpha(M_{16}-1/r_1)e_{01}=\alpha(M_{16}-1/r_1)e_{1-1}=0.
$$
Similarly vice versa, from (A.4) to (A.3).
}\end{remark}
\begin{proposition}
 \label{prop:A2}
The problem is solvable with $f_0=f_1=0$ iff the matrices $M_{\nu6}$ are commuting,
\begin{equation}
 \label{eq:A.7}
M_{06}M_{16}=M_{16}M_{06},
\end{equation}
and don't satisfy Equations {\rm(\ref{eq:A.3})} and {\rm(\ref{eq:A.3})}. The general solution of
Equation~{\rm(\ref{eq:A.2})} can be presented as
$$
K=k^{\sigma_3}K_0,\qquad k\in\mathbb{C}\setminus\{0\},\qquad K_0^{-1}=(e_{01},e_{0-1})
$$
where $e_{0\epsilon}$ are the column eigenvectors of $M_{06}$ defined in Remark~{\rm\ref{rem:A1}}.
\end{proposition}
\begin{proof}
Multiplying l.-h.s. of Equation~(\ref{eq:A.7}) by $K_0$ and the r.-h.s. by $K_0^{-1}$ and taking into
account that $K_0M_{06}K^{-1}=r_0^{\sigma_3}$ one finds that $K_0M_{16}K^{-1}$ is a diagonal matrix
and, hence, $K_0M_{16}K^{-1}=r_1^{\delta\sigma_3}$, where $\delta=\pm1$. Thus, if $\delta=+1$, then
$k^{\sigma_3}K_0$ is a solution of Equations~(\ref{eq:A.2}) with $f_0=f_1=0$; if $\delta=-1$, then
the matrices $M_{\nu6}$ satisfy Equations~(\ref{eq:A.3}) and (\ref{eq:A.4}) and the problem is
unsolvable.

Let $L$ be a matrix that solves Equations~(\ref{eq:A.2}) for some $f_0$ and $f_1$ . Then multiplying
both parts of Equation~(\ref{eq:A.7}) by $L$ from l.-h.s. and by $L^{-1}$ from r.-h.s. and using
Equations~(\ref{eq:A.2}) one finds that $f_0=f_1=0$. Conversely, suppose that we can present
$M_{\nu6}$ as $L^{-1}r_\nu^{\sigma_3}L$ and hence the commutator $[M_{06},M_{16}]$ vanishes and
Equations~(\ref{eq:A.3}) and (\ref{eq:A.4}) are not valid.
\end{proof}
\begin{proposition}
 \label{prop:A3}
If noncommuting matrices $M_{\nu6}$ $([M_{06},M_{16}]\neq0)$ satisfy one of the following equations:
\begin{gather}
 \label{eq:A.9}
(M_{06}-r_0)(M_{16}-1/r_1)=0,\\
 \label{eq:A.10}
(M_{16}-r_1)(M_{06}-1/r_0)=0,
\end{gather}
then the problem is solvable. For all solutions:
\begin{equation}
 \label{eq:A.11}
f_1=0,\qquad f_0\in\mathbb{C}\setminus\{0\}.
\end{equation}
Conversely, if the problem is solvable with $f_1=0$, then Equations~{\rm(\ref{eq:A.9})} and
{\rm(\ref{eq:A.10})} are valid. For any pair $(f_0,f_1)$ satisfying {\rm(\ref{eq:A.11})} there are
two solutions $K$ of System~{\rm(\ref{eq:A.2})}. These solutions differ only by a choice of the sign
and can be written as follows:
\begin{equation}
 \label{eq:A.12}
K^{-1}=(e_{11},e_{1-1})\varkappa_1^{\sigma_3},\qquad\varkappa_1=\sqrt{\frac{\alpha_0(r_0-1/r_0)}{f_0}},
\end{equation}
where the column eigenvectors $e_{\nu\epsilon}$, $\nu=0,1,$ $\epsilon=1,-1,$ are defined as
$(M_{\nu6}-r_\nu^\epsilon)e_{\nu\epsilon}=0$, $\det(e_{\nu1},e_{\nu-1})=1$, $e_{11}=e_{01}$,
the last equation is always possible due to Equations~{\rm(\ref{eq:A.9})} and {\rm(\ref{eq:A.10})},
and $\alpha_0,\beta_0\in\mathbb{C}\setminus\{0\}$ are coefficients in the expansion
$e_{1-1}=\alpha_0 e_{01}+\beta_0 e_{0-1}$.
\end{proposition}
\begin{remark}
 \label{rem:A2}{\rm
Equations~(\ref{eq:A.9}) and (\ref{eq:A.10}) are equivalent. The proof can be done in the same manner
as the one in Remark~\ref{rem:A1} for the equivalence of (\ref{eq:A.3}) and (\ref{eq:A.4}).
}\end{remark}
\begin{remark}
 \label{rem:A3}{\rm
The eigenvalues $(e_{\nu1},e_{\nu-1})$ are defined up to a parameter
$$
a_\nu\in\mathbb{C}\setminus\{0\}:\quad(e_{\nu1},e_{\nu-1})\to(e_{\nu1},e_{\nu-1})a_\nu^{\sigma_3}.
$$
Under this transformation $\varkappa_\nu$ is also changing as $\varkappa_\nu\to\varkappa_\nu/a_\nu$,
so that $\varkappa_1$ in (\ref{eq:A.12}), by a proper choice of the basis, can be defined with
$\alpha_0=f_0=1$ with the only ambiguity in the branch of the square root.
}\end{remark}
\begin{proof}
Consider the following equation:
\begin{gather*}
K(M_{06}-r_0)(M_{16}-1/r_1)K^{-1}\\
=\left(\begin{array}{cc}
0&f_0\\
0&1/r_0-r_0
\end{array}\right)
\left(\begin{array}{cc}
r_1-1/r_1&0\\
f_1&0
\end{array}\right)=
\left(\begin{array}{cc}
f_0f_1&0\\
f_1(1/r_0-r_0)&0
\end{array}\right).
\end{gather*}
It proves both: if the problem is solvable and $f_1=0$, then Equation~(\ref{eq:A.9}) holds; and if
Equation~(\ref{eq:A.9}) is valid, then $f_1=0$.

Suppose (\ref{eq:A.9}) holds, define matrix $K$  by the first equation in (\ref{eq:A.12}) with a
parameter $\varkappa_1\in\mathbb{C}\setminus\{0\}$. This matrix $K$ solves Equation~(\ref{eq:A.2})
with $\nu=1$ and $f_1=0$. Substituting thus defined matrix $K$ into Equation~(\ref{eq:A.2}) with
$\nu=0$ and $f_0\in\mathbb{C}\setminus\{0\}$ one proves that it becomes an identity iff $\varkappa_1$
is given by the second equation (\ref{eq:A.12}). Note that $(f_0,f_1)\equiv(0,0)$ cannot be a solution
of the problem since $[M_{06},M_{16}]\neq0$.
\end{proof}
\begin{proposition}
 \label{prop:A4}
If noncommuting matrices $M_{\nu6}$ $([M_{06},M_{16}]\neq0)$ satisfy one of the following equations:
\begin{gather}
 \label{eq:A.13}
(M_{06}-1/r_0)(M_{16}-r_1)=0,\\
 \label{eq:A.14}
(M_{16}-1/r_1)(M_{06}-r_0)=0,
\end{gather}
then the problem is solvable. For all solutions:
\begin{equation}
 \label{eq:A.15}
f_0=0,\qquad f_1\in\mathbb{C}\setminus\{0\}.
\end{equation}
Conversely, if the problem is solvable with $f_0=0$, then Equations~{\rm(\ref{eq:A.13})} and
{\rm(\ref{eq:A.14})} are valid. For any pair $(f_0,f_1)$ satisfying {\rm(\ref{eq:A.15})} there are
two solutions $K$ of System~{\rm(\ref{eq:A.2})}. These solutions differ only by a choice of the sign
and can be written as follows:
\begin{equation}
 \label{eq:A.16}
K^{-1}=(e_{01},e_{0-1})\varkappa_0^{\sigma_3},\qquad\varkappa_0=\sqrt{\frac{f_1}{\beta_1(1/r_1-r_1)}},
\end{equation}
where the column eigenvectors $e_{\nu\epsilon}$, $\nu=0,1,$ $\epsilon=1,-1,$ are defined as
$(M_{\nu6}-r_\nu^\epsilon)e_{\nu\epsilon}=0$, $\det(e_{\nu1},e_{\nu-1})=1$, $e_{1-1}=e_{0-1}$,
the last equation is always possible due to Equations~{\rm(\ref{eq:A.13})} and {\rm(\ref{eq:A.14})},
and $\alpha_1,\beta_1\in\mathbb{C}\setminus\{0\}$ are coefficients in the expansion
$e_{01}=\alpha_1 e_{11}+\beta_1 e_{1-1}$.
\end{proposition}
\begin{remark}{\rm
 \label{rem:A4}
Equations~(\ref{eq:A.13}) and (\ref{eq:A.14}) are equivalent; the proof is analogous to the one in
Remark~\ref{rem:A1}. As for the ambiguity in Equation~(\ref{eq:A.16}), see Remark~\ref{rem:A3}.
}\end{remark}
\begin{proof}
The proof is analogous to the one for Proposition~\ref{prop:A3}. The details are left to the interested
reader.
\end{proof}
\begin{theorem}
 \label{th:A}
The solution of the problem is described by the following cases:
\begin{enumerate}
\item
If ${\rm tr}(M_{06}-1/r_0)(M_{16}-r_1)=0$ and $[M_{06},M_{16}]\neq0$, then the problem is solvable.
In this case either Equations~{\rm(\ref{eq:A.9}), (\ref{eq:A.10})} or
{\rm(\ref{eq:A.13}), (\ref{eq:A.14})} are valid and the solution is is given in
Proposition~{\rm\ref{prop:A3}} or {\rm\ref{prop:A4}}, respectively;
\item
If ${\rm tr}(M_{06}-1/r_0)(M_{16}-r_1)=0$ and $[M_{06},M_{16}]=0$, then the problem is solvable.
The solution is given in Proposition~{\rm\ref{prop:A2}};
\item
If ${\rm tr}(M_{06}-1/r_0)(M_{16}-1/r_1)=0$, then there are two subcases:
\subitem{3a.}
Equations~{\rm(\ref{eq:A.3})} and {\rm(\ref{eq:A.4})} are valid. This is unsolvable case: apply
Proposition~{\rm\ref{prop:A1}};
\subitem{3b.}
$(M_{06}-1/r_0)(M_{16}-1/r_1)=0$. The problem is solvable. All solutions are the pairs $(f_0,f_1)$,
such that $f_0,f_1\in\mathbb{C}\setminus\{0\}$ and
\begin{equation}
 \label{eq:A.17}
f_0f_1=-(r_0-1/r_0)(r_1-1/r_1).
\end{equation}
For any pair $(f_0,f_1)$ satisfying Equation {\rm(\ref{eq:A.17})}, there are two solutions, $K$ and
$-K$, of System~{\rm(\ref{eq:A.2})}. Formula~{\rm(\ref{eq:A.19})} for different choices of the
branches of the square root in {\rm(\ref{eq:A.20})} gives both, $K$ and $-K$; and
\item
If Conditions {\rm1--3} do not valid, then the problem is solvable: all solutions are the pairs
$(f_0,f_1)$ such that $f_0,f_1\in\mathbb{C}\setminus\{0\}$ and
\begin{equation}
 \label{eq:A.18}
\begin{aligned}
f_0f_1&={\rm tr}(M_1-r_1)(M_0-1/r_0)={\rm tr}(M_1-1/r_1)(M_0-r_0)\\
&={\rm tr}(M_0-1/r_0)(M_1-1/r_1)-(r_0-1/r_0)(r_1-1/r_1)\\
&=(r_0+r_1)(1/r_0+1/r_1)-\det(M_0+M_1)
\end{aligned}
\end{equation}
For a fixed pair $(f_0,f_1)$ satisfying Equation~{\rm(\ref{eq:A.18})} there are only two solutions
of System~{\rm(\ref{eq:A.2})}; both are given by the formula
\begin{equation}
 \label{eq:A.19}
K^{-1}=(\gamma_0e_{01},\gamma_1e_{1-1})
\end{equation}
and specified by the choice of a branch of the square root:
\begin{equation}
 \label{eq:A.20}
\gamma_0=\sqrt{\frac{f_1}{qs(r_1-1/r_1)}},\qquad\gamma_1=\frac1{\gamma_0s}.
\end{equation}
The column eigenvectors $e_{\nu\epsilon}$, $\nu=0,1$, $\epsilon=-1,1$ are defined as follows
$(M_{\nu6}-r_{\nu}^\epsilon)e_{\nu6}=0$ and $\det(e_{\nu1},e_{\nu-1})=1$. The complex numbers
$p$, $q$, $r$, $s$, satisfying the relation $ps-qr=1$ define transformation between the basis
$(e_{01},e_{0-1})$ and $(e_{11},e_{1-1})$ in $\mathbb{C}^2$:
\begin{equation}
 \label{eq:A.21}
\begin{aligned}
e_{11}&=pe_{01}+qe_{0-1},\\
e_{1-1}&=re_{01}+se_{0-1}.
\end{aligned}
\end{equation}
\end{enumerate}
\end{theorem}
\begin{remark}
 \label{rem:A5}{\rm
The numbers $q,s,f_1\neq0$ see the proof below. The comment concerning an ambiguity in
Equation~(\ref{eq:A.19}) is analogous to that in Remark~\ref{rem:A3}.
}\end{remark}
\begin{proof}
Consider Transformation~(\ref{eq:A.21}). If $s=0$, then one proves that Equations~(\ref{eq:A.3}) and
(\ref{eq:A.4}) are valid, so that we have the unsolvable case of the problem, see Item $3a$ in
Theorem~\ref{th:A}.

If $q=r=0$, then $[M_{06},M_{16}]=0$. It is a solvable case of the problem, see
Proposition~\ref{prop:A2} and Item $2$ of Theorem~\ref{th:A}.

If $q=0, r\neq0$, then Equations~(\ref{eq:A.9}) and (\ref{eq:A.10}) are valid: it is a solvable case
of the problem, see Proposition~\ref{prop:A3} and Item~$1$ of Theorem~\ref{th:A}.

If $q\neq0, r=0$, then Equations~(\ref{eq:A.13}) and (\ref{eq:A.14}) are valid: it is a solvable case
of the problem, see Proposition~\ref{prop:A4} and Item~$1$ of Theorem~\ref{th:A}.

If $p=0$ , then it is also a solvable case of the problem: see Item $3b$ of Theorem~\ref{th:A}.
This case can be also treated as Item $4$ of Theorem~\ref{th:A} for the parameters $qr=-1$.

Consider now the general case Items $3b$ and $4$ of Theorem~\ref{th:A}, i.e., $qrs\neq0$. If the
problem is solvable, then $K^{-1}$ can be written in the form (\ref{eq:A.19}) with
$\gamma_0,\gamma_1\in\mathbb{C}\setminus\{0\}$. Since $K^{-1}\in{\rm SL}(2,\mathbb{C})$, one
using Transformation~(\ref{eq:A.21}) proves the second condition in (\ref{eq:A.20}). Substitution
of Equation~(\ref{eq:A.19}) into System~(\ref{eq:A.2}) yields:
\begin{equation}
 \label{eq:A.22}
\gamma_1r(r_0-1/r_0)=f_0\gamma_0,\qquad q\gamma_0(r_1-1/r_1)=f_1\gamma_1.
\end{equation}
Multiplying Equations~(\ref{eq:A.21}) one obtains
\begin{equation}
 \label{eq:A.23}
f_0f_1=qr(r_0-1/r_0)(r_1-1/r_1).
\end{equation}
The second equations (\ref{eq:A.20}) and (\ref{eq:A.22}) yield the first equation in (\ref{eq:A.20}).

To prove that Equation~(\ref{eq:A.23}) is equivalent to Equations~(\ref{eq:A.18}) one should notice
that the matrix $(M_{06}-1/r_0)(M_{16}-r_1)$ has the following eigenvectors: $e_{11}$ and
$e_{1-1}-\frac{s}{q}e_{11}$ corresponding  to the eigenvalues $0$ and $qr(r_0-1/r_0)(r_1-1/r_1)$,
respectively.

Conversely, by substituting Equation~(\ref{eq:A.19}) into System~(\ref{eq:A.2}) and taking into account
Equations~(\ref{eq:A.20}) and (\ref{eq:A.21}) one finds equations (\ref{eq:A.18}).
\end{proof}
\begin{remark}
 \label{rem:A6}{\rm
Suppose the problem is unsolvable, i.e., Equations (\ref{eq:A.3}) and (\ref{eq:A.4}) hold. Define
the following matrices:
$$
\tilde M_{0\Delta}=\left(\begin{array}{cc}
1/r_0&\tilde f_0\\
0&r_0
\end{array}\right),\qquad
\tilde M_{1\Delta}=\left(\begin{array}{cc}
1/r_1&0\\
\tilde f_1&r_1
\end{array}\right)
$$
According to Propositions~(\ref{prop:A3}) and (\ref{prop:A4}) System~(\ref{eq:A.2}) for
the pair $M_{0\Delta}, \tilde M_{1\Delta}$ with $\tilde f_1=0$ and for the pair
$\tilde M_{0\Delta}, M_{1\Delta}$ with $\tilde f_0=0$.
}\end{remark}

\end{document}